\documentclass[natbib]{elsarticle}

\usepackage{subcaption}
\usepackage[utf8]{inputenc}
\usepackage[margin=1in]{geometry}
\usepackage{tikz-dimline}
\usepackage{morefloats}
\usepackage{setspace}
\usepackage{diagbox}
\usepackage{amssymb}
\usepackage{stmaryrd} 
\usepackage{algorithm,algorithmicx}
\usepackage{algpseudocode}
\usepackage{upgreek}
\usepackage{algpseudocode}
\usepackage{mathrsfs}
\usepackage{soul}
\usepackage{color}
\usepackage{url} 
\usepackage{graphicx} 
\usepackage{tikz}
\usetikzlibrary{positioning}
\usetikzlibrary{snakes}
\usetikzlibrary{arrows.meta}
\usetikzlibrary{intersections}
\usetikzlibrary{shapes}
\usetikzlibrary{calc,patterns,
                decorations.pathmorphing,
                decorations.markings}
\usetikzlibrary{arrows.meta}
\usepackage[EULERGREEK]{sansmath}
\usepackage{amsmath}

\usepackage{array}
\newcolumntype{P}[1]{>{\centering\arraybackslash}p{#1}}
\newcolumntype{M}[1]{>{\centering\arraybackslash}m{#1}}
\usepackage{tabularx}

\newlength\figurewidth
\newlength\figureheight

\usepackage{everypage}

\makeatletter
\newcommand{\subalign}[1]{%
    \vcenter{%
        \Let@ \restore@math@cr \default@tag
        \baselineskip\fontdimen10 \scriptfont\tw@
        \advance\baselineskip\fontdimen12 \scriptfont\tw@
        \lineskip\thr@@\fontdimen8 \scriptfont\thr@@
        \lineskiplimit\lineskip
        \ialign{\hfil$\m@th\scriptstyle##$&$\m@th\scriptstyle{}##$\crcr
            #1\crcr
        }%
    }
}
\newcommand{\norm}[1]{\left\lVert#1\right\rVert}
\makeatother

\begin{document}
\pagenumbering{arabic}
\newcommand{\rojo}[1]{#1}
\newcommand{\adapt}[1]{#1}
\renewcommand{\hl}[1]{#1}

\begin{frontmatter}
    \title{\rojo{Another source of mesh dependence in topology optimization.}}
    \author[LLNL]{Miguel A. Salazar de Troya\corref{correspondingauthor}}
    \author[LLNL]{Geoffrey M. Oxberry}
    \author[LLNL]{Cosmin G. Petra}
    \author[LLNL]{Daniel A Tortorelli}
    \address[LLNL]{Lawrence Livermore National Laboratory, Livemore, CA, USA, 94550}
    \cortext[correspondingauthor]{Corresponding author. }


    \begin{abstract}
        The topology optimization community has regularly employed nonlinear programming (NLP) algorithms  from the operations research community.
        However, these algorithms are implemented in the real vector space $\mathbb{R}^n$ instead of the proper function space where the design variable resides.
        \rojo{In this article, we show how the volume fraction variable discretization on non-uniform meshes affects the convergence of $\mathbb{R}^n$ based NLP algorithms.}
        We do so by first \rojo{summarizing} the functional analysis tools necessary to understand why convergence is affected by the mesh.
        Namely, the distinction between derivative and gradient definitions and the role of the mesh-dependent inner product within the NLP algorithm.
        These tools are then used to make the Globally Convergent Method of Moving Asymptotes (GCMMA), a popular NLP algorithm in the topology optimization community, converge in a mesh independent fashion when starting from the same initial design.
        We then benchmark our algorithm with different examples of topology optimization problems. Namely, ill-conditioned and large scale problems and designs problems with adaptive mesh refinement.
    \end{abstract}

\end{frontmatter}

\newcommand*{\paperarticle}{article}
\vspace{5mm}
\newcommand*{\MyPath}{.}%
\section{Introduction}
\label{sec:introduction}
Topology optimization finds the optimal distribution of material in a given design domain $D$ to minimize a cost function and satisfy constraint function inequalities.
In the classic element-wise uniform density based method, which in this paper we refer to as element-wise uniform volume fraction method\footnote{We opt to use volume fraction over density to avoid confusion with the physical quantity density, which is often used in topology optimization applications such as elastodynamics, fluids, etc.}, the optimization algorithm places material in individual elements of a background mesh to define the geometry of the optimal design.
Regions devoid of material are meaningless, hence the motivation to coarsen the mesh in these regions.
By the same token, regions which contain material require a higher mesh resolution.
This cost saving strategy that distributes elements with different sizes within the mesh is known as the Adaptive Mesh Refinement (AMR).

If the elements have different sizes, as in AMR, it is intuitively wrong to think that all design variables have the same contribution to the design.
In particular, we explain why it is necessary to accommodate the element size when calculating the inner products involving the design variables within the NLP algorithms.
However, this is not done within most NLP algorithms used in the topology optimization community as they assume the design is a vector
in the real vector space $\mathbb{R}^n$, i.e., simply a vector of length equal to the number of elements in the mesh, $n$, cf. IPOPT \citep{wachter2006}, SNOPT \citep{snoptmanual},  MMA \citep{mma}, FMINCON \citep{matlab} and Optimality Criteria \citep{ocmethod}.
On the other hand, the NLP libraries Optizelle \citep{optizelle}, Moola \citep{moola}, ROL \citep{ridzal} and TAO \citep{tao-user-ref} contain to various extents the capacity of treating design fields as elements of their underlying function spaces.
Related work by \cite{funke_book} compares mesh-independent and dependent versions of the steepest descent algorithms for an unconstrained problem and estimates their rates of convergence.
It is also necessary to mention that in the PDE-constrained optimization community, NLP algorithms are inherently mesh-independent as they are implemented in the corresponding function space.
For instance, \cite{ulbrich2009primal} implements an infinite-dimensional (inf-dim) primal-dual interior-point method with a Newton solver,
\cite{ziems2011adaptive} implements an inexact sequential quadratic programming method with an adaptive multilevel mesh refinement scheme and
\rojo{\cite{Blank2017} solves a phase field based topology optimization with a projected gradient method for cases where the cost function is only differentiable in $L^{\infty}$.}

We choose the $L^2$ function space to represent the set of possible designs on the design domain $D$.
It is equipped with an inner product and discretized to be piecewise uniform\footnote{We use uniform to describe functions that do not change in space and constant to describe functions that do not change in time.} over the finite elements.
Other infinite-dimensional function spaces choices, e.g $H^1$, are possible and should be addressed in the future.
The inconsistency of using a design field in $L^2$ with an NLP algorithm formulated in $\mathbb{R}^n$ is generally not a problem because most of the topology optimization studies use uniform meshes.
However, when using meshes with different element sizes, as is the case in AMR, the $\mathbb{R}^n$ viewpoint yields mesh-dependent designs, whereas the $L^2$ approach does not.
\hl{An immediate corollary is that restriction by filtration alone does not ensure mesh-independent designs, as is commonly accepted in the topology optimization community.}

This work is laid out as follows: Section \ref{sec:mathprelim} presents the mathematical tools we need to implement mesh-independent NLP algorithms.
We use these tools in Section \ref{sec:mmasec} to make one of the most popular NLP algorithms in the topology optimization community, the Globally Convergent Method of Moving Asymptotes (GCMMA/MMA) algorithm, mesh-independent.
In Section \ref{sec:examples}, \hl{we first validate the NLP algorithm by solving three common problems in topology optimization with contrived meshes specifically built to increase the ill-conditioning of the optimization problem.
We then apply the algorithm to a three dimensional problem with a uniform mesh and different levels of refinement to showcase the importance of our algorithm for large scale problems.
Finally, we solve two design problems with different physics and AMR applied during the optimization.}
Section \ref{sec:conclusions} briefly summarizes our findings and presents conclusions.

The function spaces concepts used in this article require a rigorous mathematical discussion to be absolutely precise.
Namely, it is necessary to show the differentiability of the convex approximation within the GCMMA to ensure that applying the Newton's method is mathematically sound.
These details would quickly obscure our main focus.
Our intent is to convey the differences between the $L^2$ and $\mathbb{R}^n$ NLP algorithms as simply as possible and to demonstrate their differences.
We therefore opt to take a more pragmatic approach in our discussions at the expense of glossing over important mathematical details.


\section{Mathematical Preliminaries}
\label{sec:mathprelim}
A topology optimization algorithm converges in a mesh-independent fashion by treating the design as a field, here a field in the $L^2$ space, using concepts from functional analysis.
The design field is then discretized in a consistent manner using the finite element basis, i.e. \rojo{resulting in} the \rojo{widely} used element-wise uniform volume fraction field\footnote{
    \rojo{\label{ft:elemwise}In this work, for simplicity, we focus on the most popular topology optimization approach, where the design is defined
        by an element-wise uniform material volume fraction in the Hilbert space $L^2$, but it could be extended to other parametrizations
        such as an $H^1$ nodal-based material volume fraction.}}.
Notably, the norms in the NLP algorithm that check for convergence are discretized in this finite element space.

To illustrate the proper discretization, consider the unconstrained minimization problem
\begin{equation}
    \begin{aligned}
        \begin{split}
            & \underset{\nu \in V}{\text{min}}
            & & \theta(\nu)\,, \\
        \end{split}
        \label{eq:orig_problem}
    \end{aligned}
\end{equation}
with the functional
\begin{align}
    \theta & :V\rightarrow \mathbb{R}\,,
\end{align}
where $\nu: V \rightarrow \mathbb{R}$ is our volume fraction design field that belongs to $V$, a Hilbert space on domain $D$, equipped with an inner product $(\cdot, \cdot)_V$, which induces the primal norm $\norm{\cdot}_V$.
For our topology optimization, $V=L^2(D)$, which is equipped with the norm
\begin{align}
    \norm{\nu}_{L^2} = \sqrt{(\nu, \nu)_{L^2}} =  \left( \int_{D} \nu^2 ~ dV  \right)^{1/2}\,.
\end{align}
The space of all bounded linear functionals that map $V$ to $\mathbb{R}$ is the dual space $V^*$, which is a subset of $\mathscr{L}(V, \mathbb{R})$, i.e. the space of linear operators from $V$ to $\mathbb{R}$, i.e. $V^* \subset \mathscr{L}(V, \mathbb{R})$.
Both the primal $\norm{\cdot}_V$ and dual $\norm{\cdot}_{V^*}$ norms can be used to check for convergence in NLP algorithms.



To formulate NLP algorithms on the function space $V$, we need the Riesz map from the Riesz representation theorem:
Let $V$ be a Hilbert space with inner product $(\cdot, \cdot)_V$ and dual space $V^*$.
For every $\varphi \in V^*$ there is a unique element $u\in V$ such that $\varphi (v) = \left( u, v \right)_V$ for all $v\in V$.
This one-to-one map is the Riesz map $\Phi:V \rightarrow V^*$ defined such that $\Phi(u) = \varphi$; it is an isometry between $V$ and $V^*$.

The discretization of the primal and dual spaces follow from \cite{funke_book}.
We approximate the volume fraction field $\nu \in V$ with $\nu_h \in V_h$, where $V_h$ is the span of basis functions $\mathscr{P} = \left\{\phi_1, ..., \phi_n \right\}$, $\phi_i \in \rojo{V_h}$ and $n$ is the dimension of $V_h$.
Our approximation now reads
\newcommand{\sumonn}[1]{\sum_{#1=1}^{n}}
\begin{align}
    \nu(\mathbf{x}) \approx \nu_h(\mathbf{x}) & = \sumonn{i} \nu^i \phi_i(\mathbf{x}) = \boldsymbol{\nu}^T \boldsymbol{\phi}(\mathbf{x})\,.
    \label{eq:l2_disc}
\end{align}
We similarly approximate $\iota\in V$ with $\iota_h \in V_h$ so that the inner product definition yields
\begin{equation}
    \begin{aligned}
        (\nu_h,\iota_h)_{V_h} & = \int_{D} (\boldsymbol{\nu}^T \boldsymbol{\phi}) (\boldsymbol{\iota}^T \boldsymbol{\phi}) dV \\
                              & = \boldsymbol{\nu}^T \int_{D} \boldsymbol{\phi} \boldsymbol{\phi}^T ~dV ~ \boldsymbol{\iota}  \\
                              & = \boldsymbol{\nu}^T \mathbf{M}  \boldsymbol{\iota} \,,\label{eq:disc_inner}
    \end{aligned}
\end{equation}
where
\begin{align}
    \mathbf{M} & = \int_{D} \boldsymbol{\phi} \boldsymbol{\phi}^T ~dV
    \label{eq:mass_matrix}
\end{align}
is the mass matrix that reflects the mesh discretization.
By construction, $\mathbf{M}$ is symmetric and invertible.

The discretized design field $\nu_h$ is in the Hilbert space $V_h=\left(\mathbb{R}^n, (\cdot,\cdot)_{\mathbf{M}}\right)$, i.e. it is a vector in $\mathbb{R}^n$ of dimension $n$ with an $\mathbf{M}$ inner product.
This inner product induces the norm $\norm{\nu_h}_{V_h} = \norm{\boldsymbol{\nu}}_{\mathbf{M}}= (\boldsymbol{\nu}^T \mathbf{M} \boldsymbol{\nu})^{1/2}$.

Clearly the $L^2$ norm $\norm{\nu_h}_{L^2} = (\boldsymbol{\nu}^T \mathbf{M} \boldsymbol \nu)^{1/2}$ differs from the $\mathbb{R}^n$ norm $\norm{\nu_h}_{\mathbb{R}^n} = (\boldsymbol{\nu}^T  \boldsymbol \nu)^{1/2}$.
In topology optimization, $\nu_h$ is usually discretized via piecewise uniform functions over the individual elements so $\mathbf{M} = diag\left( |\Omega_1|, ..., |\Omega_n|\right)$  where $|\Omega_e|$ is the volume of the element $\Omega_e$.
So if the mesh is uniform, $\mathbf{M} = |\Omega_e| \mathbf{I}$ and hence $\norm{\nu_h}_{L^2} = \sqrt{|\Omega_e|} \norm{\nu_h}_{\mathbb{R}^n}$.

The basis $\mathscr{P}=\left\{\phi_1, ..., \phi_n \right\}$ induces a unique dual basis $\mathscr{P}^* = \left\{\phi^{*1}, ..., \phi^{*n} \right\}$ for $V_h^*$ defined such that $\phi^{*i} \in \rojo{V_h^*}$ and $\phi^{*i}(\phi_j)= \delta_{j}^i ~\forall ~i,j=1,...,n$.
This $\mathscr{P}^*$ basis is used to discretize any $F\in V^*$ as $F_h \in V_h^*$ such that for all $\iota_h \in V_h$
\begin{align}
    F_h(\iota_h) = \sum^{n}_{i=1} F_i \phi^{*i} (\iota_h)\,,
\end{align}
\rojo{where $F_i=F(\phi_i)$ for $i=1,...,n$, i.e. the vector components $F_i$ are interpolated from $F\in V^*$.}
In this way, $F_h(\iota_h)$ is computed as
\begin{equation}
    \begin{aligned}
        F_h(\iota_h) & = \sumonn{i} F_i \phi^{*i} \left( \sumonn{j} \iota^j \phi_j \right)  \,, \\
                     & = \sumonn{i} F_i   \sumonn{j} \iota^j \phi^{*i}(\phi_j)   \,,            \\
                     & = \sumonn{i} F_i   \iota^i    \,,                                        \\
                     & =  \mathbf{F}^T   \boldsymbol{\iota}   \,,\label{eq:dual_norm_disc}
    \end{aligned}
\end{equation}
where we used the orthonormal property between the bases $\mathscr{P}$ and $\mathscr{P}^*$ and the linearity of $\phi^{*i}$.
From the Riesz representation theorem, there exists $\nu_h\in V_h$ such that $\Phi(\nu_h) = F_h$ or $\Phi^{-1}(F_h) = \nu_h$, where \rojo{$\Phi^{-1}: V_h^* \rightarrow V_h$}.
Therefore,
\begin{equation}
    \begin{aligned}
        F_h(\iota_h) & = (\overbrace{\Phi^{-1}(F_h)}^{\nu_h}, \iota_h)_{V} \,,                           \\
                     & = (\nu_h, \iota_h)_{V} \,,                                                        \\
                     & = \boldsymbol{\nu}^T \mathbf{M}  \boldsymbol{\iota} \,.\label{eq:dual_norm_disc2}
    \end{aligned}
\end{equation}
From Equations \eqref{eq:dual_norm_disc} and \eqref{eq:dual_norm_disc2} we can see that
\begin{align}
    \mathbf{F}^T   \boldsymbol{\iota}  =  \boldsymbol{\nu}^T \mathbf{M}  \boldsymbol{\iota}\,.
\end{align}
Therefore,
\begin{align}
    \mathbf{F} = \mathbf{M} \boldsymbol{\nu} \,,
\end{align}
and the discrete Riesz map and its inverse are defined such that
\begin{align}
    \Phi_h(\boldsymbol{\nu})=\mathbf{M}\boldsymbol{\nu}
\end{align}
and
\begin{align}
    \Phi_h^{-1}(\mathbf{F}) & = \mathbf{M}^{-1} \mathbf{F} \,.
    \label{eq:rieszmap}
\end{align}
Recalling that the Riesz map is an isometry between the spaces $V_h$ and $V_h^*$, we can now define and calculate the norm of an object $F_h \in V_h^*$ as
\begin{equation}
    \begin{aligned}
        \norm{F_h}_{V^*_h} & = \norm{\Phi^{-1} (F_h)}_{V_h} \,,                                               \\
                           & = \norm{\mathbf{M}^{-1} \mathbf{F}}_{\mathbf{M}} \,,                             \\
                           & = \sqrt{\mathbf{F}^T \mathbf{M}^{-1} \mathbf{M} \mathbf{M}^{-1} \mathbf{F} } \,, \\
                           & = \norm{\mathbf{F}}_{\mathbf{M}^{-1}} \,,
    \end{aligned}
\end{equation}
where we used the definition of the discrete inner product in Equation \eqref{eq:disc_inner}.


In this work, we use the Fr\'echet derivative $D\theta(\nu) \in V^*$ of the function $\theta : V \rightarrow \mathbb{R}$ at $\nu$.
If it exists, this derivative is defined such that\footnote{The ``little-$o$ notation'' $o(\norm{h}_V)$ for a functional $q: V \rightarrow \mathbb{R}$ means
    $\lim_{\norm{h}_V \to 0} \frac{q(h)}{ \norm{h}_V} =0$}
\begin{align}
    \theta(\nu + h) - \theta(\nu) - D\theta(\nu)[h] = o( \norm{h}_V)
    \label{eq:frechet_deriv}
\end{align}
for all $h\in V$.
By definition, $D\theta(\nu)\in V^*$, and hence the Riesz representation theorem tells us there is an object in $V$ that we will denote $\nabla \theta (\nu) \in V$, i.e. the gradient of $\theta$ at $\nu$ such that
$D\theta(\nu)[h]=(\nabla \theta(\nu), h)_V$ for all $h\in V$.
Using the Riesz map, $\Phi(\nabla \theta(\nu)) = D\theta(\nu)$ and because the Riesz map depends on the inner product $(\cdot, \cdot)_V$, so does $\nabla \theta(\nu) \in V$.
This inner product dependence is crucial in our NLP algorithm as the inner product depends on the mesh discretization, notably from \eqref{eq:rieszmap} we have
\begin{align}
    \boldsymbol{\nabla \theta} = \mathbf{M}^{-1} \boldsymbol{D \theta}
    \label{eq:rieszmap_grad}
\end{align}
where $\boldsymbol{\nabla \theta}$ and $\boldsymbol{D \theta}$ are the discrete counterparts of $\nabla \theta(\nu)$ and $D\theta(\nu)$.

We are now in position to show how these functional analysis concepts apply to NLP algorithms in the inf-dim space.
We start by examining the most basic NLP algorithm for the solution of the simple unconstrained minimization problem of \eqref{eq:orig_problem}.
i.e, the steepest descent algorithm, for which the iterate $\nu^{(k)}$ is updated as
\begin{align}
    \nu^{(k+1)} = \nu^{(k)} - \gamma \nabla \theta(\nu^{(k)}) \,.
    \label{eq:steepestdesc}
\end{align}
where $\gamma \geq 0$ is the step length.
The discretized Equation \eqref{eq:steepestdesc} becomes
\begin{equation}
    \begin{aligned}
        \boldsymbol{\nu}^{(k+1)} & = \boldsymbol{\nu}^{(k)} - \gamma \boldsymbol{\nabla \theta} \,.            \\
                                 & = \boldsymbol{\nu}^{(k)} - \gamma \mathbf{M}^{-1} \boldsymbol{D \theta} \,.
        \label{eq:steepestdescdisc}
    \end{aligned}
\end{equation}
When calculating the norm to check for convergence, we use $\norm{\nabla \theta (\nu)}_V$, which upon discretization is $\norm{\boldsymbol{\nabla \theta}}_{\mathbf{M}}$.

It seems intuitive that $\nu$ and $\nabla \theta(\nu)$ must be in the same function space since they are added together.
This motivates us to use the gradient $\nabla \theta(\nu)$ and not the derivative $D\theta(\nu)$ in \eqref{eq:steepestdesc}, which is contrary to most topology optimization algorithms.
For a uniform mesh $\mathbf{M}=|\Omega_e|\mathbf{I}$ so $\boldsymbol{\nabla \theta} = \frac{1}{|\Omega_e|} \mathbf{I} \boldsymbol{D\theta} $ \rojo{
    and hence $\boldsymbol{\nabla \theta}$
    and $\boldsymbol{D\theta}$ are
    parallel and there is no difference in the search direction.}
\rojo{However}, the number of iterations to convergence will be different due to the difference in the inner product.

The second NLP algorithm uses Newton's method, wherein we iterate to find $\nu$ such that
\begin{align}
    D\theta(\nu)[\delta \nu]= 0~ \forall \delta \nu \in V.
\end{align}
To do so, we linearize around $\nu^{(k)}$ and solve for the update $\Delta \nu^{(k)}$ via
\begin{align}
    D^2\theta(\nu^{(k)})[\Delta \nu^{(k)},\delta \nu] =-D\theta(\nu^{(k)})[\delta \nu]~~ \forall \delta \nu \in V \,,
    \label{eq:newtonstep}
\end{align}
where $D^2\theta(\nu^{(k)})[ \cdot, \cdot]: V \rightarrow  \mathscr{L} (V \times V, \mathbb{R})$ is the Hessian, i.e. second derivative of $\theta$ at $\nu^{(k)}$; it is a bilinear map from $(V \times V)$ to $\mathbb{R}$.
The difference here is that we need to supply the NLP algorithm with the derivative $D\theta(\nu)$ (and the Hessian $D^2\theta(\nu)$) and not the gradient $\nabla \theta(\nu)$ as in the steepest descent algorithm.
Upon discretization, Equation \eqref{eq:newtonstep} becomes
\begin{align}
    \boldsymbol{D}^2\boldsymbol{\theta} \boldsymbol{\Delta \nu}^{(k)}  =-\boldsymbol{D\theta}
\end{align}
When calculating the norm to check for convergence, we use $\norm{D \theta (\nu)}_{V^*}$, whose discretization is $\norm{\boldsymbol{D \theta}}_{\mathbf{M}^{-1}}$

{
Inspired by \cite{funke_book}, we showcase the difference between $\nabla \theta(\nu)$ and $D\theta(\nu)$ with the following one dimensional convex unconstrained optimization problem\label{page:exampleproblem}

\newcommand{\sinfunc}{\text{sin}\left(\frac{x}{4}\right)}
\begin{align}
    \underset{\nu\in V}{\text{min}} ~ \theta(\nu) = \frac{1}{2}( \nu, c \, \nu )_V - ( \nu, b )_V \,,
\end{align}
where $b$ and $c$ are given functions on $V=L^2(D)$ with $D=[1, 10]$, $c(x) = \sinfunc$ and $b(x) = x$.
The solution is trivially calculated by the stationary condition
\begin{equation}
    \begin{aligned}
        D\theta(\nu) [v] = 0 & =  ( \nabla \theta, v )_V  \,,  \\
                             & =      ( c \, \nu - b,v )_V \,.
    \end{aligned}
    \label{eq:derivexample}
\end{equation}
which must hold for all $\nu \in V$ and hence $\nu(x)=\frac{b(x)}{c(x)}$.
We proceed to discretize the function $\theta$ with piecewise uniform elements, resulting in the expression
\newcommand{\nubold}{\boldsymbol{\nu}}
\begin{align}
    \theta(\nubold) =  \frac{1}{2} \nubold^T \mathbf{H} \nubold - \nubold^T \mathbf{M b}  \,,
\end{align}
whose derivative is
\newcommand{\gradient}{\mathbf{H}\nubold^{(k)} -\mathbf{M}\mathbf{b}}
\begin{align}
    \boldsymbol{D \theta}=  \gradient \,,
    \label{eq:discderiv}
\end{align}
where the Hessian matrix $\mathbf{H}$ is calculated as
\begin{align}
    \mathbf{H} = \int_D c(x) \boldsymbol \phi(x)  \boldsymbol \phi(x)^T ~dV \,,
\end{align}
\rojo{using one point quadrature per element.
    The vectors $\boldsymbol \nu$ and $\bf b$ represent the values of the functions $\nu(x)$ and $b(x)$ at the quadrature points.}
Note that Equation \eqref{eq:discderiv} is the discretization of $D\theta(\nu)$ in Equation \eqref{eq:derivexample}.
Applying the steepest descent algorithm in $\mathbb{R}^n$ yields
\begin{equation}
    \begin{aligned}
        \nubold^{(k+1)} & = \nubold^{(k)} - \gamma \boldsymbol{D \theta} \,, \\
                        & = \nubold^{(k)} - \gamma(\gradient)\,.
        \label{eq:steepestexample}
    \end{aligned}
\end{equation}
The optimal step size $\gamma$ is calculated with the closed-form expression
\newcommand{\bbold}{\mathbf{b}}
\newcommand{\stepsize}{\frac{(\gradient)^T(\gradient)}{(\gradient)^T \mathbf{H} (\gradient)}}
\begin{align}
    \gamma = \stepsize \,,
\end{align}
which we use in Equation \eqref{eq:steepestexample} to obtain the fixed point iteration
\begin{align}
    \nubold^{(k+1)} = \nubold^{(k)} - \stepsize (\gradient) \,.
    \label{eq:fixedsteepest}
\end{align}

\newcommand{\newgradient}{\mathbf{M}^{-1}\mathbf{H}\nubold^{(k)} - \bbold}
On the other hand, in the $L^2$ reformulation,
we replace the descent direction $\boldsymbol{D \theta}$ with the gradient $\boldsymbol{\nabla \theta} = \mathbf{M}^{-1} \boldsymbol{D\theta}(\nubold) =  \newgradient$ in Equation \eqref{eq:steepestexample}.a and calculate the optimal step size
\newcommand{\newstepsize}{\frac{(\mathbf{M}^{-1} \mathbf{H}\boldsymbol\nu^k - \mathbf{b})^T\mathbf{M} (\newgradient)}{(\newgradient)^T \mathbf{M} (\newgradient)}}
\begin{align}
    \gamma = \newstepsize = 1
\end{align}
to obtain the fixed point iteration
\begin{align}
    \nubold^{(k+1)} = \nubold^{(k)} - (\newgradient)\,.
    \label{eq:newfixedsteepest}
\end{align}

We run both fixed point iteration Equations \eqref{eq:newfixedsteepest} and \eqref{eq:fixedsteepest} starting from $\boldsymbol \nu^{(0)}=\boldsymbol 0$.
Convergence is declared when the error $e \leq 10^{-7}$, where $e=\norm{\boldsymbol \nu - \mathbf{b}\oslash \mathbf{c}}$ for $\mathbb{R}^n$ (the operator $\oslash$ is the \rojo{Hadamard} division) and $e=\norm{\nu - b / c}_{L^2}$ for $L^2$.

Table \ref{tab:steeptable} denotes the iteration history of both methods with $n$ design variables over a one-dimensional mesh with nodes at positions $x_r=10^{y_r}$ where $y_r=\{\frac{r}{n+1}\}_{r=0}^{n+1}$.
The convergence in the $\mathbb{R}^n$ NLP algorithm deteriorates with the number of elements, whereas the number of iterations in the $L^2$ NLP algorithm remains nearly constant.
The slower convergence of the $\mathbb{R}^n$ NLP algorithm is attributed to the fact that we are adding members from different spaces ($\boldsymbol\nu \in V_h$ and $\boldsymbol{D\theta}\in V_h^*$).

\begin{table}[h!]
    \begin{center}
        \begin{tabular}{c|c|c|c} 
            \diagbox{Method}{$n$}
                           & $10^1$ & $10^3$ & $10^5$ \\
            \hline
            $\mathbb{R}^n$ & 195    & 270    & 302    \\
            $L^2$          & 52     & 56     & 56     \\
        \end{tabular}
        \caption{Iteration count for the discrete steepest descent of algorithms Equations \eqref{eq:fixedsteepest} and \eqref{eq:newfixedsteepest}.}
        \label{tab:steeptable}
    \end{center}
\end{table}
}

\newcommand{\sumonm}{\sum_{i=1}^m}
\section{GCMMA in function space}
\label{sec:mmasec}
We are now motivated to formulate the first-order GCMMA algorithm \citep{gcmma} in the $L^2$ space.
GCMMA \citep{gcmma} and its non globally-convergent version MMA \citep{mma} are widely used NLP algorithms in the topology optimization community.
Following the implementation given in \cite{svanbergmma}, we highlight here the necessary changes to make the GCMMA algorithm converge in a mesh-independent fashion.
To begin, we consider the optimization problem
\begin{equation}
    \begin{aligned}
        \label{eq:minproblem}
        \begin{split}
            & \underset{\nu \in V}{\text{min}}
            & & \theta_0(\nu) \,,\\
            & \text{s.t.}
            & & \theta_i(\nu) \leq 0, \; i = 1, \ldots, m\,, \\
            &&& \nu_{\text{min}} \leq \nu \leq \nu_{\text{max}} \; \text{a.e}\,.
        \end{split}
    \end{aligned}
\end{equation}
First, the artificial optimization variables $\mathbf{y} = (y_1,...,y_m)$ \hl{are added to ensure feasibility} and $z$ \hl{is added make certain subclasses of problems, like least squares or minmax problems, easier to formulate, i.e.}
\begin{equation}
    \begin{aligned}
                       & \underset{
            \subalign{
        \nu            & \in V                                                                                                                                                                               \\
        \boldsymbol{y} & \in \mathbb{R}^m                                                                                                                                                                    \\
        z              & \in \mathbb{R}
            }
        }
        {\text{min}
        }
                       &                  & \theta_0(\nu) + a_0z + \sum\limits_{i=1}^m \left( c_i y_i + \frac{1}{2}d_i y_i^2 \right) \,,                                                                     \\
                       & \text{s.t.}
                       &                  & \theta_i(\nu) - a_iz - y_i\leq 0, \; i = 1, \ldots, m\,,                                                                                                         \\
                       &                  &                                                                                              & \nu_{\text{min}} \leq \nu  \leq \nu_{\text{max}}\; \text{a.e.}\,, \\
                       &                  &                                                                                              & \boldsymbol{y} \geq 0 \,,                                         \\
                       &                  &                                                                                              & z\geq 0 \,.
    \end{aligned}
    \label{eq:mmaproblem}
\end{equation}
where $a_0, a_i, c_i$ and $d_i$ are real numbers which satisfy $a_0 > 0, a_i \geq 0, c_i \geq 0, d_i \geq 0$ and $c_i + d_i > 0$ for all $i$, and also $a_i c_i > a_0$ for all $i$ \citep{gcmma}.
Note that we recover the original NLP algorithm for $z=0$ and $\boldsymbol{y}=0$.

For each optimization iteration $k$, we solve the following convex approximate subproblem based on Equation \eqref{eq:mmaproblem}, the cost and constraint functions and
their derivatives and the values at the current iterate $(\nu^{(k)},\boldsymbol{y}^{(k)},z^{(k)})$.
Ultimately, we iterate by solving
\begin{equation}
    \begin{aligned}
        (\nu^{(k+1)}, \boldsymbol{y}^{(k+1)}, & z^{(k+1)}) =
        \underset{
            \subalign{
        \nu                                   & \in V                                                                                                                                                                 \\
        \boldsymbol{y}                        & \in \mathbb{R}^m                                                                                                                                                      \\
        z                                     & \in \mathbb{R}
            }
        }{\text{arg min}}
                                              &                  & \tilde{\theta}_0(\nu) + a_0z + \sum\limits_{i=1}^m \left( c_i y_i + \frac{1}{2}d_i y_i^2 \right) \,,                                               \\
                                              & \text{s.t.}
                                              &                  & \tilde{\theta}_i(\nu) - a_iz - y_i\leq 0, \; i = 1, \ldots, m. \,,                                                                                 \\
                                              &                  &                                                                                                      & \alpha \leq \nu \leq \beta\; \text{a.e} \,, \\
                                              &                  &                                                                                                      & \boldsymbol{y} \geq 0 \,,                   \\
                                              &                  &                                                                                                      & z\geq 0 \,.
    \end{aligned}
    \label{eq:mmasubproblem}
\end{equation}
where the newly introduced functions $\tilde{\theta}_0$ and  $\tilde{\theta}_i$ and bounds $\alpha$ and $\beta$ are defined momentarily.

In our formulation of the above subproblem, we replace the summations of the approximating functionals $\tilde{\theta}_i$ in \cite{svanbergmma} with integrals over the domain.
\begin{align}
    \tilde{\theta}_i(\nu) & =\int_D \left(\frac{p_{i}}{U^{(k)} - \nu}+\frac{q_{i}}{\nu - L^{(k)}}\right)~dV + r_i , \; i = 0,1, \ldots, m \,,      \label{eq:convex_approx} \\
    r_i                   & = \theta_i(\nu^{(k)}) - \int_D \left(\frac{p_{i}}{U^{(k)} - \nu^{(k)}}+\frac{q_{i}}{\nu^{(k)} - L^{(k)}}\right)~dV  \,,
    \label{eq:mma_integrals}
\end{align}
where
\begin{align}
    p_{i} & = (U^{(k)} - \nu^{(k)})^2\left(1.001\left(\nabla \theta_i(\nu^{(k)})\right)^+
    +0.001\left(\nabla \theta_i(\nu^{(k)})\right)^- +\frac{\rho_i^{(k,j)}}{\nu_{\text{max}} - \nu_{\text{min}}}\right) \label{eq:pfunc} \,, \\
    q_{i} & = (\nu^{(k)} - L^{(k)})^2\left(0.001\left(\nabla \theta_i(\nu^{(k)})\right)^+
    +1.001\left(\nabla \theta_i(\nu^{(k)})\right)^- +\frac{\rho_i^{(k,j)}}{\nu_{\text{max}} - \nu_{\text{min}}}\right) \label{eq:qfunc} \,.
\end{align}
and $U^{(k)}$ and $L^{(k)}$ are the soon to be defined moving upper and lower asymptotes; they are all elements of $V$.
We emphasize here that the original GCMMA implementation does not make a distinction between gradients and derivatives when building the convex approximation in Equation \eqref{eq:convex_approx}.
As shown in this paper, it is vital to use gradients $\nabla \theta_i(\nu^{(k)})$ for $i=0,1, \ldots,m$.
It is important to warn readers that the convex approximations in Equation \eqref{eq:convex_approx} are not Frech\'et differentiable where $\nu=0$ in $L^2$.
We do not allow this $\nu=0$ situation, but a more mathematically rigorous rederivation of the method to accept such cases should be considered in the future.
To ensure the subproblem is convex, we use the ramp like functions $\left(a\right)^+ = \text{max} (0, a)$, and $\left(a\right)^-= \text{max} (0 ,-a)$.
The bounds (now fields in $V$) $\alpha$ and $\beta$ are taken as
\begin{equation}
    \begin{aligned}
        \alpha & = \text{max} \{\nu_{\text{min}}, L^{(k)} + 0.1 (\nu^{(k)} - L^{(k)}), \nu^{(k)} - 0.5 (\nu_{\text{max}} -\nu_{\text{min}} )\} \\
        \beta  & = \text{min} \{\nu_{\text{max}}, U^{(k)} - 0.1 (U^{(k)} - \nu^{(k)}), \nu^{(k)} + 0.5 (\nu_{\text{max}} -\nu_{\text{min}} )\}
    \end{aligned}
    \label{eq:alphabeta}
\end{equation}

\begin{figure*}
    \centering
    \tikzstyle{block} = [rectangle, draw, fill=blue!20,
    text width=3cm, text centered, rounded corners, minimum height=0.5cm]
    \tikzstyle{input} = [draw, ellipse, fill=green!20,
    text width=3cm, text centered, minimum height=0.8cm]
    \tikzstyle{output} = [draw, ellipse, fill=red!20,
    text width=3cm, text centered, minimum height=0.8cm]
    \tikzstyle{decision} = [diamond, aspect= 3, draw, fill=orange!20, inner sep=-6pt,
    text width=5cm, text centered, minimum height=0.5cm]
    \tikzstyle{decision2} = [diamond, aspect= 3, draw, fill=orange!20, inner sep=-6pt,
    text width=5cm, text centered, minimum height=0.4cm]
    \tikzstyle{line} = [draw, -latex]
    \begin{tikzpicture}[node distance = 1.2cm, auto]
        \sffamily
        \node [input] (nu_k) {$\nu^{(0)}$};
        \node [block, below of = nu_k] (obj_grad) {Calculate $\theta_i(\nu^{(k)}), ~\nabla \theta_i(\nu^{(k)})$};
        \node [block, below of = obj_grad] (param_gc) {Obtain $\rho^{(k,0)}$};
        \node [block, below of = param_gc, node distance = 1.2cm] (fapp) {Build $\tilde \theta_i^{(k,j)}(\nu^{(k)})$};
        \node [block, text width=5cm, below of = fapp] (nu_new) {Solve for $\nu^{(k,j)}$ in subproblem Equation \eqref{eq:mmasubproblem}};
        \node [block, text width=5cm, below of = nu_new] (new_f) {Calculate $\theta_i (\nu^{(k,j)})
                ,\tilde{\theta}_i(\nu^{(k,j)})$};
        \node [decision, below of = new_f, node distance = 2.5cm] (f_condition) {
            \begin{align*}
                \text{If} ~ \tilde \theta_i(\nu^{(k,j)})
                                  & \geq \theta_i(\nu^{(k,j)}) \\
                \forall i=0,...,m &
            \end{align*} };
        \node [output, below of = f_condition, node distance = 3cm] (new_opti) {$\nu^{(k+1, 0)} = \nu^{(k,j)}$};
        \node [block, right of = f_condition, node distance = 6cm] (repeat) {Obtain $\rho^{(k,j+1)}$};
        \path [line] (repeat) |- node[near start]{$j=j+1$} (fapp);
        \path [line] (f_condition) -- node{No} (repeat);
        \path [line] (f_condition) -- node{Yes} (new_opti);
        \path [line] (nu_k) -- (obj_grad);
        \path [line] (obj_grad) -- (param_gc);
        \path [line] (param_gc) -- node {$j=0$} (fapp);
        \path [line] (fapp) -- (nu_new);
        \path [line] (nu_new) -- (new_f);
        \path [line] (new_f) -- (f_condition);
        \node [decision2, below of = new_opti, node distance = 1.5cm] (stopping) {Termination criteria};
        \path [line] (stopping) -- node{No} ($(stopping) + (-5,0)$) |- node[near start]{$k=k+1$}  (obj_grad);
        \path [line] (new_opti) -- (stopping);
        \node [block, below of = stopping, node distance = 1.5cm] (final) {Solution};
        \path [line] (stopping) -- node{Yes} (final);
    \end{tikzpicture}
    \caption{GCMMA algorithm.}
    \label{fig:gcmmainner}
\end{figure*}
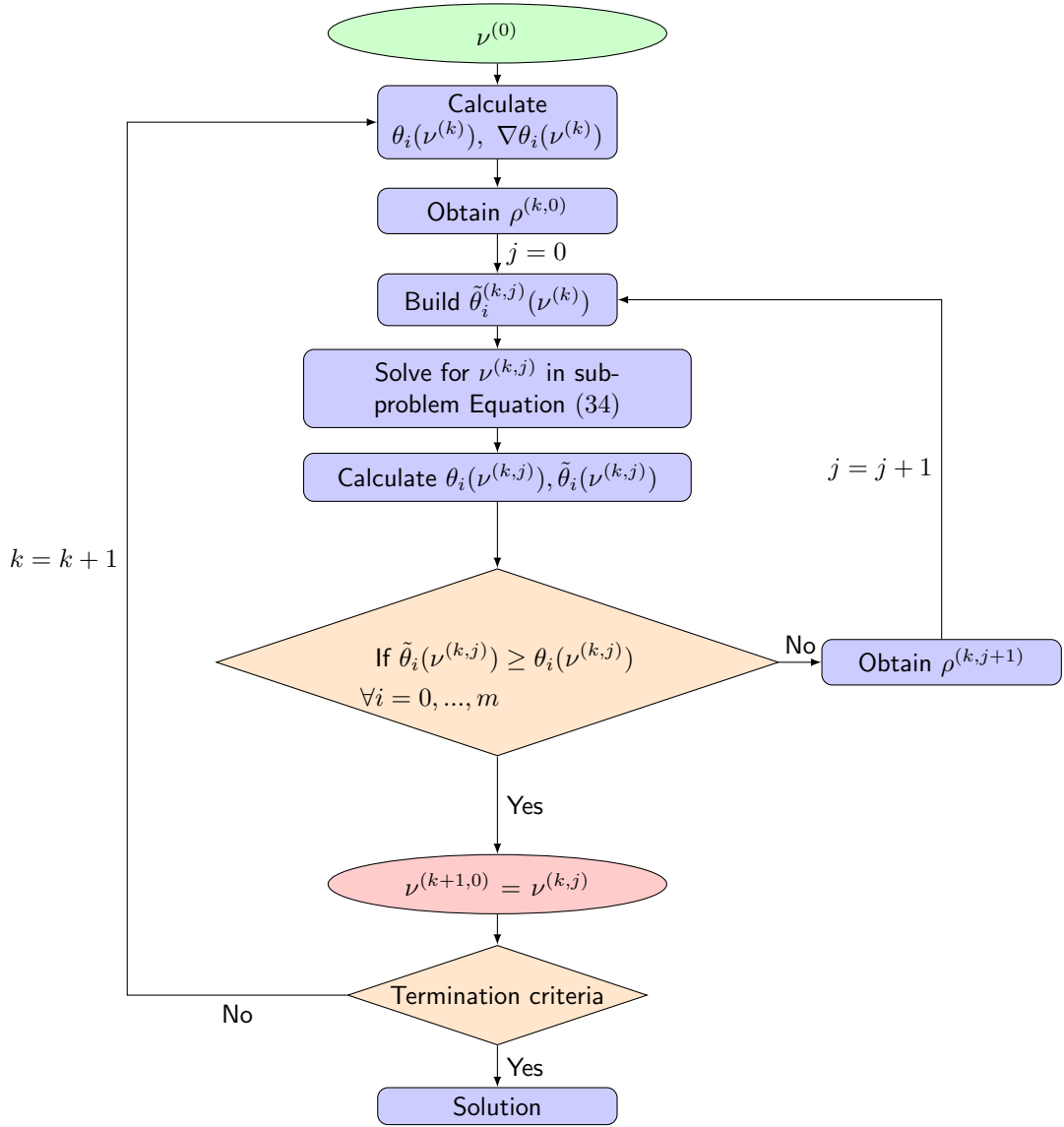

The GCMMA differs from the MMA in its attempt to achieve global convergence by controlling the parameter $\rho_i^{(k,j)}$ (which in the MMA is a fixed small positive value, usually lower than $10^{-5}$) in Equations \eqref{eq:pfunc} and \eqref{eq:qfunc}.
Here, the added superscript $j$ corresponds to the inner iteration within the GCMMA.
For the initial $j=0$ inner iteration, the solution $(\nu^{(k,0)},\boldsymbol{y}^{(k,0)},z^{(k,0)})$ of the Equation \eqref{eq:mmasubproblem} subproblem, whose details are explained later, is accepted if
\begin{align}
    \tilde \theta_i(\nu^{(k,j)}) \geq \theta_i(\nu^{(k,j)})  ;~ ~ i=0,...,m  \,,
    \label{eq:condition_check}
\end{align}
whereupon the outer iteration $k+1$ commences with the initial iterate $(\nu^{(k+1)},\boldsymbol{y}^{(k+1)},z^{(k+1)}) = (\nu^{(k,0)},\boldsymbol{y}^{(k,0)},z^{(k,0)})$.
Otherwise, the $j+1$ subproblem \eqref{eq:mmasubproblem} is solved with a more conservative convex approximation by replacing $\rho_i^{(k,j)}$ with $\rho_i^{(k,j+1)} > \rho_i^{(k,j)}$ and the Equation \eqref{eq:condition_check} inequality is rexamined.
If Equation \eqref{eq:condition_check} is satisfied, we begin the outer iteration $k+1$ with the initial iterate $(\nu^{(k+1,0)},\boldsymbol{y}^{(k+1,0)},z^{(k+1,0)}) = (\nu^{(k, j+1)}, \boldsymbol{y}^{(k, j+1)}, z^{(k, j+1)})$, otherwise, the inner $j+2$ subproblem is solved and so on, cf. Figure \ref{fig:gcmmainner}.
The termination criteria will be explained in detail later.

The parameters $\rho_i^{(k,j)}$ are calculated following \cite{gcmma}, but with integrals over $D$ replacing summations.
For subproblem $(k,0)$
\begin{align}
    \rho^{(k,0)}_i = \frac{0.1}{\hat V} \int_D \lvert \nabla \theta_i (\nu^{(k,0)}) \rvert \left( \nu_{\text{max}}  - \nu_{\text{min}}  \right) ~dV ~\text{for}\; i = 0,1, \ldots, m \,,
    \label{eq:rho_mechanism}
\end{align}
where $\hat V$ is the volume (area) of $D$.
For the subsequent $(k,j+1)$ subproblems
\begin{equation}
    \begin{aligned}
        \rho^{(k, j+1)}_i & = \text{min} \left\{ 1.1 \left(\rho^{(k, j)}_i + \delta_i^{(k,j)} \right) , 10 \rho_i^{(k,j)} \right\} & \text{if} ~\delta_i^{(k,j)} > 0 \,,   \\
        \rho^{(k, j+1)}_i & = \rho^{(k, j)}_i                                                                                      & \text{if} ~\delta_i^{(k,j)} \leq 0\,,
    \end{aligned}
\end{equation}
where
\begin{align}
    \delta^{(k,j)}_i = \frac{\theta_i\left(\nu^{(k,j)}\right) - \tilde{\theta}_i\left(\nu^{(k,j)}\right)}{d\left(\nu^{(k,j)}\right)} \,.
\end{align}
with
\begin{align}
    d (\nu) = \int_D \frac{\left( U^{(k)} - L^{(k)} \right) \left(\nu - \nu^{(k)} \right)^2}{\left(U^{(k)} - \nu\right) \left( \nu - L^{(k)} \right) \left( \nu_{\text{max}}  - \nu_{\text{min}} \right)} ~dV \,,
    \label{eq:d_mechanism}
\end{align}

The moving asymptote fields $L \in V$ and $U \in V$ are updated via heuristic rules.
For iterations $k=1$ and $k=2$
\begin{equation}
    \begin{aligned}
        L^{(k)} = \nu^{(k)} - 0.5(\nu_{\text{max}} - \nu_{\text{min}}) \,, \\
        U^{(k)} = \nu^{(k)} + 0.5(\nu_{\text{max}} - \nu_{\text{min}}) \,.
    \end{aligned}
    \label{eq:LUupdate1}
\end{equation}
For iterations $k \geq 3 $,
\begin{equation}
    \begin{aligned}
        L^{(k)} = \nu^{(k)} - \gamma^{(k)}(\nu^{(k-1)} - L^{(k-1)}) \,, \\
        U^{(k)} = \nu^{(k)} + \gamma^{(k)}(U^{(k-1)} - \nu^{(k-1)}) \,.
    \end{aligned}
    \label{eq:asympt}
\end{equation}
The field $\gamma$ that appears in Equations \eqref{eq:asympt} is determined by the values of $\nu$ in the last three outer iterations.
When there is no oscillation in $\nu$, we reduce the convexity by pushing the asymptotes further apart by choosing a larger $\gamma$ to accelerate convergence.
Otherwise, we use a smaller value to move the asymptotes closer together.
Specifically, we assign
\begin{equation}
    \begin{aligned}
        \gamma^{(k)} =
        \begin{cases}
            0.7 & \text{if}~(\nu^{(k)} - \nu^{(k-1)})(\nu^{(k-1)} - \nu^{(k-2)}) < 0 \,,  \\
            1.2 & \text{if} ~(\nu^{(k)} - \nu^{(k-1)})(\nu^{(k-1)} - \nu^{(k-2)}) > 0 \,, \\
            1   & \text{if} ~(\nu^{(k)} - \nu^{(k-1)})(\nu^{(k-1)} - \nu^{(k-2)}) = 0\,,
        \end{cases}
    \end{aligned}
    \label{eq:gammaupdate}
\end{equation}
subject to the inequalities
\begin{equation}
    \begin{aligned}
        L^{(k)} & \leq  \nu^{(k)} - 0.01(\nu_{\text{max}} - \nu_{\text{min}}) \,, \\
        L^{(k)} & \geq  \nu^{(k)} - 10(\nu_{\text{max}} - \nu_{\text{min}}) \,,   \\
        U^{(k)} & \geq  \nu^{(k)} + 0.01(\nu_{\text{max}} - \nu_{\text{min}}) \,, \\
        U^{(k)} & \leq  \nu^{(k)} + 10(\nu_{\text{max}} - \nu_{\text{min}}) \,.
    \end{aligned}
    \label{eq:LUrestrict}
\end{equation}

\rojo{We remark that the pointwise Equations \eqref{eq:alphabeta}, \eqref{eq:LUupdate1} - \eqref{eq:LUrestrict} are discretized directly by using their element-wise counter parts, which is consistent with our $L^2$ element-wise piecewise uniform parameterization of $\nu$.
    We do not worry here about the existence of these discretized counterparts as it is outside of the scope of this paper.}

\rojo{From here, one can solve the MMA subproblem following similar steps to those in the original article \citep{svanbergmma} with the exception of the calculation
    of the norms in their corresponding function spaces.
    We provide these details in \ref{app:appendixA}, where we also summarize all the necessary changes to the original GCMMA.
}

\section{Numerical examples}
\label{sec:examples}
To illustrate the effectiveness of incorporating the $L^2$ function space approach in the GCMMA, we solve six topology optimization problems, five in two (2D) and one in three (3D) dimensions.
For the 2D cases, we use triangular elements, and hexahedral elements for 3D,
both with first order Lagrange basis functions to represent the displacement.
The volume fraction field discretization uses the typical topology optimization 
approach with element-wise uniform basis functions.
All of the examples were solved using the finite element library Firedrake \citep{Rathgeber2016} \citep{Luporini2016} \citep{Homolya2017},
which uses PETSc \citep{petsc-efficient} \citep{ petsc-user-ref} \citep{Dalcin2011} as the backend for the linear algebra.
We use the direct solver MUMPS \citep{MUMPS01, MUMPS02} in 2D and the PETSc GAMG preconditioner in 3D.
The 2D optimizations ran on a single 2.60 GHz Intel XeonE5-2670 processor.
We employed up to 36 processors for the 3D cases.
The modified GCMMA is an adaptation of a Python implementation of the original MMA algorithm from the GetDP finite element library \citep{dular1998general}.
It was rewritten for better performance in parallel and to include an interface for the Firedrake-adjoint library \citep{Mitusch2019}.
We use the MMA parameters $a_0 = 1$, $c_i = 10000$ and $a_i = d_i = 0$ for all $i \geq 1$.
All results are visualized with ParaView \citep{paraview} and the graphs are plotted with Matplotlib \citep{Hunter:2007}.
To launch all the simulations, we use Signac \citep{signac_commat} and Signac-flow \citep{signac_scipy_2018}.

\subsection{Ill-conditioned meshes}
\hl{In this subsection, we deliberately use meshes with highly refined regions which at first glance,
can be deemed as cherry-picked to validate our approach.
    However, these meshes render ill-conditioned optimization problems and it is precisely this issue
    that we wish to highlight and resolve.}
We first solve three common topology optimization problems in linear elasticity.
The topology optimization problem is formulated as
\begin{equation}
    \begin{aligned}
        \underset{\nu \in \rojo{V}}{\text{min}}~ \theta_0(\nu)                                      & = \int_{D} \pi(\hat{\nu}, \mathbf{u})~ dV \,,                                 \\
        \text{such that} ~ \mathbf{u} \in \rojo{W} \text{satisfies} ~ a(\nu; \mathbf{u},\mathbf{v}) & = L(\mathbf{v}) ~ \text{for all}~ \mathbf{v} \in \rojo{W}  \,,                \\
        \theta_i(\nu)                                                                               & =  \int_{D} g_i(\hat{\nu}, \mathbf{u}) dV \leq 0               & i=1,2...m\,,
    \end{aligned}
    \label{eq:elastic_optimization}
\end{equation}
where
\begin{equation}
    a(\nu; \mathbf{u},\mathbf{v}) =\int_{D} r(\hat\nu) \mathbb{C}[\nabla \mathbf{u}] \cdot \nabla \mathbf{v}  ~dV  \,,
    \label{eq:weak_form_simp_mma}
\end{equation}
and
\begin{equation}
    L(\mathbf{v}) =\int_{\Gamma_N} \mathbf{t} \cdot \mathbf{v} ~da  \,.
\end{equation}

\rojo{
    The function spaces used are
    \begin{equation}
        V = \{\nu \in L^2(D) ~|~ 0 \leq \nu \leq 1\}
    \end{equation}
    and
    \begin{equation}
        W = \{\mathbf{u} \in [H^1(D)]^3 ~|~ \mathbf{u}|_{\Gamma_D} = 0\}
    \end{equation}
}

\rojo{The domain boundary $\Gamma$ is comprised of three complementary regions: $\Gamma_D$, $\Gamma_N$ and
    $\Gamma_F$ over which the Dirichlet, non-homogeneous Neumann and the homogeneous Neumann boundary conditions are applied.
    The functionals $\theta_i:L^2 \to \mathbb{R}$, $i=0,1,...$ are assumed to be Fr\'echet differentiable. \label{page:costfunction}}
The filtered volume fraction $\hat\nu$ in the above is obtained from the PDE-based filter \citep{boyanfilter} to generate a well-posed topology optimization problem
\rojo{
    \begin{equation}
        \begin{aligned}
            -\kappa \nabla^2 \hat{\nu} + \hat\nu     & = \nu ~ & \text{in} ~  D      & \,, \\
            \kappa \nabla \hat{\nu} \cdot \mathbf{n} & = 0~    & \text{on} ~  \Gamma &
        \end{aligned}
        \label{eq:pdefiltermma}
    \end{equation}
}
where $\kappa$ determines the minimum length scale of the design such that a small (large) 
$\kappa$ allows for fine (coarse) scale design fluctuations.
We solve Equation \eqref{eq:pdefiltermma} with a finite volume scheme to maintain the 
volume fraction values between 0 and 1.
We use the SIMP penalization \citep{bendsoe2013topology} to encourage 0-1 designs,
i.e. designs where 0 or $\nu=1$ almost everywhere.
As such,
\begin{align}
    r(\hat\nu) = \epsilon_{\nu} + (1 - \epsilon_{\nu})\hat{\nu}^3 \,,
    \label{eq:simp}
\end{align}
where $\epsilon_{\nu}=10^{-5}$ ensures the stiffness matrix in the finite element analysis is nonsingular.

Finally, $\mathbb{C}$ is the elasticity tensor corresponding to an isotropic material with Young modulus $E=1$ and Poisson ratio $\upnu=0.3$ and $\mathbf{t}$ is the applied traction on the surface $\Gamma_N$.

As usual, a reduced space approach is taken wherein we account for dependence of $\mathbf{u}$ on $\nu$, i.e. $\mathbf{u} \rightarrow \mathbf{u}(\nu)$ and the adjoint method is used to calculate the derivatives of the cost and constraint functions $\theta_i$.

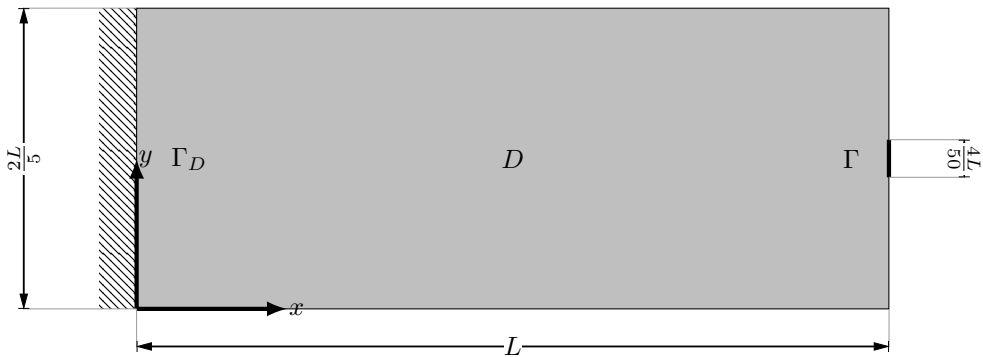
\begin{figure}[!h]
    \center
    \tikzset{>=latex}
    \begin{tikzpicture}[  spring/.style = {decorate,decoration={zigzag,amplitude=6pt,segment length=4pt}}
            , scale=0.1, every node/.style={inner sep=0pt}]
        \node (SW) at (0,0) {};
        \node (SE) at (100,0) {};
        \node (NE) at (100,40) {};
        \node (NW) at (0,40) {};
        \draw[fill=gray!50!white] (0,0) -- (100,0) -- (100,40) -- (0,40) -- (0,0);
        \node (D) at (50, 20) {$D$};
        \fill[pattern=north west lines] (SW) rectangle ($(NW) - (5.0, 0.0)$);
        \node (GammaD) at ($(NW) !.5! (SW) + (7.0, 0)$) {$\Gamma_D$};
        \node (loadP) at ($0.5*(NE) + 0.5*(SE)$) {};
        \draw[ultra thick] ($(loadP) + (0.0, 2.5)$) -- ($(loadP) + (0.0, -2.5)$);
        \node (Gamma) at ($(loadP) - (5, 0)$) {$\Gamma$};
        \draw[->, ultra thick] (0,0) -- (0, 20) node[right]{$y$};
        \draw[->, ultra thick] (0,0) -- (20, 0) node[right]{$x$};
        \node (A) at ($(SW) - (0.0, 5.0)$) {};
        \node (B) at ($(SE) - (0.0, 5.0)$) {};
        \dimline[line style = {line width=0.7},
            extension start length=-5cm,
            extension end length=-5cm]{(A)}{(B)}{$L$};
        \node (C) at ($(NW) - (15.0, 0.0)$) {};
        \node (D) at ($(SW) - (15.0, 0.0)$) {};
        \dimline[line style = {line width=0.7},
            extension start length=15cm,
            extension end length=15cm]{(D)}{(C)}{$\frac{2L}{5}$};
        \dimline[line style = {line width=0.7},
            extension start length=10cm,
            extension end length=10cm]{($(loadP) + (10, 2.5)$)}{($(loadP) + (10, -2.5)$)}{$\frac{4L}{50}$};
    \end{tikzpicture}
    \caption{Compliance domain. $L=100$}
    \label{fig:compliancedomain}
\end{figure}

The first problem we study is the proverbial structural compliance minimization subject to a \rojo{maximum} volume constraint $\hat V = 0.3 |D|$ i.e.
\begin{align}
    \theta_0 & = \int_{\Gamma} \mathbf{t} \cdot \mathbf{u} ~da \label{eq:compliance} \,, \\
    \theta_1 & = \int_{D} \hat \nu ~dV  - \hat V \label{eq:volume_comp}  \,.
\end{align}
The design domain $D$, cf. Figure \ref{fig:compliancedomain}, is subject to the traction
$\mathbf{t}=-1.0\mathbf{e}_2$ on $\Gamma_N$, the length scale parameter is $\kappa = 0.2$ and the initial design is a uniform field $\nu(\mathbf{x}) = 0.1$.
We perform four different optimizations corresponding to
uniform and nonuniform meshes with optimizations in $\mathbb{R}^n$ and $L^2$.
The uniform mesh contains 128,000 elements.
Our non-uniform mesh is illustrated in Figure \ref{fig:amr_compliance} and contains 80,577 elements.
It is important to have meshes that are sufficiently refined so the infinite-dimensional response $\mathbf{u}$ is well approximated.
We also include a highly refined arbitrary region on the top to clearly illustrate the deficiency of the NLP algorithm in $\mathbb{R}^n$.
Both meshes are included as \textit{beam\_uniform.geo} and \textit{beam\_amr.geo} in the files to reproduce the results.
In all cases, we run the optimization problems until the number of iterations reaches 200, although some designs converge sooner.

\begin{figure}
    \centering
    \input{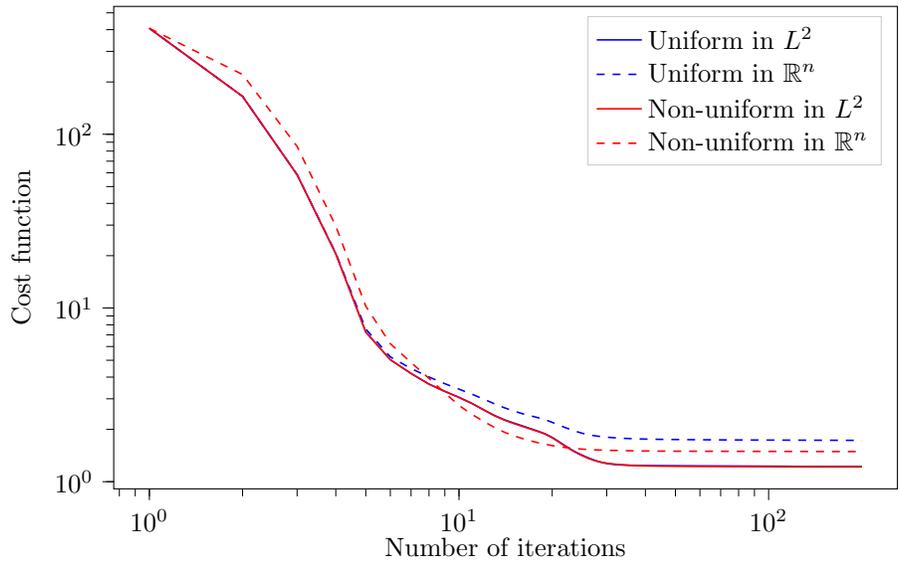}
    \caption{Cost function evolution for the compliance problem.}
    \label{fig:compliancecostfunction}
\end{figure}

\begin{figure}
    \centering
    \input{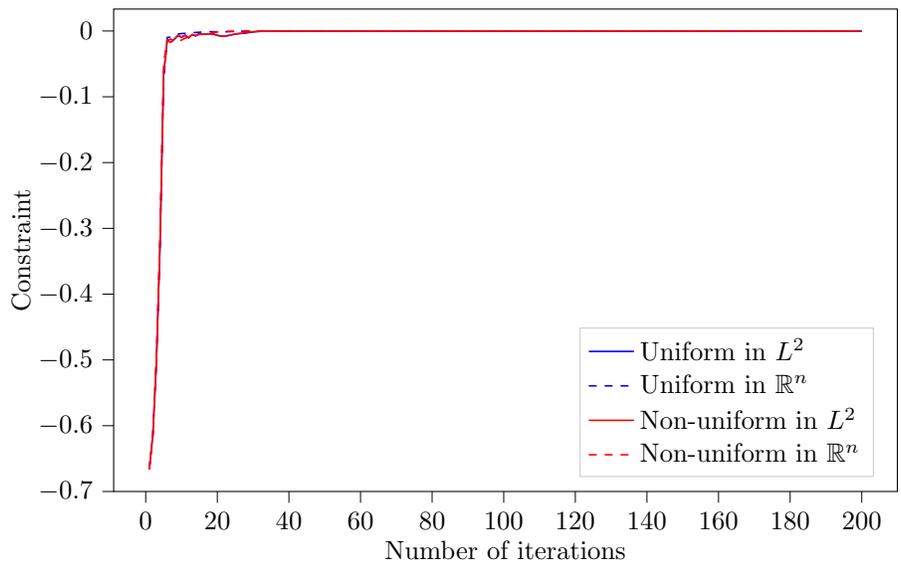}
    \caption{Constraint function evolution for the compliance problem.}
    \label{fig:complianceconstraintfunction}
\end{figure}

\begin{figure}
    \centering
    \input{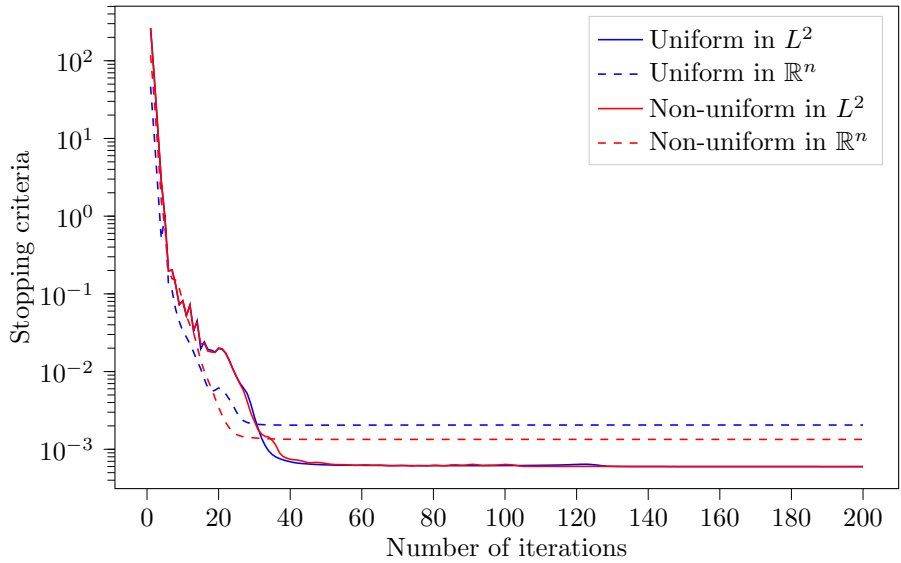}
    \caption{Convergence metric evolution for the compliance problem.}
    \label{fig:compliancekktfunction}
\end{figure}

The optimized designs in Table \ref{tab:complianceresults} show that the $\mathbb{R}^n$ NLP algorithm is mesh dependent as opposed to the $L^2$ NLP algorithm.
Figures \ref{fig:compliancecostfunction} - \ref{fig:compliancekktfunction} show the evolution of the cost and constraint functions and the convergence metric, cf. Equations \eqref{eq:kkt_stopping_final}.
It is well known that the topology optimization problem is not convex and hence different initial designs might lead to different local minima.
The examples in Table \ref{tab:complianceresults} start from the same initial design.
We therefore attribute the difference in the optimized designs to the mesh dependency of the NLP algorithm in $\mathbb{R}^n$.
\newlength{\RUno}
\newlength{\BeamW}
\setlength{\BeamW}{5cm}
\newlength{\BeamL}
\setlength{\BeamL}{2.875\BeamW}

\definecolor{naranja}{RGB}{255, 165,53}
\definecolor{rojo}{RGB}{193,64,77}
\definecolor{granate}{RGB}{139,0,63}

\begin{table}[h!]
    \centering
    \begin{tabular}{|M{0.5cm}|M{5.2cm}|M{5.2cm}|}
        \hline
                                         & Optimization in $\mathbb{R}^n$ & Optimization in $L^2$ \\ \hline
        \rotatebox{90}{Uniform mesh}     &
        \includegraphics[scale=0.12]{\MyPath/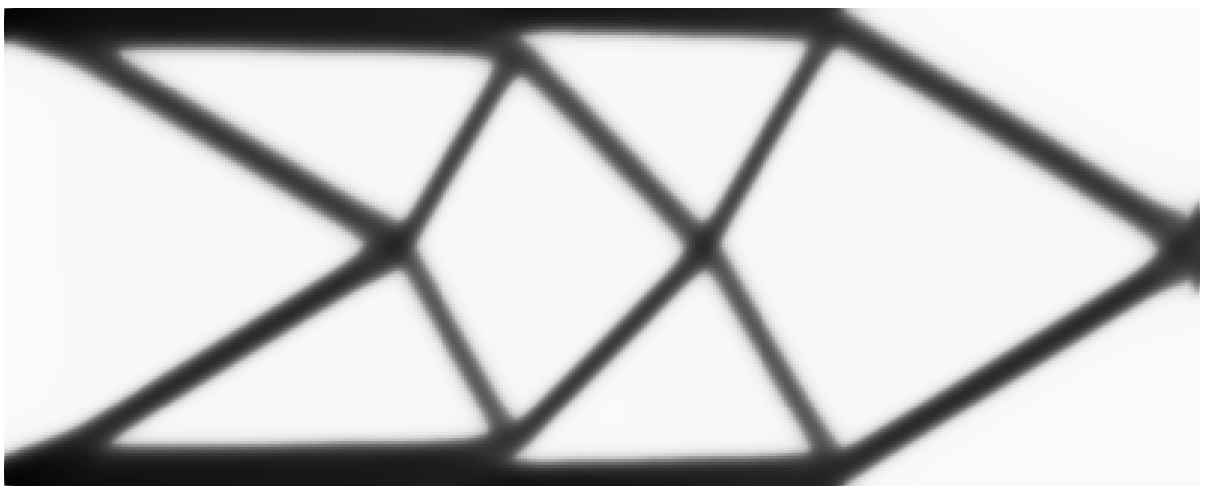}
                                         &
        \includegraphics[scale=0.12]{\MyPath/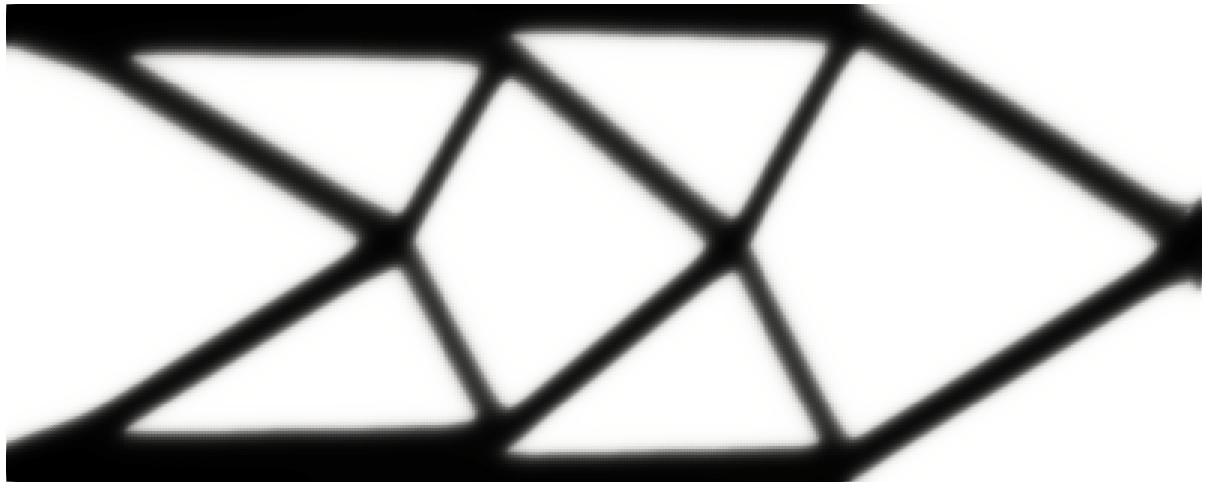}                   \\ \hline
        \rotatebox{90}{Non-uniform mesh} &
        \includegraphics[scale=0.12]{\MyPath/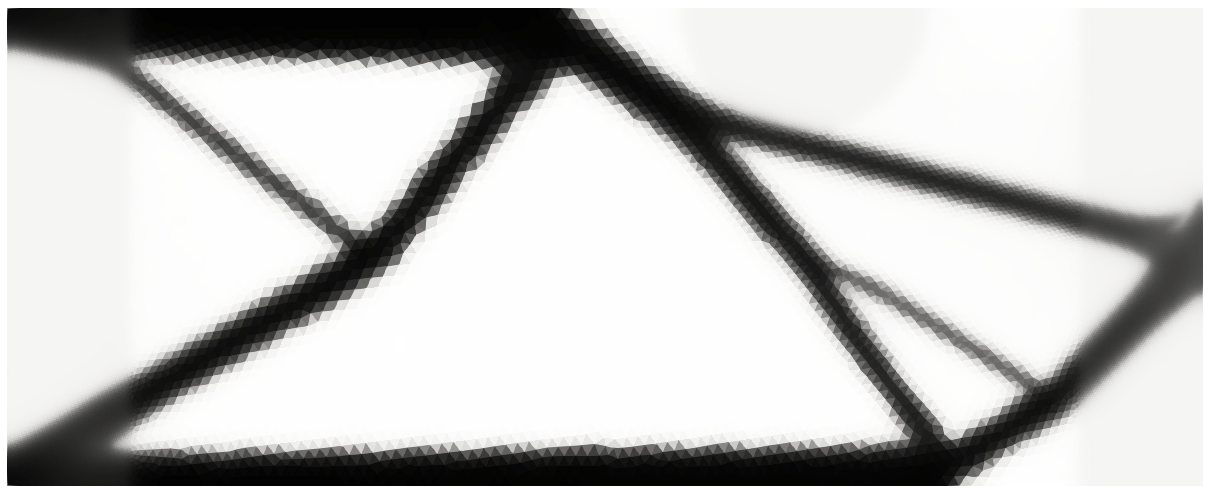}
                                         &
        \includegraphics[scale=0.12]{\MyPath/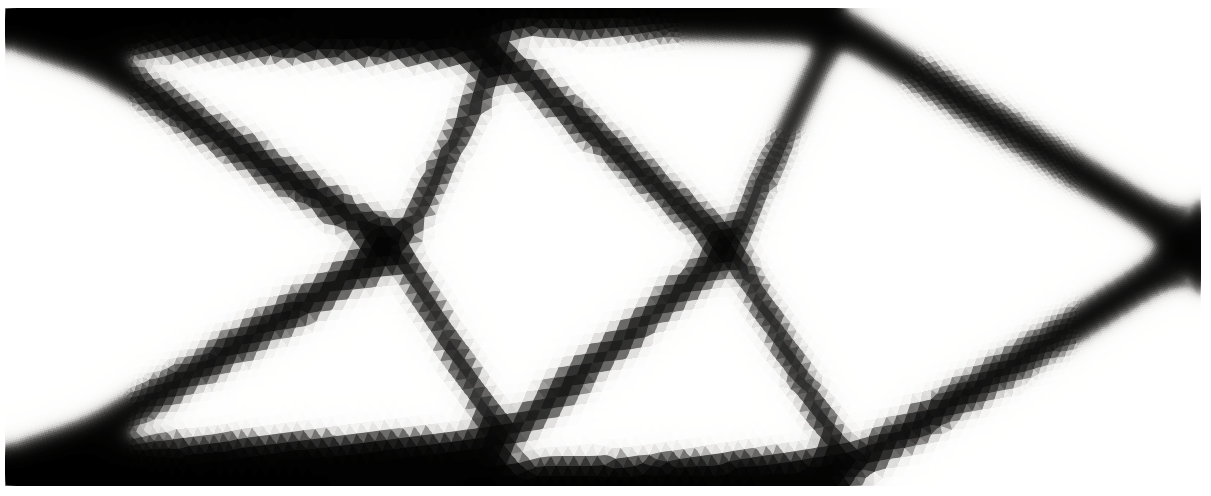}               \\ \hline
    \end{tabular}
    \caption{Optimized designs for the compliance problem.}
    \label{tab:complianceresults}
\end{table}

\setlength\figureheight{8cm}
\setlength\figurewidth{12cm}

\begin{figure}[!h]
    \center
    \tikzset{>=latex}
    \begin{tikzpicture}[  spring/.style = {decorate,decoration={zigzag,amplitude=6pt,segment length=4pt}}
            , scale=0.07, every node/.style={inner sep=0pt}]
        \small
        \node (NW) at (0,50) {};
        \draw[fill=gray!50!white] (0,0) -- (100,0) -- (100,50) -- (0,50) -- (0,0);
        \foreach \x in {0,4,...,100}
        \draw [thick] (NW) ++ (\x,2.0cm) circle (2.0cm);
        \draw[thick]  (0,54.5) -- (100, 54.5);
        \fill[pattern=north west lines]  (0,54.5) rectangle (100, 57);
        \draw ($(NW) - (0.0, 5.0)$) -- ($(NW) - (5.0, 5.0)$);
        \draw[spring] ($(NW) - (5.0, 5.0)$) -- ($(NW) - (15.0, 5.0)$);
        \draw ($(NW) - (15.0, 5.0)$) -- ($(NW) - (20.0, 5.0)$);
        \node (Gamma1) at ($(NW) - (5, 10)$) {$\Gamma_1$};
        \draw ($(NW) - (20.0, 10.0)$) -- ($(NW) - (20.0, 0.0)$);
        \draw[very thick] (NW) -- ($(NW) + (0.0, -10.0)$);
        \fill[pattern=north west lines] ($(NW) - (25.0, 10.0)$) rectangle ($(NW) - (20.0, 0.0)$);
        \node (NE) at (100,50) {};
        \draw ($(NE) - (0.0, 5.0)$) -- ($(NE) + (5.0, -5.0)$);
        \draw[spring] ($(NE) + (5.0, -5.0)$) -- ($(NE) + (15.0, -5.0)$);
        \draw ($(NE) + (15.0, -5.0)$) -- ($(NE) + (20.0, -5.0)$);
        \node (Gamma2) at ($(NE) + (5, -10)$) {$\Gamma_2$};
        \draw ($(NE) + (20.0, 0.0)$) -- ($(NE) + (20.0, -10.0)$);
        \draw[very thick] (NE) -- ($(NE) + (0.0, -10.0)$);
        \fill[pattern=north west lines] ($(NE) + (20.0, 0.0)$) rectangle ($(NE) + (25.0, -10.0)$);
        \node (SW) at (0,0) {};
        \node (SE) at (100,0) {};
        \node (GammaD) at ($(SW) + (5, 2)$) {$\Gamma_D$};
        \fill[pattern=north west lines] ($(SW) + (0.0, 2.0)$) rectangle ($(SW) - (5.0, 0.0)$);
        \dimline[line style = {line width=0.7},
            extension start length=10cm,
            extension end length=10cm]{($(SW) - (10cm, 0.0)$)}{($(SW) + (-10cm, 2)$)}{$2$};
        \foreach \y in {0, 2.5, 5.0, 7.5, 10.0}
        \draw [-{Latex[width=1mm]}] ($(NW) - (0.0, \y)$) -- ($(NW) + (15.0, -\y)$);
        \normalsize
        \node at ($(NW) + (20.0, -5.0)$) {$f_{x}$};
        \foreach \y in {0, 2.5, 5.0, 7.5, 10.0}
        \draw [color=red, -{Latex[width=1mm]}] ($(NE) - (0.0, \y)$) -- ($(NE) - (15.0, \y)$);
        \node at ($(NE) - (25.0, 5.0)$) {$u_{out}$};
        \small
        \node (A) at ($(SW) - (0.0, 5.0)$) {};
        \node (B) at ($(SE) - (0.0, 5.0)$) {};
        \dimline[line style = {line width=0.7},
            extension start length=-5cm,
            extension end length=-5cm]{(A)}{(B)}{120};
        \node (C) at ($(NE) + (50.0, 0.0)$) {};
        \node (D) at ($(SE) + (50.0, 0.0)$) {};
        \dimline[line style = {line width=0.7},
            extension start length=-50cm,
            extension end length=-50cm]{(D)}{(C)}{60};
        \node (E) at ($(C) + (-10.0, 0.0)$) {};
        \node (F) at ($(E) + (0.0, -10.0)$) {};
        \dimline[line style = {line width=0.7},
            extension start length=-40cm,
            extension end length=-40cm]{(F)}{(E)}{10};
        \node (G) at ($(NW) + (-30.0, 0.0)$) {};
        \node (H) at ($(G) + (0.0, -10.0)$) {};
        \dimline[line style = {line width=0.7},
            extension start length=-30cm,
            extension end length=-30cm]{(G)}{(H)}{10};
    \end{tikzpicture}
    \caption{Design domain and boundary conditions for the compliant mechanism problem.
        Domain symmetry is used whereby only the lower half of the structure is analyzed.}
    \label{fig:mechanismdomain}
    \normalsize
\end{figure}
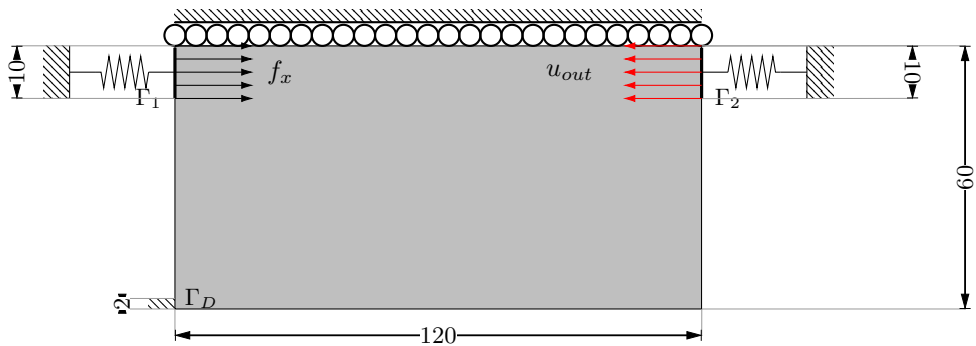

We next benchmark our $L^2$ algorithm with the compliant mechanism design problem, cf. Figure \ref{fig:mechanismdomain} where we minimize the horizontal displacement $u=\mathbf{u}\cdot\mathbf{e}_1$ in the output port $\Gamma_2$, subject to a maximum volume constraint $\hat V = 0.3 |D|$, viz.
\begin{align}
    \theta_0 & = \int_{\Gamma_2} u ~da \label{eq:compliant} \,,           \\
    \theta_1 & = \int_{D} \hat \nu ~dV - \hat V \,.\label{eq:volume_mech}
\end{align}
Consistent with \cite{de2015stress}, we introduce Robin boundary conditions into the formulation \eqref{eq:elastic_optimization}.
\begin{equation}
    \begin{aligned}
        \mathbf{n} \cdot \mathbb{C}[\nabla \mathbf{u}]  \mathbf{n} & = -k_{in} (\mathbf{u}\cdot \mathbf{n}) + f_x & ~ \text{on} ~\Gamma_1 \,,   \\
        \mathbf{n} \cdot \mathbb{C}[\nabla \mathbf{u}]  \mathbf{n} & = -k_{out} (\mathbf{u}\cdot \mathbf{n})      & ~ \text{on} ~\Gamma_2   \,,
    \end{aligned}
\end{equation}
where $f_x=10$ and the spring coefficients are $k_{in}=\frac{1}{3}$ and $k_{out}=\frac{0.001}{3}$.
Figures \ref{fig:mechanismdomain} and \ref{fig:amr_mechanism} illustrate the design domain and the non-uniform mesh.
The number of elements in the uniform (\textit{mechanism\_uniform.geo}) and non-uniform (\textit{mechanism\_amr.geo}) meshes are 57,600 and 199,404.
The initial design is again $\nu(\mathbf{x}) = 0.1$.
We run all the optimization problems until the number of iterations reaches 200, for which almost all the problems converged to optimized designs.

\begin{figure}
    \centering
    \input{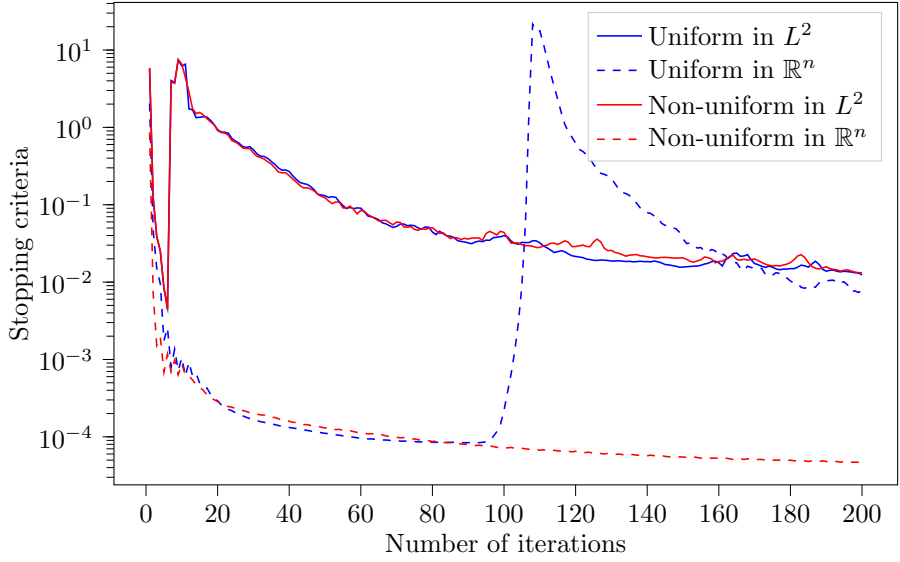}
    \caption{Convergence metric evolution for the compliant mechanism problem.}
    \label{fig:mechanismkktfunction}
\end{figure}

The optimized designs for the length scale $\kappa=0.8$, summarized and illustrated
in Table \ref{tab:mechanismresults} and Figures \ref{fig:mechanismcostfunction} - \ref{fig:mechanismkktfunction},
again show that the original $\mathbb{R}^n$ NLP algorithm is mesh dependent.
Notably, the design on the non-uniform mesh with the $\mathbb{R}^n$ algorithm does not even converge in 200 iterations.

\newlength{\Long}
\setlength{\Long}{12cm}
\newlength{\RWidth}
\setlength{\RWidth}{0.5\Long}

\begin{table}[h]
    \centering
    \begin{tabular}{|M{0.5cm}|M{6.0cm}|M{6.0cm}|}
        \hline
                                         & Optimization in $\mathbb{R}^n$ & Optimization in $L^2$ \\ \hline
        \rotatebox{90}{Uniform mesh}     &
        \includegraphics[scale=0.12]{\MyPath/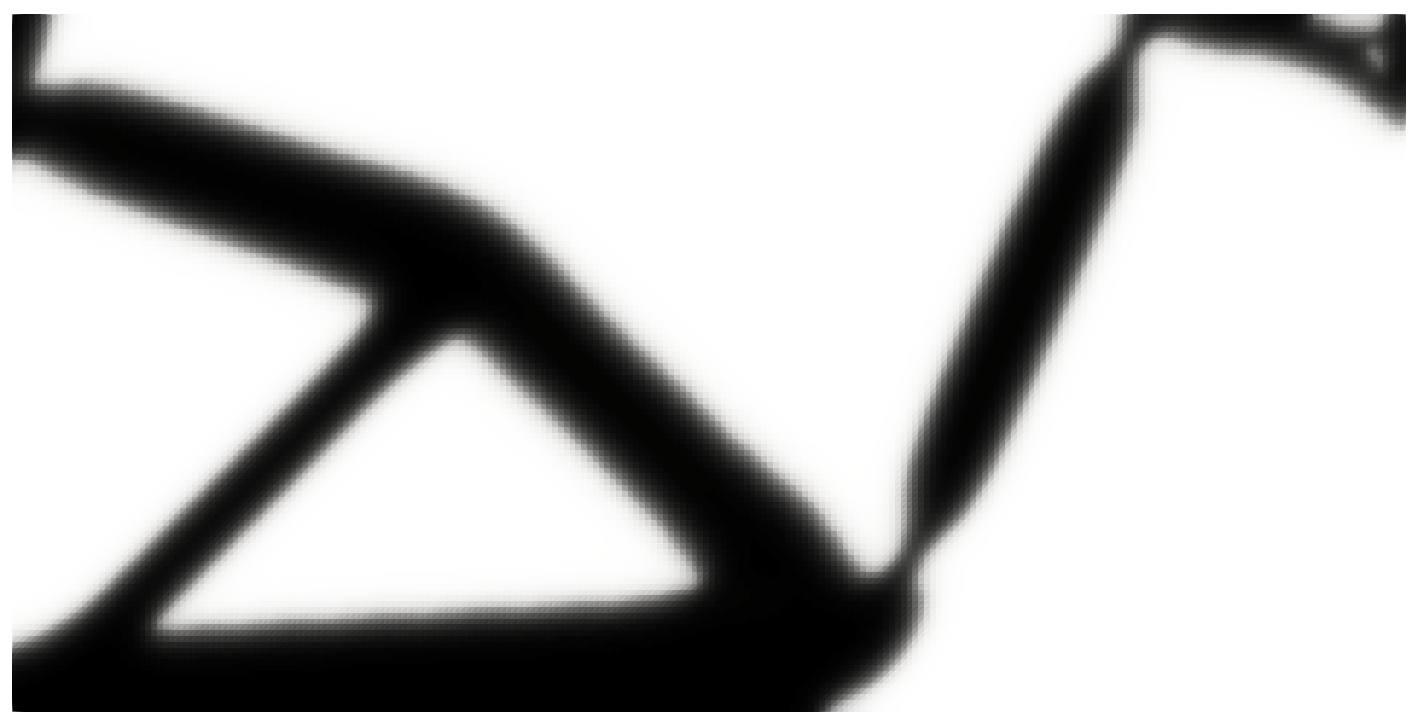}
                                         &
        \includegraphics[scale=0.12]{\MyPath/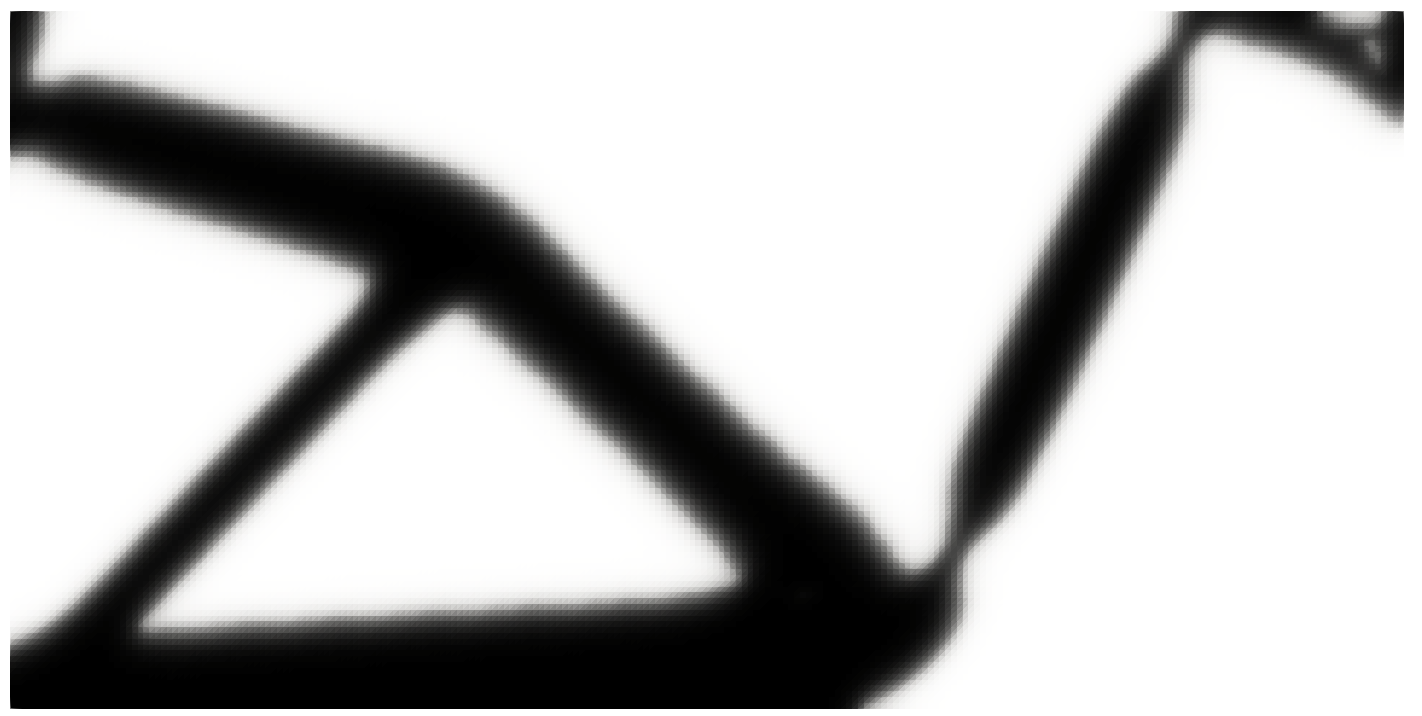}                    \\ \hline
        \rotatebox{90}{Non-uniform mesh} &
        \includegraphics[scale=0.12]{\MyPath/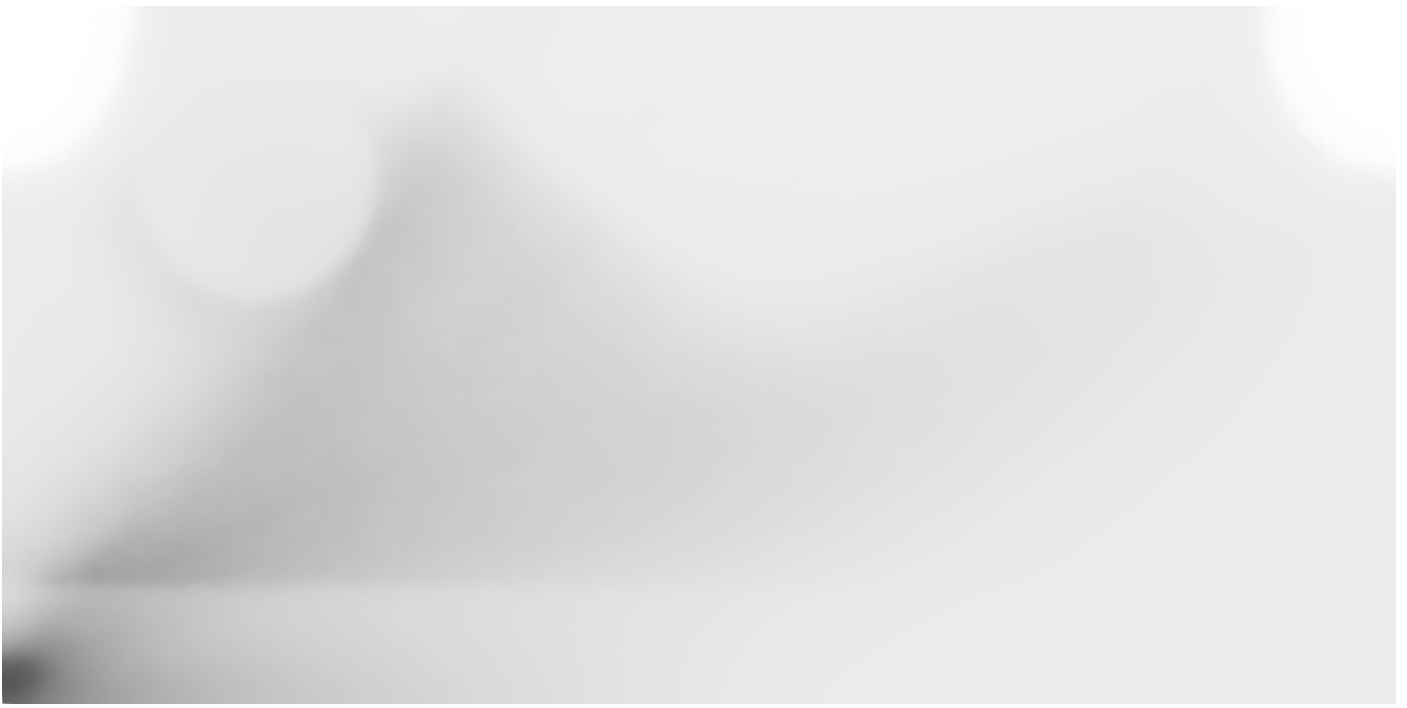}
                                         &
        \includegraphics[scale=0.12]{\MyPath/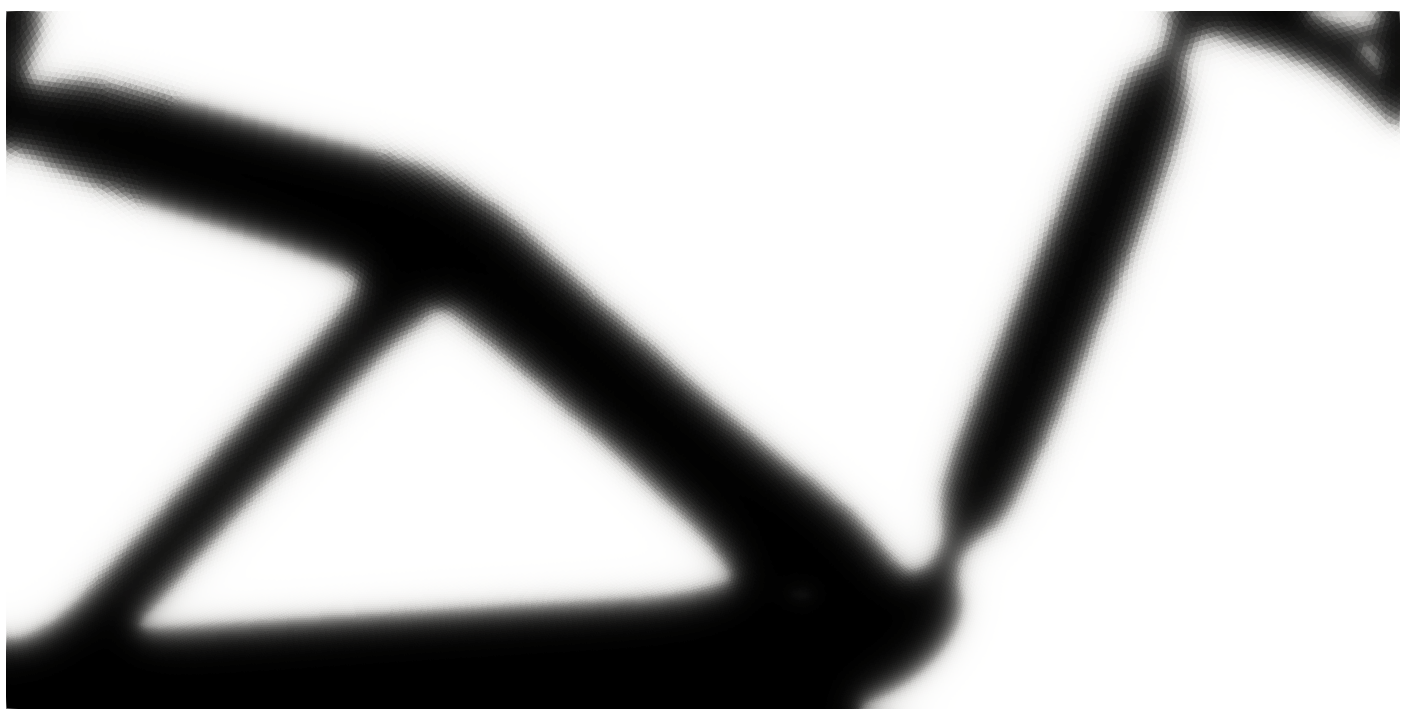}                \\ \hline
    \end{tabular}
    \caption{Optimized designs for the compliant mechanism problem.}
    \label{tab:mechanismresults}
\end{table}

\setlength\figureheight{8cm}
\setlength\figurewidth{12cm}
\begin{figure}
    \centering
    \input{\MyPath/tex_figures-mechanism-convergence_plots1.tex}
    \caption{Cost function evolution for the compliant mechanism problem.}
    \label{fig:mechanismcostfunction}
\end{figure}

\begin{figure}
    \centering
    \input{\MyPath/tex_figures-mechanism-convergence_plots2.tex}
    \caption{Constraint function evolution for the compliant mechanism problem.}
    \label{fig:mechanismconstraintfunction}
\end{figure}

Our next example is the stress-constrained problem where the goal is to minimize the volume of an L-bracket subject to a \rojo{maximum}
pointwise constraint in the Von Mises stress field \rojo{$\sigma_{VM} \leq \sigma_y$}, cf. Figure \ref{fig:stressdomain}.
We follow the formulation in \cite{novotny}, where the stress constraint is imposed via a penalty method, i.e.
our unconstrained problem uses the cost function
\begin{align}
    \theta_0 & = \int_{D} \hat \nu ~da+ \gamma  \norm{\sigma_{VM} - \sigma_y}_+  \,,
\end{align}
\adapt{where the penalty parameter is $\gamma=10$,}
\begin{align}
    \norm{\sigma_{VM} - \sigma_y}_+ = \int_{D} R_p\left(\frac{\sigma_{VM}}{\sigma_y}\right)dV \,,
\end{align}
\begin{align}
    R_p(x) = (1 + (x)^p)^{\frac{1}{p}} - 1\,,
    \label{eq:penalty_func}
\end{align}
\adapt{$\sigma_y=1.5$ and $p=8$, for the first 300 iterations and $p=20$ for the remaining 100.}
The relaxed stress formulation \cite{chau} uses
\begin{align}
    \sigma_{\eta} = \eta_c \mathbb{C}[\nabla \mathbf{u}] \,,
\end{align}
to calculate the Von Mises stress $\sigma_{VM} = \sqrt{\frac{2}{3} \sigma_{\eta}^{'}:\sigma_{\eta}^{'}}$, where $\sigma_{\eta}^{'} = \sigma_{\eta} - \frac{2}{3} \text{tr}(\sigma_{\eta}) \textbf{I}$ and
\begin{align}
    \eta_c (\hat\nu) = \hat\nu^{0.5}\,,
\end{align}
The filter parameter is $\kappa=1.2$ and the initial design is $\nu(\mathbf{x}) = 0.5$.
To obtain a better design we extend the domain by adding a region $\Omega_i$ of finite elements near the reentrant corner, cf. Figure \ref{fig:amr_stress}.
This region however, is excluded in the design by enforcing the constraint $\int_{\Omega_i} \hat{\nu}~dV \leq 0$.
This added region lessens the boundary effect of the filter operation in the reentrant corner region which otherwise adversely affects our results \cite{wallin2020consistent}.
The non-uniform mesh (file \textit{lbracket\_amr.geo}) contains 53,122 elements, whereas the uniform mesh (\textit{lbracket\_uniform.geo}) contains 139,264 elements.

\begin{figure}
    \centering
    \input{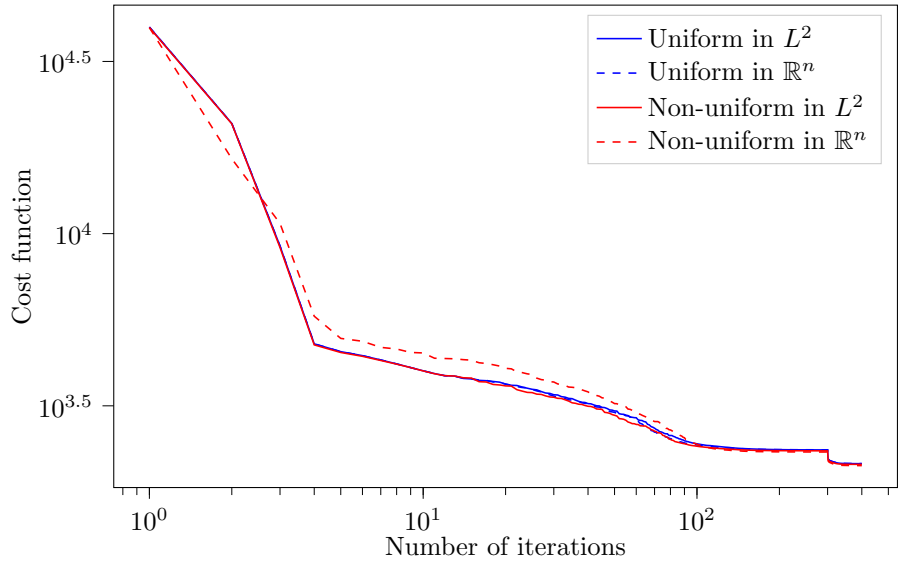}
    \caption{Cost function evolution for the stress-constrained problem.}
    \label{fig:stresscostfunction}
\end{figure}

The optimized designs in Table \ref{tab:stressresults} again illustrate the mesh dependence of the original $\mathbb{R}^n$ NLP algorithm.
This time however, the difference in the cost function values between the $\mathbb{R}^n$ and $L^2$ designs is barely noticeable, cf. the log plot in Figure \ref{fig:stresscostfunction}.
Due to the high oscillation in the convergence metric, we plot them separately in Figures \ref{fig:stresskktfunction1} - \ref{fig:stresskktfunction4}.

\newlength{\ChauWidth}
\setlength{\ChauWidth}{10cm}
\newlength{\ChauHeight}
\setlength{\ChauHeight}{\ChauWidth}
\newlength{\ChauCorner}
\setlength{\ChauCorner}{0.4\ChauWidth}

\begin{figure}[!h]
    \center
    \tikzset{>=latex}
    \begin{tikzpicture}[every node/.style={inner sep=0pt}]
        \fill[pattern=north west lines]  (-0.05\ChauCorner,\ChauHeight) rectangle (1.05\ChauCorner, 1.05\ChauHeight);
        \draw[very thick] (-0.05\ChauCorner,\ChauHeight) -- (1.05\ChauCorner,\ChauHeight);
        \draw[fill=gray!50!white] (0,0) -- (\ChauWidth,0) -- (\ChauWidth,\ChauCorner) -- (\ChauCorner,\ChauCorner) -- (\ChauCorner,\ChauHeight) -- (0,\ChauHeight) -- (0,0);
        \node (GammaD) at (0.5\ChauCorner, 0.95\ChauHeight) {$\Gamma_D$};
        \node (A) at (0, -0.5) {};
        \node (B) at (\ChauWidth, -0.5) {};
        \dimline[line style = {line width=0.7},
            extension start length=-1cm,
            extension end length=-1cm]{(A)}{(B)}{$100$};
        \node (C) at (-0.5, 0) {};
        \node (D) at (-0.5, \ChauHeight) {};
        \dimline[line style = {line width=0.7},
            extension start length=1cm,
            extension end length=1cm]{(C)}{(D)}{$100$};
        \node (F) at (0, 1.1\ChauHeight) {};
        \node (E) at (\ChauCorner, 1.1\ChauHeight) {};
        \dimline[line style = {line width=0.7},
            extension start length=1cm,
            extension end length=1cm]{(F)}{(E)}{$40$};
        \node (G) at (1.05\ChauWidth, 0) {};
        \node (H) at (1.05\ChauWidth, \ChauCorner) {};
        \dimline[line style = {line width=0.7},
            extension start length=-1cm,
            extension end length=-1cm]{(G)}{(H)}{$40$};
        \node (Gamma) at (0.97\ChauWidth, 0.95\ChauCorner) {$\Gamma_N$};
        \node (L1) at (\ChauWidth, \ChauCorner) {};
        \node (L1p) at (\ChauWidth, 1.2\ChauCorner) {};
        \node (L2) at (0.97\ChauWidth, \ChauCorner) {};
        \node (L2p) at (0.97\ChauWidth, 1.2\ChauCorner) {};
        \node (L3) at (0.94\ChauWidth, \ChauCorner) {};
        \node (L3p) at (0.94\ChauWidth, 1.2\ChauCorner) {};
        \draw [thick,->] (L1p) -- (L1);
        \draw [thick,->] (L2p) -- (L2);
        \draw [thick,->] (L3p) -- (L3);
        \node (L4p) at (0.94\ChauWidth, 1.4\ChauCorner) {};
        \node (L5p) at (\ChauWidth, 1.4\ChauCorner) {};
        \dimline[line style = {line width=0.7},
            extension start length=1.6cm,
            extension end length=1.6cm]{(L4p)}{(L5p)}{$5$};
    \end{tikzpicture}
    \caption{Intended design domain for the stress-constrained problem.}
    \label{fig:stressdomain}
\end{figure}
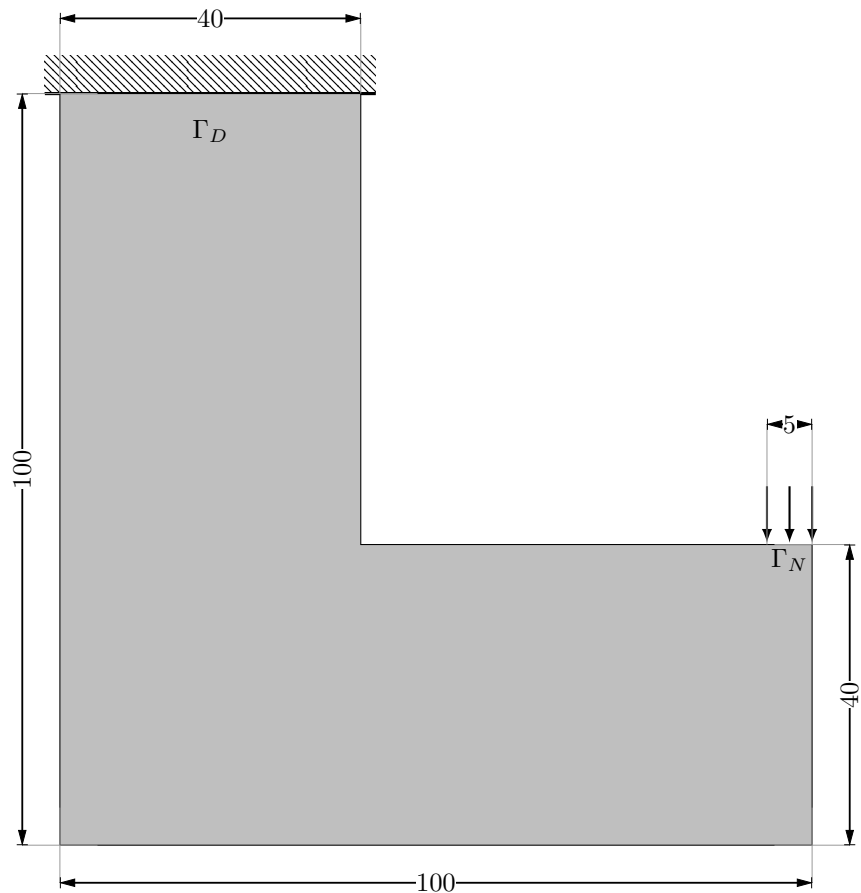

\setlength\figureheight{8cm}
\setlength\figurewidth{12cm}

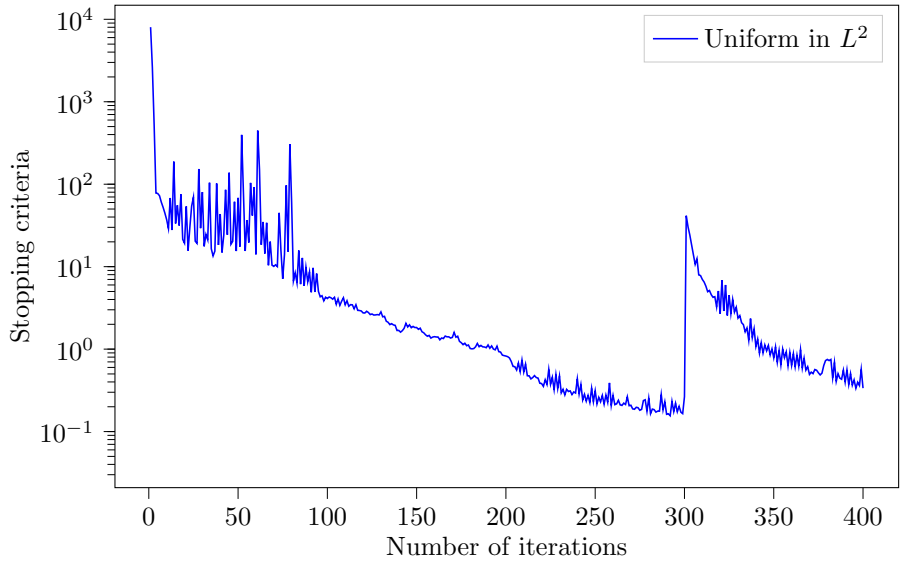
\begin{figure}
    \centering
\begin{tikzpicture}

\begin{axis}[
height=\figureheight,
legend cell align={left},
legend style={fill opacity=0.8, draw opacity=1, text opacity=1, draw=white!80!black},
log basis y={10},
tick align=outside,
tick pos=left,
width=\figurewidth,
x grid style={white!69.0196078431373!black},
xlabel={Number of iterations},
xmin=-18.95, xmax=419.95,
xtick style={color=black},
y grid style={white!69.0196078431373!black},
ylabel={Stopping criteria},
ymin=0.0208917309728883, ymax=14811.3521698619,
ymode=log,
ytick style={color=black}
]
\addplot [semithick, blue]
table {%
1 8028.925399
2 2605.876592
3 517.704124
4 78.114975
5 77.163357
6 72.468086
7 60.628845
8 52.512107
9 44.693758
10 37.159462
11 29.086105
12 68.129654
13 27.802263
14 188.042896
15 33.408214
16 55.545779
17 31.546016
18 75.350929
19 21.426239
20 19.306069
21 53.872393
22 15.64046
23 30.094505
24 54.396599
25 69.463212
26 20.383019
27 19.180867
28 153.071904
29 29.282711
30 80.276813
31 17.640667
32 24.571812
33 21.182625
34 104.348453
35 16.48841
36 13.490231
37 15.448761
38 102.14621
39 18.349419
40 43.510091
41 14.853069
42 24.219231
43 85.650718
44 24.379259
45 138.011072
46 18.72104
47 20.316216
48 61.306352
49 15.456495
50 68.351734
51 17.486731
52 396.795034
53 78.802729
54 15.47769
55 36.499389
56 19.612948
57 103.836182
58 41.63418
59 92.457986
60 14.089972
61 447.07827
62 149.558677
63 18.430473
64 34.896081
65 14.430171
66 34.065151
67 10.390352
68 20.12827
69 10.426098
70 10.078999
71 10.528418
72 10.031845
73 44.958104
74 15.688254
75 7.121874
76 14.11998
77 97.037983
78 15.213889
79 305.597791
80 59.483203
81 6.714559
82 8.428346
83 6.355715
84 15.712138
85 6.136489
86 12.776625
87 5.860516
88 9.838498
89 6.582891
90 8.601384
91 4.889627
92 9.677203
93 4.929704
94 8.246828
95 5.018433
96 4.310619
97 4.439844
98 3.854522
99 4.260818
100 4.107506
101 4.285648
102 4.196578
103 4.033456
104 4.249594
105 3.495005
106 4.058536
107 3.39528
108 3.813777
109 4.198821
110 3.451609
111 3.859816
112 3.309526
113 3.463856
114 3.445533
115 3.083011
116 3.453301
117 2.953243
118 2.939106
119 2.902791
120 2.749595
121 2.741131
122 2.878371
123 2.781228
124 2.622418
125 2.668152
126 2.58992
127 2.599111
128 2.623643
129 2.597262
130 2.828106
131 2.46573
132 2.496729
133 2.180317
134 2.116824
135 1.982713
136 2.03667
137 1.96574
138 1.94232
139 1.68334
140 1.68562
141 1.605996
142 1.674884
143 1.788178
144 2.053559
145 1.857149
146 1.970203
147 1.81851
148 1.885963
149 1.827418
150 1.828967
151 1.738145
152 1.795363
153 1.624288
154 1.576996
155 1.491486
156 1.433835
157 1.467368
158 1.360841
159 1.400778
160 1.415929
161 1.403168
162 1.401721
163 1.296099
164 1.364431
165 1.343773
166 1.440052
167 1.409298
168 1.410306
169 1.362356
170 1.38745
171 1.593338
172 1.396347
173 1.43359
174 1.241812
175 1.196207
176 1.13274
177 1.178639
178 1.09972
179 1.121405
180 1.01727
181 1.006986
182 1.015628
183 1.066323
184 1.171545
185 1.06929
186 1.109584
187 1.066158
188 1.061324
189 1.040121
190 1.124079
191 1.025206
192 1.09352
193 1.001695
194 0.97997
195 1.086019
196 0.932175
197 0.92103
198 0.846354
199 0.833824
200 0.824019
201 0.809924
202 0.78127
203 0.707679
204 0.617763
205 0.612135
206 0.558319
207 0.677934
208 0.534067
209 0.737518
210 0.528324
211 0.656298
212 0.475567
213 0.474297
214 0.435205
215 0.451669
216 0.477447
217 0.451973
218 0.445323
219 0.385502
220 0.382354
221 0.353503
222 0.4254
223 0.372994
224 0.567037
225 0.374755
226 0.455794
227 0.322966
228 0.491027
229 0.334022
230 0.470555
231 0.294708
232 0.325588
233 0.278065
234 0.324934
235 0.306374
236 0.311574
237 0.281142
238 0.302583
239 0.293526
240 0.435743
241 0.285303
242 0.371568
243 0.241913
244 0.282798
245 0.233702
246 0.274141
247 0.22539
248 0.334831
249 0.230562
250 0.312207
251 0.216198
252 0.25981
253 0.229721
254 0.273957
255 0.214665
256 0.277507
257 0.227109
258 0.3896
259 0.207068
260 0.274521
261 0.211999
262 0.216674
263 0.239436
264 0.210708
265 0.207529
266 0.221391
267 0.213171
268 0.264251
269 0.20601
270 0.208574
271 0.188427
272 0.186423
273 0.197112
274 0.194841
275 0.181262
276 0.18581
277 0.237755
278 0.243018
279 0.176659
280 0.253936
281 0.165071
282 0.187381
283 0.184688
284 0.172214
285 0.176918
286 0.176939
287 0.280679
288 0.182559
289 0.22973
290 0.162017
291 0.16418
292 0.155756
293 0.232467
294 0.175078
295 0.220294
296 0.174061
297 0.202073
298 0.171862
299 0.165155
300 0.264995
301 41.480786
302 30.359866
303 24.166548
304 18.351932
305 13.958846
306 10.68379
307 12.467426
308 7.948391
309 7.796709
310 7.004506
311 6.511009
312 5.796985
313 4.968729
314 5.137951
315 4.574233
316 4.246677
317 4.321928
318 3.257201
319 5.074061
320 2.686534
321 6.839474
322 2.899874
323 6.000629
324 2.557372
325 4.515358
326 2.691833
327 3.986457
328 2.815915
329 3.219459
330 2.377221
331 2.546087
332 2.089996
333 1.98268
334 1.621375
335 1.80452
336 1.280453
337 2.346408
338 1.379377
339 1.762154
340 1.077709
341 1.322033
342 0.895643
343 1.194628
344 0.898163
345 1.123129
346 0.954668
347 1.106451
348 0.848781
349 1.007196
350 0.784543
351 1.061388
352 0.693106
353 0.960008
354 0.642527
355 0.999997
356 0.786851
357 0.925016
358 0.626505
359 0.909936
360 0.63579
361 0.877348
362 0.610355
363 0.842244
364 0.621551
365 0.982636
366 0.588993
367 0.776431
368 0.556137
369 0.612359
370 0.493222
371 0.525737
372 0.504273
373 0.565808
374 0.554882
375 0.516538
376 0.486463
377 0.511039
378 0.63217
379 0.726446
380 0.7502
381 0.726073
382 0.749376
383 0.448669
384 0.688826
385 0.406232
386 0.503513
387 0.44513
388 0.431299
389 0.555237
390 0.397824
391 0.559683
392 0.371686
393 0.502958
394 0.36521
395 0.423651
396 0.335156
397 0.396884
398 0.36224
399 0.567301
400 0.336291
};
\addlegendentry{Uniform in $L^2$}
\end{axis}

\end{tikzpicture}
    \caption{Convergence metric evolution for the stress-constrained problem with uniform mesh in $L^2$.}
    \label{fig:stresskktfunction1}
\end{figure}

\begin{figure}
    \centering
    \begin{tikzpicture}

\begin{axis}[
height=\figureheight,
legend cell align={left},
legend style={fill opacity=0.8, draw opacity=1, text opacity=1, draw=white!80!black},
log basis y={10},
tick align=outside,
tick pos=left,
width=\figurewidth,
x grid style={white!69.0196078431373!black},
xlabel={Number of iterations},
xmin=-18.95, xmax=419.95,
xtick style={color=black},
y grid style={white!69.0196078431373!black},
ylabel={Stopping criteria},
ymin=0.0208917309728883, ymax=14811.3521698619,
ymode=log,
ytick style={color=black}
]
\addplot [semithick, blue, dashed]
table {%
1 1774.965649
2 577.143567
3 123.290146
4 17.898987
5 52.02411
6 56.124282
7 48.640839
8 42.989018
9 34.804679
10 27.059813
11 16.198127
12 15.836332
13 15.600422
14 36.163132
15 8.829455
16 17.655051
17 7.867892
18 22.09521
19 5.449375
20 6.157203
21 6.11664
22 8.217662
23 12.893207
24 3.710212
25 5.690136
26 11.993539
27 7.894414
28 10.044579
29 8.237891
30 7.505265
31 6.212948
32 26.485747
33 4.517381
34 17.163318
35 3.371413
36 7.312246
37 2.572619
38 5.138973
39 6.687416
40 6.668997
41 17.412789
42 9.020162
43 15.110429
44 2.462375
45 5.934328
46 1.978539
47 8.909177
48 52.775806
49 10.262199
50 2.741393
51 18.493472
52 5.565242
53 29.636072
54 6.176587
55 18.262801
56 6.201801
57 3.499386
58 10.075467
59 7.304276
60 2.743677
61 9.06676
62 9.308697
63 7.377299
64 3.664874
65 6.433109
66 4.44691
67 6.771151
68 3.701246
69 5.281009
70 4.430718
71 13.043699
72 2.804418
73 2.249316
74 2.160042
75 3.189115
76 2.611501
77 2.747663
78 3.734907
79 3.096839
80 1.969964
81 3.295751
82 1.738519
83 2.272624
84 2.745892
85 7.348082
86 1.88783
87 1.400156
88 3.229889
89 1.38942
90 2.422396
91 1.417842
92 1.267134
93 1.419443
94 1.235471
95 1.068865
96 1.148703
97 1.076969
98 1.293754
99 0.908646
100 1.056786
101 1.235544
102 1.227389
103 1.340976
104 1.080009
105 0.893224
106 0.756924
107 0.744024
108 1.051629
109 0.850984
110 0.748944
111 0.720906
112 0.841403
113 0.926158
114 0.83831
115 0.726709
116 0.746905
117 0.820519
118 0.814333
119 0.715067
120 0.822503
121 0.807092
122 0.742277
123 0.640344
124 0.698329
125 0.656379
126 0.649246
127 0.542856
128 0.570271
129 0.507302
130 0.526083
131 0.476399
132 0.564841
133 0.521337
134 0.495447
135 0.510457
136 0.538666
137 0.508926
138 0.524051
139 0.519979
140 0.532825
141 0.503188
142 0.504978
143 0.559216
144 0.492493
145 0.466367
146 0.47493
147 0.465524
148 0.454954
149 0.446992
150 0.395206
151 0.388337
152 0.398939
153 0.375056
154 0.370309
155 0.353935
156 0.367926
157 0.333265
158 0.311421
159 0.289425
160 0.281904
161 0.275167
162 0.261915
163 0.246197
164 0.231134
165 0.218005
166 0.225156
167 0.212459
168 0.225298
169 0.190328
170 0.221349
171 0.184509
172 0.221844
173 0.17963
174 0.197051
175 0.155446
176 0.181031
177 0.152269
178 0.19686
179 0.147292
180 0.169812
181 0.130336
182 0.137121
183 0.12127
184 0.140676
185 0.121825
186 0.126595
187 0.117561
188 0.115063
189 0.107986
190 0.112081
191 0.099189
192 0.136026
193 0.102226
194 0.151904
195 0.101113
196 0.122319
197 0.103922
198 0.173607
199 0.101188
200 0.13376
201 0.082426
202 0.088277
203 0.083104
204 0.07768
205 0.081199
206 0.078662
207 0.074208
208 0.079519
209 0.07574
210 0.099716
211 0.079198
212 0.106318
213 0.074357
214 0.081876
215 0.06817
216 0.065119
217 0.074474
218 0.108354
219 0.067539
220 0.1287
221 0.066918
222 0.096349
223 0.066004
224 0.073905
225 0.071053
226 0.059402
227 0.061497
228 0.063067
229 0.057185
230 0.061087
231 0.057501
232 0.07152
233 0.053318
234 0.092362
235 0.055198
236 0.077693
237 0.056183
238 0.059778
239 0.08719
240 0.056934
241 0.079423
242 0.051811
243 0.051828
244 0.053902
245 0.059039
246 0.049548
247 0.059587
248 0.050088
249 0.054359
250 0.056583
251 0.050085
252 0.05879
253 0.047983
254 0.09021
255 0.04727
256 0.081701
257 0.048648
258 0.052771
259 0.046028
260 0.054551
261 0.076627
262 0.051231
263 0.07508
264 0.052134
265 0.052106
266 0.054028
267 0.046157
268 0.048214
269 0.050001
270 0.047743
271 0.045239
272 0.053369
273 0.04553
274 0.045656
275 0.072233
276 0.04701
277 0.063342
278 0.044507
279 0.0761
280 0.047281
281 0.064928
282 0.041912
283 0.043094
284 0.049201
285 0.042165
286 0.044119
287 0.043656
288 0.041611
289 0.049141
290 0.041601
291 0.046065
292 0.042244
293 0.038892
294 0.040409
295 0.050911
296 0.052516
297 0.058366
298 0.045356
299 0.055877
300 0.03854
301 11.927856
302 7.407127
303 5.57854
304 4.173444
305 3.232833
306 2.431787
307 2.222398
308 1.869522
309 1.747111
310 1.565283
311 1.571222
312 1.254225
313 1.22785
314 0.983042
315 1.126415
316 0.778973
317 1.075075
318 0.736697
319 1.342598
320 0.726779
321 1.468895
322 0.627836
323 1.297376
324 0.773774
325 1.032
326 0.864566
327 0.795509
328 0.802198
329 0.699938
330 0.64829
331 0.609414
332 0.489104
333 0.502632
334 0.501314
335 0.407778
336 0.407245
337 0.330934
338 0.292531
339 0.287139
340 0.244215
341 0.319745
342 0.367624
343 0.434548
344 0.328329
345 0.437332
346 0.215428
347 0.324424
348 0.163287
349 0.251298
350 0.196588
351 0.30228
352 0.211924
353 0.288654
354 0.199571
355 0.233411
356 0.189132
357 0.217213
358 0.170515
359 0.191752
360 0.139975
361 0.175935
362 0.1471
363 0.206347
364 0.260522
365 0.268
366 0.20238
367 0.230207
368 0.12928
369 0.167059
370 0.11598
371 0.17849
372 0.129156
373 0.1863
374 0.124546
375 0.151595
376 0.131473
377 0.123968
378 0.117784
379 0.127042
380 0.100756
381 0.17293
382 0.102527
383 0.177614
384 0.101641
385 0.138415
386 0.087022
387 0.120693
388 0.081018
389 0.104866
390 0.084394
391 0.116106
392 0.102701
393 0.125301
394 0.099904
395 0.126353
396 0.091723
397 0.100013
398 0.081762
399 0.136297
400 0.118378
};
\addlegendentry{Uniform in $\mathbb{R}^n$}
\end{axis}
\end{tikzpicture}
    \caption{Convergence metric evolution for the stress-constrained problem with uniform mesh in $\mathbb{R}^n$.}
    \label{fig:stresskktfunction2}
\end{figure}
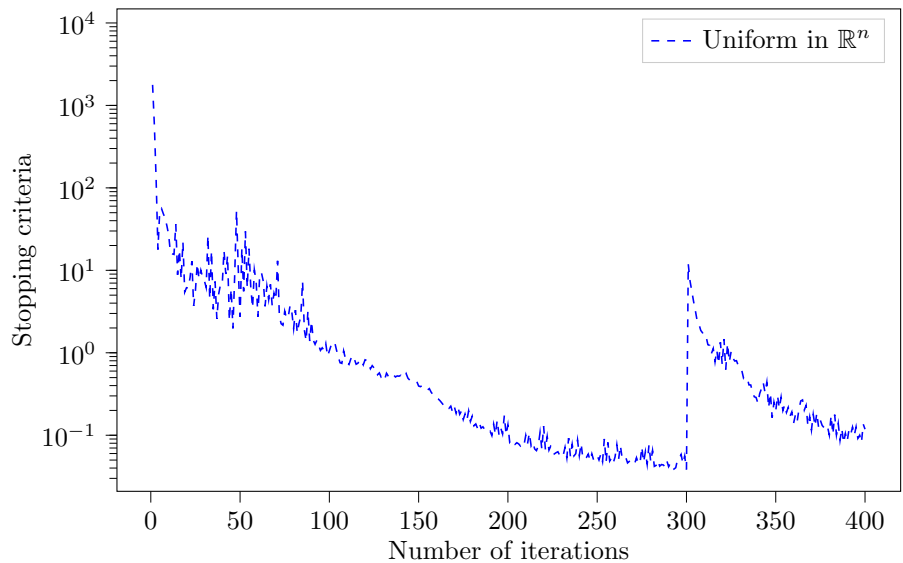

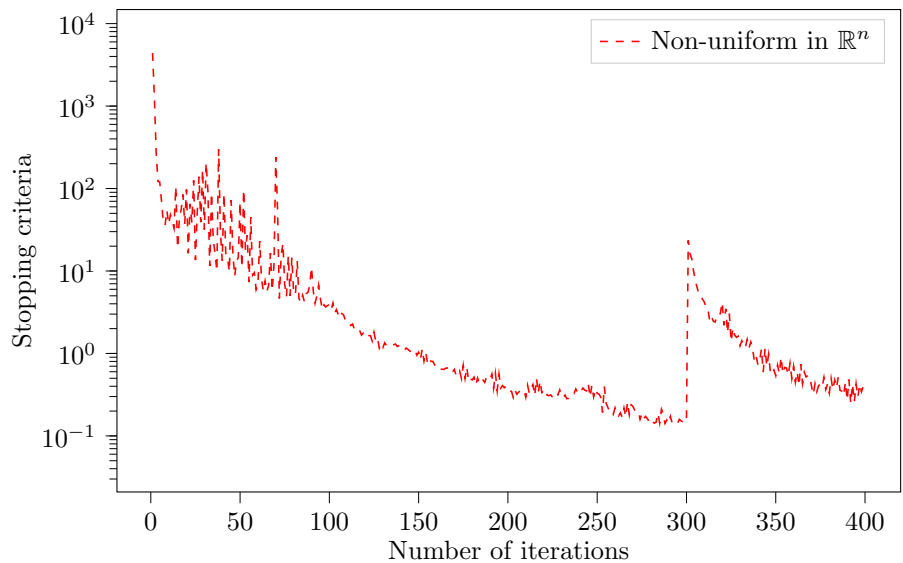
\begin{figure}
    \centering
\begin{tikzpicture}

\begin{axis}[
height=\figureheight,
legend cell align={left},
legend style={fill opacity=0.8, draw opacity=1, text opacity=1, draw=white!80!black},
log basis y={10},
tick align=outside,
tick pos=left,
width=\figurewidth,
x grid style={white!69.0196078431373!black},
xlabel={Number of iterations},
xmin=-18.95, xmax=419.95,
xtick style={color=black},
y grid style={white!69.0196078431373!black},
ylabel={Stopping criteria},
ymin=0.0208917309728883, ymax=14811.3521698619,
ymode=log,
ytick style={color=black}
]
\addplot [semithick, red, dashed]
table {%
1 4402.888506
2 1288.495009
3 275.140316
4 121.889549
5 120.080525
6 70.437532
7 39.016514
8 35.954867
9 54.304186
10 37.089141
11 47.816466
12 40.350136
13 33.01208
14 105.335566
15 18.946988
16 43.224892
17 55.176552
18 84.620102
19 37.368509
20 97.349244
21 16.423713
22 65.227054
23 38.734307
24 125.767774
25 13.567205
26 49.526326
27 138.97233
28 39.112065
29 168.892707
30 31.843184
31 205.96745
32 91.029905
33 11.495063
34 88.712549
35 23.588169
36 16.248152
37 9.921913
38 298.507722
39 29.236121
40 13.384741
41 86.843485
42 15.432836
43 14.255535
44 9.845692
45 72.166548
46 15.591429
47 8.619982
48 13.138893
49 14.607286
50 68.53739
51 11.005095
52 96.651434
53 18.098003
54 27.03617
55 7.320066
56 44.789134
57 8.878454
58 9.517755
59 5.943649
60 6.826232
61 23.003648
62 8.10363
63 5.687721
64 7.554254
65 5.648754
66 6.693429
67 16.474785
68 5.761396
69 11.896602
70 240.926314
71 53.650331
72 4.640473
73 15.435612
74 21.682828
75 6.806004
76 4.890746
77 15.258763
78 5.000723
79 14.595537
80 5.514421
81 5.163081
82 13.139984
83 4.613612
84 4.152104
85 5.075311
86 4.256304
87 5.281518
88 5.503222
89 6.677724
90 11.123719
91 5.722093
92 4.070486
93 4.986529
94 6.711483
95 4.60284
96 3.410191
97 3.882701
98 3.658399
99 3.813323
100 3.915201
101 3.381838
102 4.072312
103 3.092715
104 3.402996
105 2.944257
106 3.33436
107 3.040638
108 2.911994
109 2.512378
110 2.503359
111 2.285219
112 2.18022
113 2.274057
114 1.961602
115 1.994977
116 1.847618
117 1.879461
118 1.668896
119 1.76658
120 1.669998
121 1.775364
122 1.653404
123 1.64178
124 1.416598
125 1.800895
126 1.406878
127 1.52489
128 1.16311
129 1.314422
130 1.18174
131 1.345125
132 1.244743
133 1.251641
134 1.278232
135 1.232665
136 1.271913
137 1.301407
138 1.21431
139 1.242408
140 1.246573
141 1.14119
142 1.100272
143 1.160645
144 1.122341
145 1.071796
146 1.052686
147 1.116135
148 1.094374
149 0.950461
150 1.032494
151 0.838443
152 1.135906
153 0.810976
154 1.033667
155 0.753479
156 0.849764
157 0.799153
158 0.801188
159 0.72296
160 0.711944
161 0.638436
162 0.637432
163 0.645518
164 0.63858
165 0.659391
166 0.668694
167 0.632132
168 0.6537
169 0.596416
170 0.640449
171 0.538948
172 0.534243
173 0.60908
174 0.494864
175 0.715131
176 0.499136
177 0.682488
178 0.456147
179 0.490114
180 0.477743
181 0.515351
182 0.41194
183 0.500531
184 0.427551
185 0.471044
186 0.477103
187 0.445045
188 0.517749
189 0.493235
190 0.484854
191 0.560581
192 0.37939
193 0.570673
194 0.354638
195 0.560578
196 0.380867
197 0.407209
198 0.383672
199 0.402093
200 0.354557
201 0.34926
202 0.347819
203 0.292197
204 0.356276
205 0.319432
206 0.339413
207 0.349965
208 0.362303
209 0.36081
210 0.294842
211 0.411613
212 0.32769
213 0.37558
214 0.418649
215 0.307596
216 0.498618
217 0.329026
218 0.463976
219 0.334883
220 0.345064
221 0.364827
222 0.311243
223 0.308106
224 0.313779
225 0.294777
226 0.320544
227 0.334763
228 0.315292
229 0.360316
230 0.298724
231 0.293633
232 0.312516
233 0.282725
234 0.284549
235 0.317813
236 0.329203
237 0.325295
238 0.370516
239 0.33669
240 0.366671
241 0.35299
242 0.381111
243 0.358946
244 0.351462
245 0.432612
246 0.295218
247 0.410244
248 0.298799
249 0.373487
250 0.32938
251 0.320609
252 0.296332
253 0.192776
254 0.395217
255 0.249364
256 0.220312
257 0.199355
258 0.204259
259 0.182984
260 0.218813
261 0.186583
262 0.174785
263 0.190342
264 0.168543
265 0.21786
266 0.170516
267 0.25723
268 0.180815
269 0.238265
270 0.228153
271 0.205561
272 0.204222
273 0.202281
274 0.159092
275 0.185842
276 0.161195
277 0.1725
278 0.163418
279 0.150024
280 0.161162
281 0.147155
282 0.143081
283 0.146508
284 0.188888
285 0.143707
286 0.209075
287 0.181758
288 0.141531
289 0.156095
290 0.142675
291 0.171889
292 0.142338
293 0.137456
294 0.141819
295 0.14142
296 0.158805
297 0.152322
298 0.146871
299 0.142066
300 0.160928
301 23.490147
302 16.075371
303 13.949223
304 9.946824
305 7.550494
306 6.266581
307 5.199672
308 5.08812
309 4.531268
310 4.210297
311 3.542452
312 3.460883
313 2.508494
314 2.835869
315 2.467948
316 2.399751
317 2.779608
318 2.729819
319 3.083126
320 4.188392
321 2.196881
322 3.477037
323 1.777577
324 3.031446
325 1.39499
326 2.042543
327 1.669474
328 1.567653
329 1.621104
330 1.239495
331 1.422336
332 1.195369
333 1.185896
334 1.470671
335 1.161156
336 1.376025
337 1.147968
338 0.816009
339 0.94317
340 0.727442
341 0.823959
342 1.137688
343 0.593881
344 1.156151
345 0.69688
346 0.972191
347 0.591273
348 0.553062
349 0.627274
350 0.537831
351 0.719098
352 0.570923
353 0.903147
354 0.568273
355 0.677897
356 0.66699
357 0.544505
358 0.757454
359 0.496716
360 0.694487
361 0.450394
362 0.600675
363 0.406983
364 0.739889
365 0.43179
366 0.593128
367 0.519159
368 0.537723
369 0.533176
370 0.403397
371 0.318534
372 0.43018
373 0.344054
374 0.405896
375 0.424205
376 0.390011
377 0.509314
378 0.365748
379 0.426734
380 0.332818
381 0.498551
382 0.36271
383 0.472522
384 0.383296
385 0.515708
386 0.356712
387 0.349698
388 0.309512
389 0.449664
390 0.270304
391 0.481558
392 0.248087
393 0.426349
394 0.257968
395 0.426488
396 0.296372
397 0.423117
398 0.344731
399 0.396944
400 0.347145
};
\addlegendentry{Non-uniform in $\mathbb{R}^n$}
\end{axis}

\end{tikzpicture}
    \caption{Convergence metric evolution for the stress-constrained problem with non-uniform mesh in $\mathbb{R}^n$.}
    \label{fig:stresskktfunction3}
\end{figure}

\begin{figure}
    \centering
\begin{tikzpicture}

\begin{axis}[
height=\figureheight,
legend cell align={left},
legend style={fill opacity=0.8, draw opacity=1, text opacity=1, draw=white!80!black},
log basis y={10},
tick align=outside,
tick pos=left,
width=\figurewidth,
x grid style={white!69.0196078431373!black},
xlabel={Number of iterations},
xmin=-18.95, xmax=419.95,
xtick style={color=black},
y grid style={white!69.0196078431373!black},
ylabel={Stopping criteria},
ymin=0.0208917309728883, ymax=14811.3521698619,
ymode=log,
ytick style={color=black}
]
\addplot [semithick, red]
table {%
1 8016.807598
2 2615.009315
3 505.916728
4 80.487268
5 75.012019
6 71.249204
7 57.597568
8 50.61801
9 46.485492
10 37.630333
11 29.235547
12 99.780155
13 24.699122
14 70.851082
15 29.740275
16 93.49264
17 23.720357
18 149.758393
19 20.706552
20 82.335893
21 20.416008
22 52.769662
23 22.901529
24 58.743389
25 27.888172
26 57.640688
27 13.969021
28 16.856779
29 123.353807
30 23.11576
31 73.469276
32 16.801851
33 12.156007
34 13.770374
35 60.492022
36 23.371364
37 76.220406
38 19.12478
39 23.366283
40 37.197823
41 9.278165
42 13.419246
43 235.538871
44 60.034443
45 18.214473
46 137.884429
47 23.339691
48 9.730341
49 12.815562
50 12.18101
51 15.428415
52 244.841364
53 42.078388
54 10.722665
55 10.734741
56 21.370358
57 12.675519
58 71.840878
59 13.726128
60 29.450154
61 16.602591
62 333.451737
63 60.697348
64 67.201932
65 10.614176
66 11.25253
67 7.544703
68 19.230126
69 6.542618
70 15.744851
71 6.625239
72 111.093294
73 17.633755
74 9.468029
75 13.051381
76 6.951802
77 18.721374
78 10.842771
79 8.881461
80 6.242663
81 7.934343
82 8.355484
83 5.625644
84 4.993095
85 6.802699
86 4.537413
87 5.619994
88 4.783125
89 5.554567
90 4.634251
91 5.615949
92 4.766826
93 4.681483
94 4.203182
95 3.941893
96 3.983324
97 4.214374
98 3.871604
99 4.407729
100 3.315614
101 4.377607
102 3.207188
103 3.429711
104 3.271264
105 3.229841
106 3.273515
107 3.340479
108 3.310338
109 2.929831
110 3.563734
111 2.906065
112 3.953634
113 3.04173
114 3.060425
115 2.827017
116 2.801916
117 2.959527
118 2.634403
119 2.992815
120 2.552196
121 2.603058
122 2.802592
123 2.548021
124 2.666547
125 2.505912
126 2.624702
127 2.519531
128 2.247439
129 2.51415
130 2.198571
131 2.675771
132 2.028172
133 2.067719
134 1.9378
135 2.048903
136 1.949314
137 1.907607
138 2.252309
139 1.784598
140 2.039137
141 1.679151
142 1.788969
143 1.823842
144 1.762443
145 1.799863
146 1.780189
147 1.750861
148 1.749395
149 1.696798
150 1.774274
151 1.6644
152 1.63929
153 1.701396
154 1.538647
155 1.428092
156 1.353207
157 1.552342
158 1.277065
159 1.370067
160 1.065604
161 1.077622
162 0.972466
163 1.013253
164 0.911269
165 0.926739
166 0.813415
167 0.976184
168 0.783583
169 0.911158
170 0.784089
171 0.776685
172 0.884611
173 0.798742
174 0.793292
175 0.700341
176 0.689834
177 0.634128
178 0.667778
179 0.617087
180 0.752005
181 0.631127
182 0.74277
183 0.592966
184 0.598926
185 0.577108
186 0.588955
187 0.584802
188 0.597316
189 0.568197
190 0.581056
191 0.546584
192 0.603955
193 0.550277
194 0.595999
195 0.494109
196 0.520908
197 0.427714
198 0.507232
199 0.4414
200 0.456566
201 0.424039
202 0.465967
203 0.454098
204 0.437409
205 0.415573
206 0.386056
207 0.458411
208 0.430205
209 0.441137
210 0.528659
211 0.48765
212 0.53833
213 0.466761
214 0.496085
215 0.480615
216 0.471264
217 0.481294
218 0.488565
219 0.490996
220 0.577822
221 0.476421
222 0.563965
223 0.491729
224 0.539779
225 0.492601
226 0.514019
227 0.50918
228 0.49555
229 0.497787
230 0.47134
231 0.481395
232 0.431465
233 0.510263
234 0.403664
235 0.51383
236 0.365267
237 0.448513
238 0.348662
239 0.383053
240 0.34061
241 0.370094
242 0.323061
243 0.360324
244 0.351694
245 0.30041
246 0.356164
247 0.299683
248 0.329374
249 0.341309
250 0.324421
251 0.490672
252 0.310683
253 0.421273
254 0.348369
255 0.28868
256 0.346404
257 0.278166
258 0.31299
259 0.271351
260 0.312033
261 0.263983
262 0.298251
263 0.270161
264 0.262136
265 0.312808
266 0.274071
267 0.340129
268 0.26682
269 0.287917
270 0.269609
271 0.245625
272 0.280454
273 0.252714
274 0.271093
275 0.255952
276 0.288036
277 0.255296
278 0.259859
279 0.26818
280 0.252481
281 0.255838
282 0.268077
283 0.2646
284 0.29306
285 0.248264
286 0.305134
287 0.271641
288 0.237419
289 0.264799
290 0.254395
291 0.237943
292 0.263879
293 0.238673
294 0.244967
295 0.233195
296 0.225937
297 0.258714
298 0.220415
299 0.285766
300 0.222247
301 40.761547
302 29.72011
303 23.422885
304 17.32391
305 13.128149
306 9.535653
307 8.623955
308 7.454838
309 6.996977
310 6.834468
311 6.306321
312 4.758279
313 5.836175
314 3.461616
315 4.712749
316 3.188757
317 4.084568
318 3.31622
319 4.184964
320 2.971075
321 5.220604
322 2.79503
323 4.77582
324 2.81343
325 3.863183
326 2.974602
327 3.067487
328 2.846144
329 2.362066
330 2.238389
331 1.745376
332 2.14329
333 1.593189
334 1.850944
335 1.278195
336 1.499356
337 1.151652
338 1.328922
339 1.006147
340 1.516765
341 0.959685
342 1.021279
343 0.921127
344 1.0849
345 1.014217
346 1.136391
347 0.994776
348 1.102471
349 0.874607
350 0.9698
351 0.82339
352 1.158316
353 0.874106
354 1.061409
355 0.839338
356 0.851828
357 0.763362
358 0.754774
359 0.747698
360 0.73668
361 0.771953
362 0.849805
363 0.764681
364 0.81299
365 0.671098
366 0.862072
367 0.623777
368 0.938497
369 0.560419
370 1.077161
371 0.550328
372 1.001965
373 0.724582
374 0.569936
375 0.779829
376 0.521307
377 0.685007
378 0.502578
379 0.738311
380 0.497924
381 0.683579
382 0.520274
383 0.582501
384 0.510844
385 0.582045
386 0.490848
387 0.519196
388 0.559571
389 0.613169
390 0.458005
391 0.605457
392 0.476216
393 0.53178
394 0.444739
395 0.541108
396 0.444038
397 0.513912
398 0.439183
399 0.569013
400 0.456988
};
\addlegendentry{Non-uniform in $L^2$}
\end{axis}
\end{tikzpicture}
    \caption{Convergence metric evolution for the stress-constrained problem with non-uniform mesh in $L^2$.}
    \label{fig:stresskktfunction4}
\end{figure}
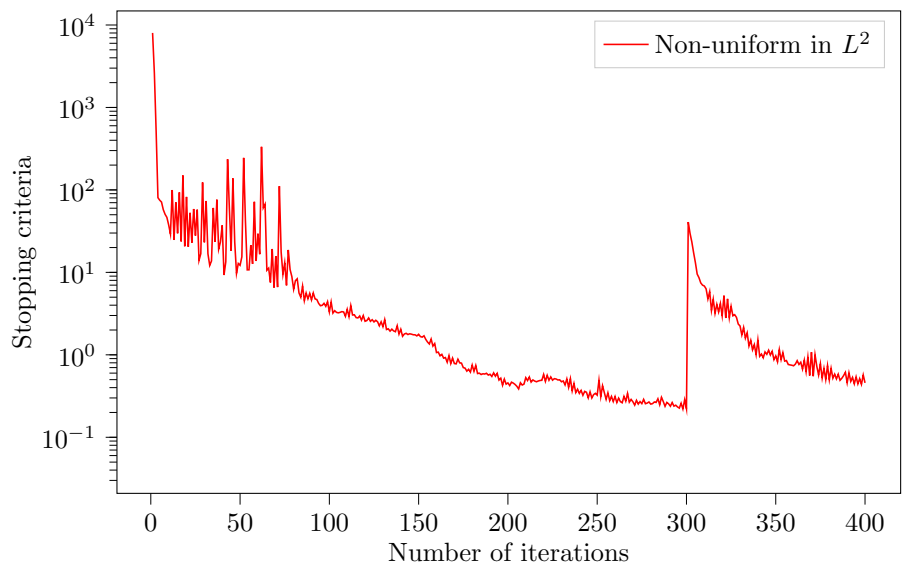

\begin{table}
    \centering
    \begin{tabular}{|M{0.5cm}|M{5.6cm}|M{5.6cm}|}
        \hline
                                         & Optimization in $\mathbb{R}^n$ & Optimization in $L^2$ \\ \hline
        \rotatebox{90}{Uniform mesh}     &
        \includegraphics[scale=0.1]{\MyPath/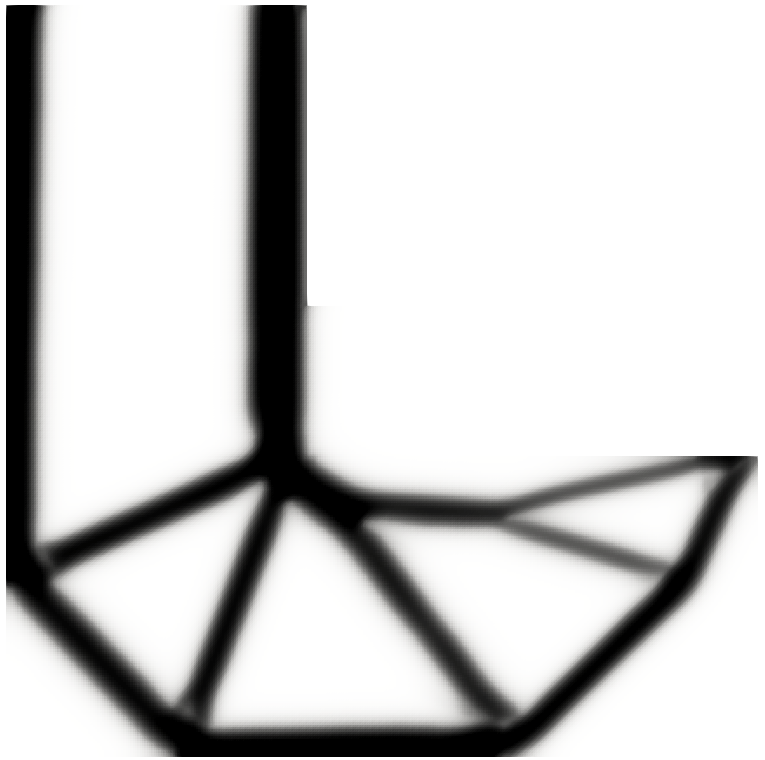}
                                         &
        \includegraphics[scale=0.1]{\MyPath/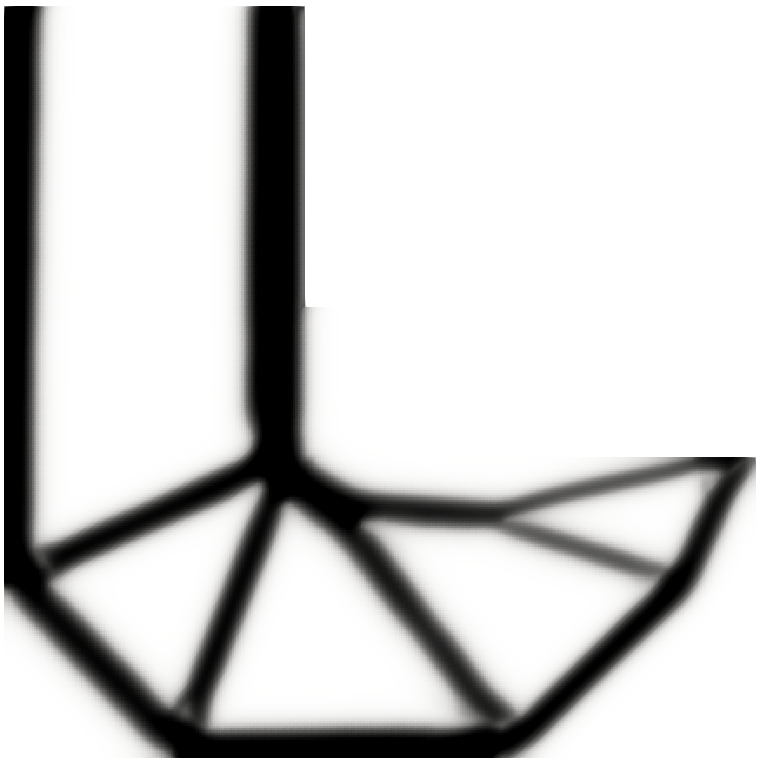}                        \\ \hline
        \rotatebox{90}{Non-uniform mesh} &
        \includegraphics[scale=0.1]{\MyPath/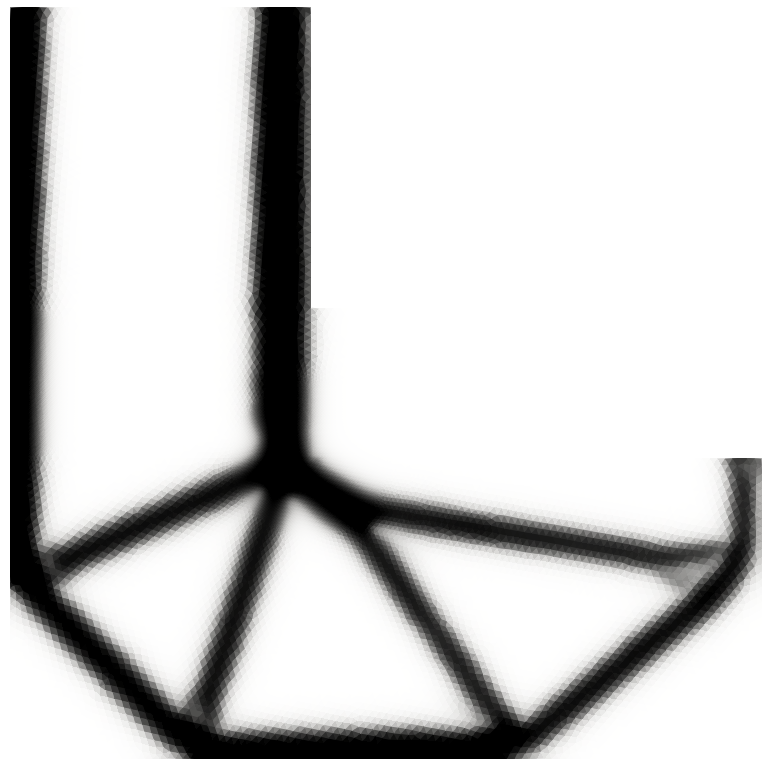}
                                         &
        \includegraphics[scale=0.1]{\MyPath/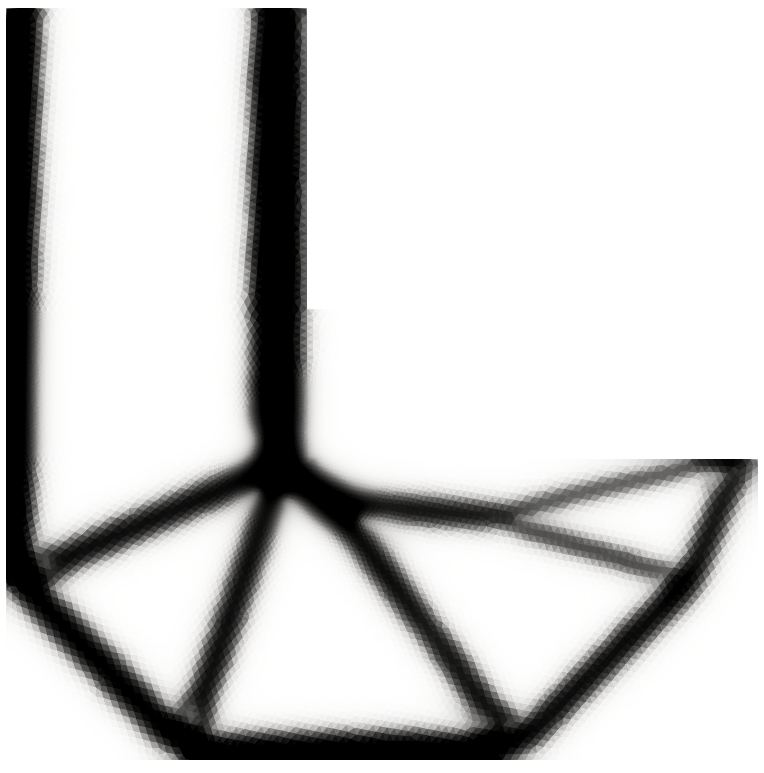}                    \\ \hline
    \end{tabular}
    \caption{Optimized designs for the stress-constrained problem.}
    \label{tab:stressresults}
\end{table}

\subsection{Uniform refinement}
\begin{figure}
    \centering
    \includegraphics[scale=0.5]{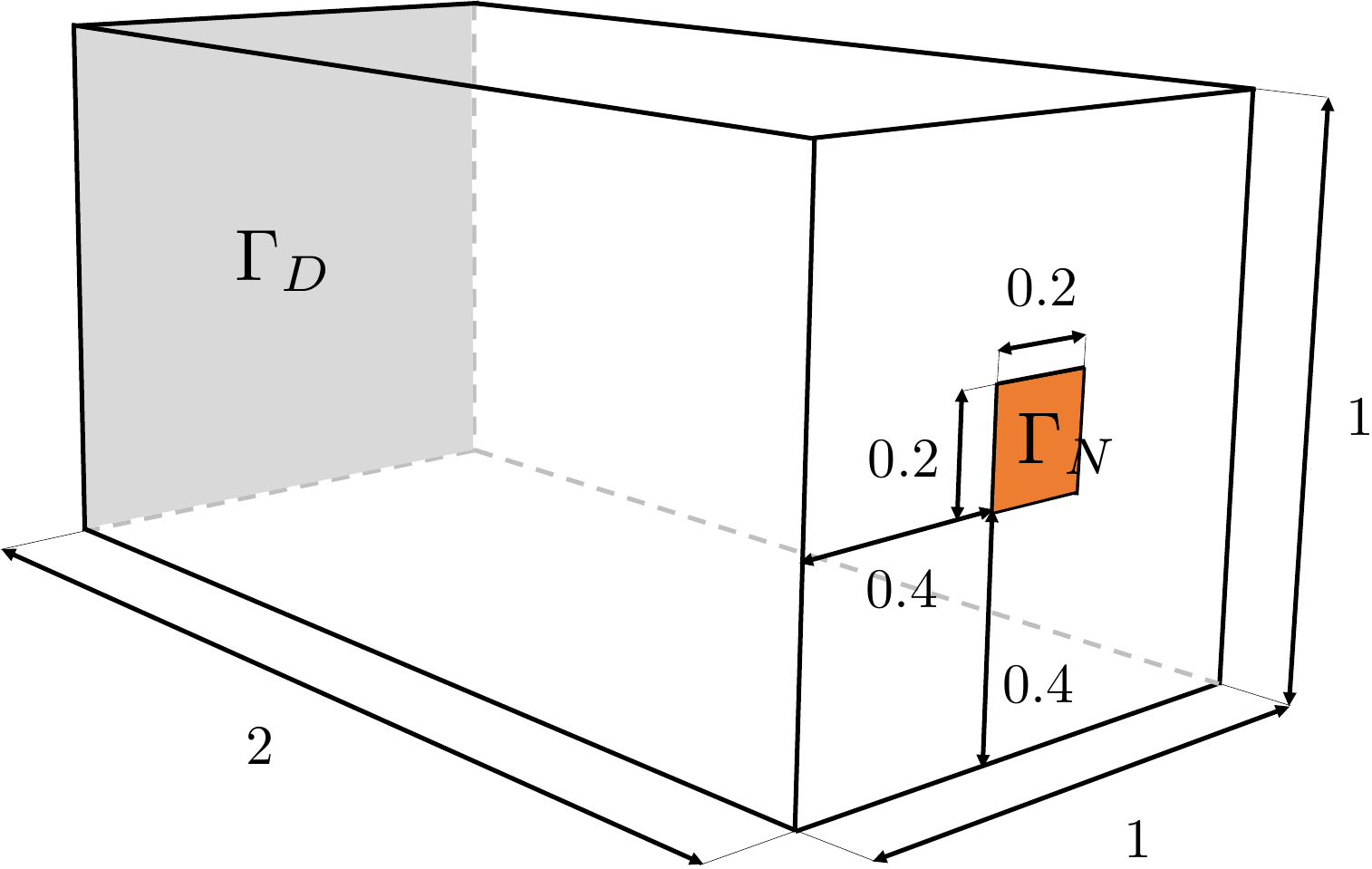}
    \caption{Design domain for the three dimensional compliance problem.}
    \label{fig:three_dim_domain}
\end{figure}

\begin{figure}
    \centering
    \includegraphics[scale=0.5]{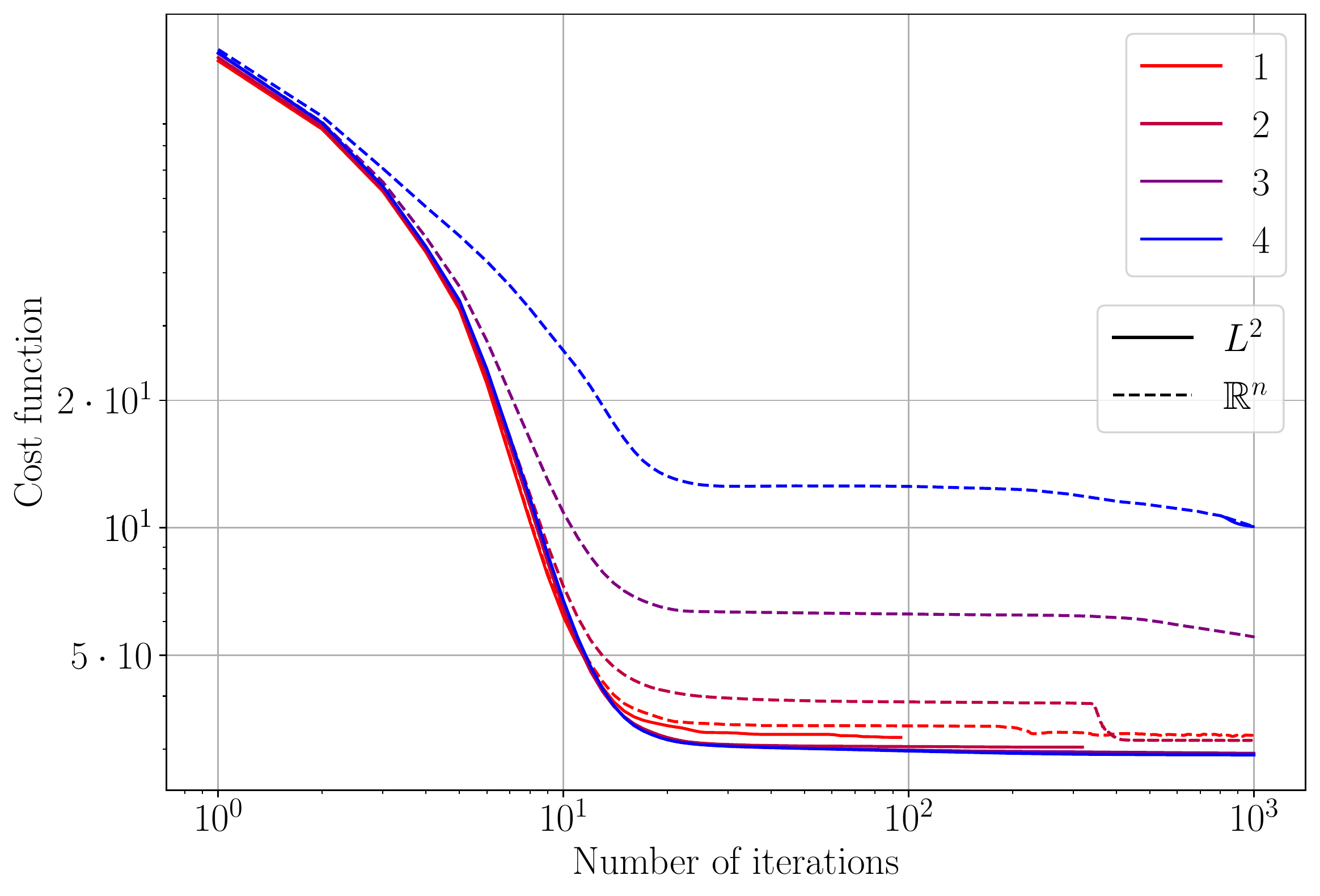}
    \caption{Cost function evolution of the compliance problem in three dimensions with 1, 2, 3 and 4 uniform refinements.
        Results plotted in log-log scale to highlight the differences.}
    \label{fig:three_dim_cost_evolution}
\end{figure}

\begin{figure}
    \centering
    \begin{subfigure}[b]{0.45\textwidth}
        \centering
        \includegraphics[scale=0.35]{\MyPath/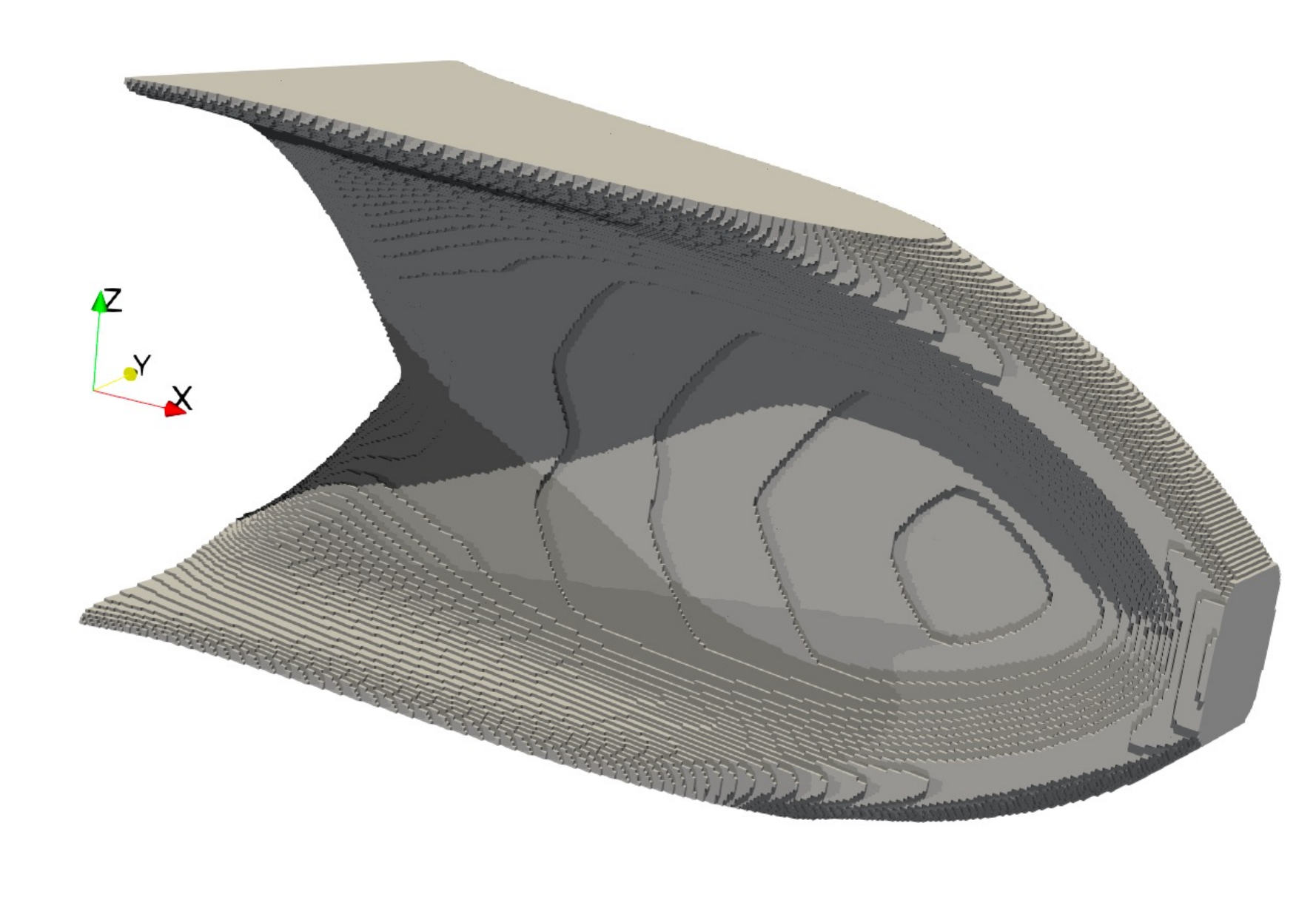}
        \caption{$L^2$}
    \end{subfigure}
    \hfill
    \begin{subfigure}[b]{0.45\textwidth}
        \centering
        \includegraphics[scale=0.35]{\MyPath/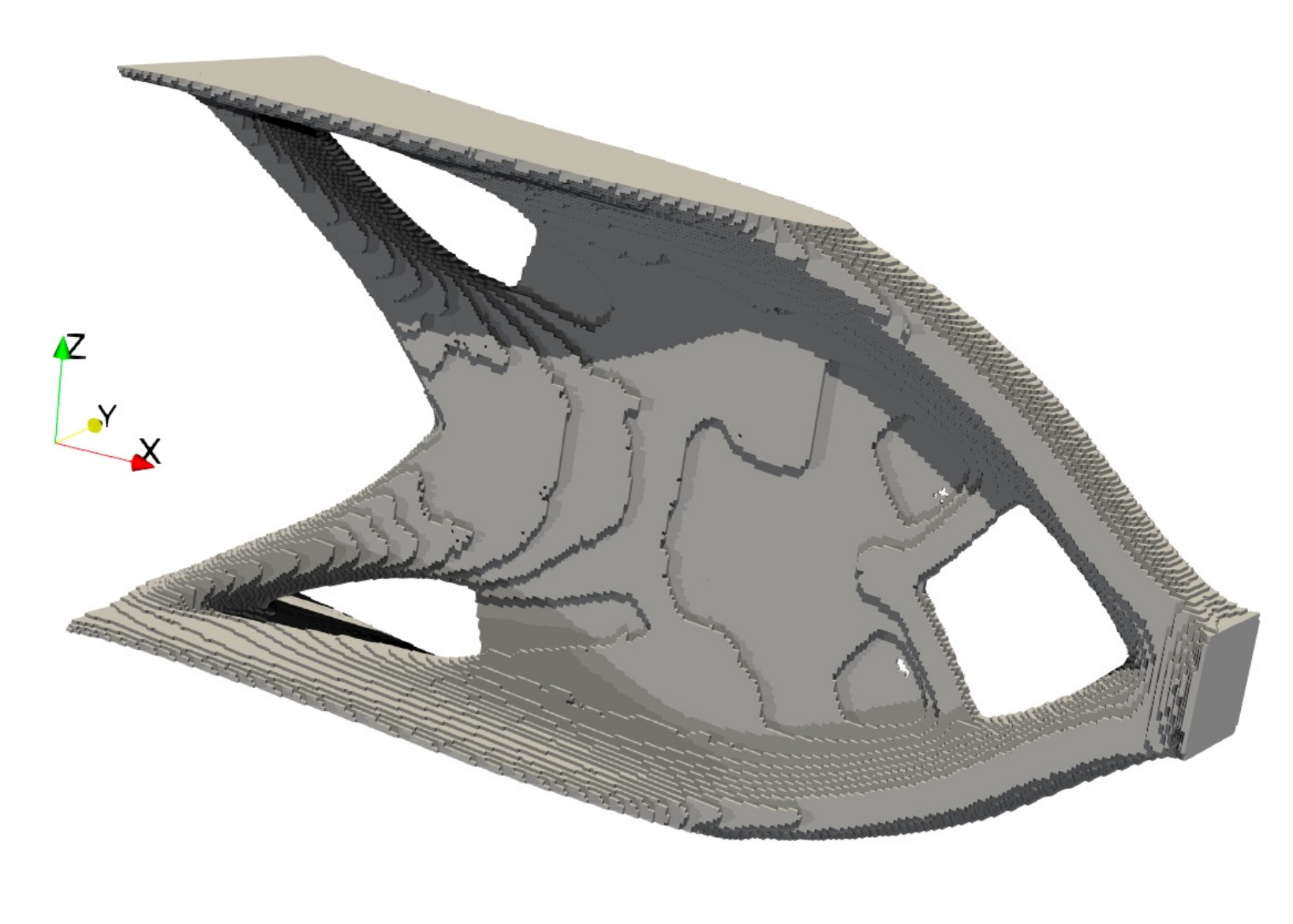}
        \caption{$\mathbb{R}^n$}
    \end{subfigure}
    \caption{Optimized designs for the compliance problem in three dimensions with four levels of refinement, thresholded with volume fraction $\hat\nu$ greater than 0.5}
    \label{fig:complianceresults3D}
\end{figure}

We return to the compliance problem for our next example, but in three dimensions, cf. Figure \ref{fig:three_dim_domain}.
The volume constraint is \hl{$\hat{V}= 0.15 |D|$},
the surface $\Gamma_N$ is subject to a traction $\mathbf{t} = -1.0 \mathbf{e}_2$, the length scale parameter $\kappa=5 \times 10^{-5}$ and the initial design is \hl{$\nu(\mathbf{x})=0.15$.}
\hl{To prevent the iterative solver from diverging, we set $\epsilon_{\nu}$ to $10^{-4}$.}
Instead of comparing two different meshes, we compare the influence of the mesh refinement 
starting from the same uniform mesh with 20x10x10 hexahedral elements.
\hl{As such, we perform an optimization with the initial mesh uniformly refined once, twice,
three and four times, for a total of 16,000, 128,000, 1,024,000 and 8,192,000 elements.}
The construction and refinement of the mesh is done by the utility meshing functions in Firedrake within the code shared along with this paper.
Each hexahedral element is split in eight at each refinement level.
We plot the evolution of the cost function \hl{for 1000 iterations} corresponding to 
the $L^2$ and $\mathbb{R}^n$ NLP algorithms for each refinement level in Figure \ref{fig:three_dim_cost_evolution}.
The behavior of the $\mathbb{R}^n$ NLP algorithm clearly depends on the refinement level, whereas the $L^2$ does not.
Notably, after 1000 iterations, the optimized design with four levels of refinement take compliances values of 10.06 and 2.91 for the $\mathbb{R}^n$ and $L^2$ NLP algorithms.
We compare the four level of refinement designs obtained with the $\mathbb{R}^n$ and the $L^2$ NLP algorithms in Figure \ref{fig:complianceresults3D}.
\hl{Upon inspection of the designs, we noticed that the $\mathbb{R}^n$ algorithm fails to reach the lower bound of the volume fraction (0);
it never dips below $\nu=10^{-2}$ for the void phase.
Due to the volume constraint and the ``heavier'' void phase, the $\mathbb{R}^n$ optimizer cannot add more
mass to the structure and therefore, the compliance is higher.
We conjecture that the reason behind this is due to the difference between the derivative and the gradient.
The $\mathbb{R}^n$ algorithm uses $D\theta$ whereas the $L^2$ uses $\nabla\theta = \mathbf{M}^{-1}D\theta$.
Thus, although $D\theta$ and $\nabla \theta$ are parallel since $\mathbf{M} = |\Omega_e| \mathbf{I}$, the sensitivity of the cost and constraint function of small elements is less influential
in the $\mathbb{R}^n$ vs $L^2$ algorithm.
We further conjecture that using an interpolation scheme with nonzero derivative values for $\nu=0$,
such as RAMP \mbox{\citep{ramp_penal}}, could alleviate this issue.}

\subsection{\hl{Adaptive mesh refinement}}
Our $L^2$ GCMMA implementation is specially suitable for use with AMR strategies during the optimization \citep{de2020three}.
In our next two examples, we show the algorithm utility when applying AMR during the optimization and compare it with the original $\mathbb{R}^n$ implementation.
There are three questions to address when applying AMR in a topology optimization problem:
where in the domain we apply AMR, what kind of AMR (coarsening, refinement or both)
and when to apply the AMR during the optimization.
We use an element-based error quantity to determine the regions subject to AMR.
The details are explained in \ref{sec:amr_appendix}.
We use two AMR strategies: only refinement or only coarsening.
We apply the AMR after a pre-determined number of iterations
during the optimization\footnote{We are not concerned with applying an optimal strategy for AMR during the optimization.
See \citep{ziems2011adaptive} for a more sophisticated scheme.}.
Table \ref{tab:ref_strategy} summarizes four AMR schemes wherein
it is seen that refinement occurs either during early iterations as the
design evolves or later once the design is well defined.
Coarsening only occurs late in the design evolution as well.
All problem/algorithm combinations are run for 600 iterations.

\begin{table}
    \centering
    \begin{tabular}{l|ll}
        \hline
             \diagbox{Strategy}{AMR type}& Coarsening & Refinement \\ \hline
        A     &
        100, 150
                                         &
        10, 80                        \\ 
        B &
        150, 200
                                         &
        100, 150                    \\ 
    \end{tabular}
    \caption{Four AMR strategies. Iteration number after which AMR is applied.}
    \label{tab:ref_strategy}
\end{table}

\begin{figure}
    \centering
    \includegraphics[scale=0.5]{\MyPath/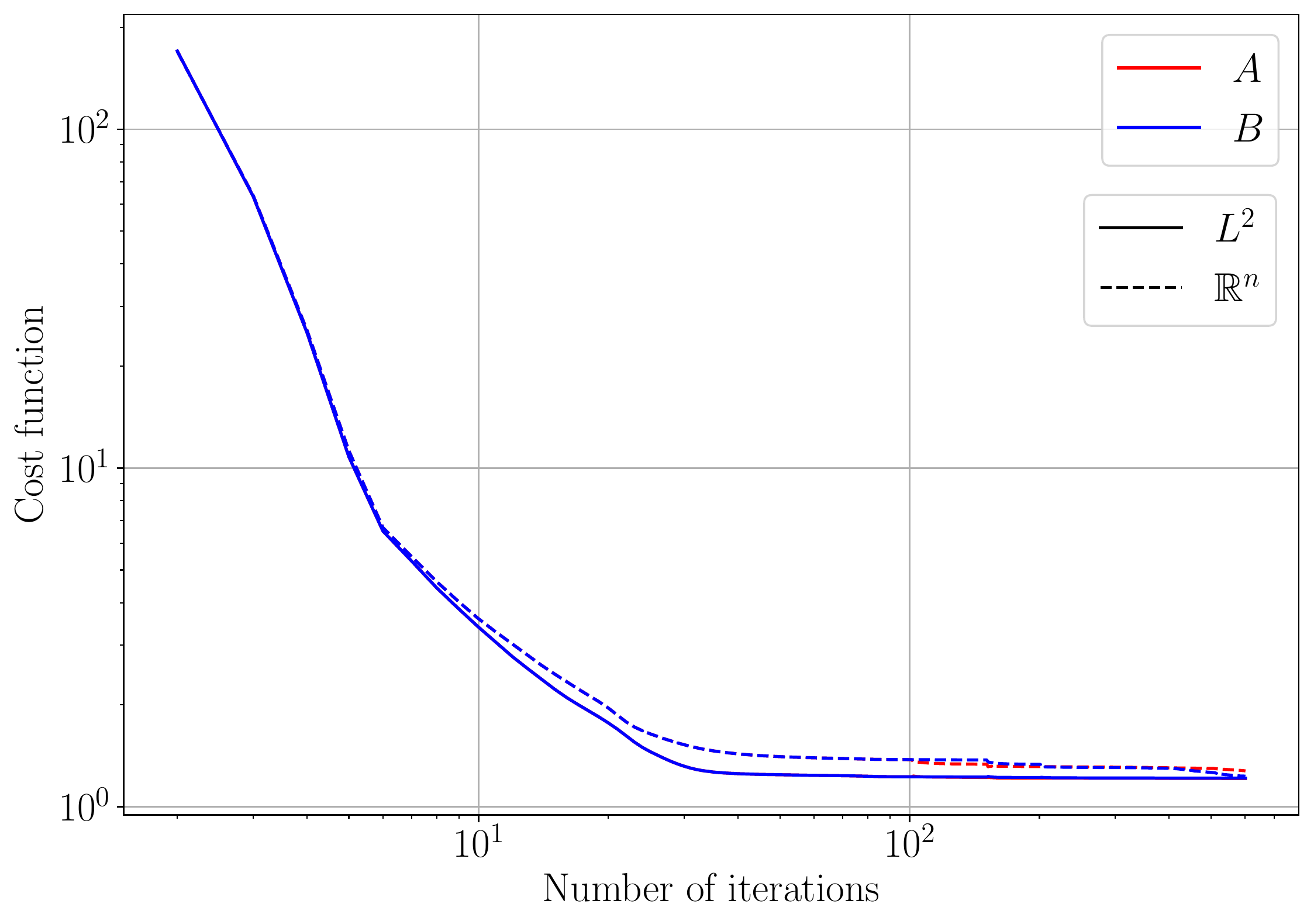}
    \caption{Cost function evolution for the compliance problem with AMR coarsening.}
    \label{fig:compliance_amr_coarsen}
\end{figure}

\begin{figure}
    \centering
    \includegraphics[scale=0.5]{\MyPath/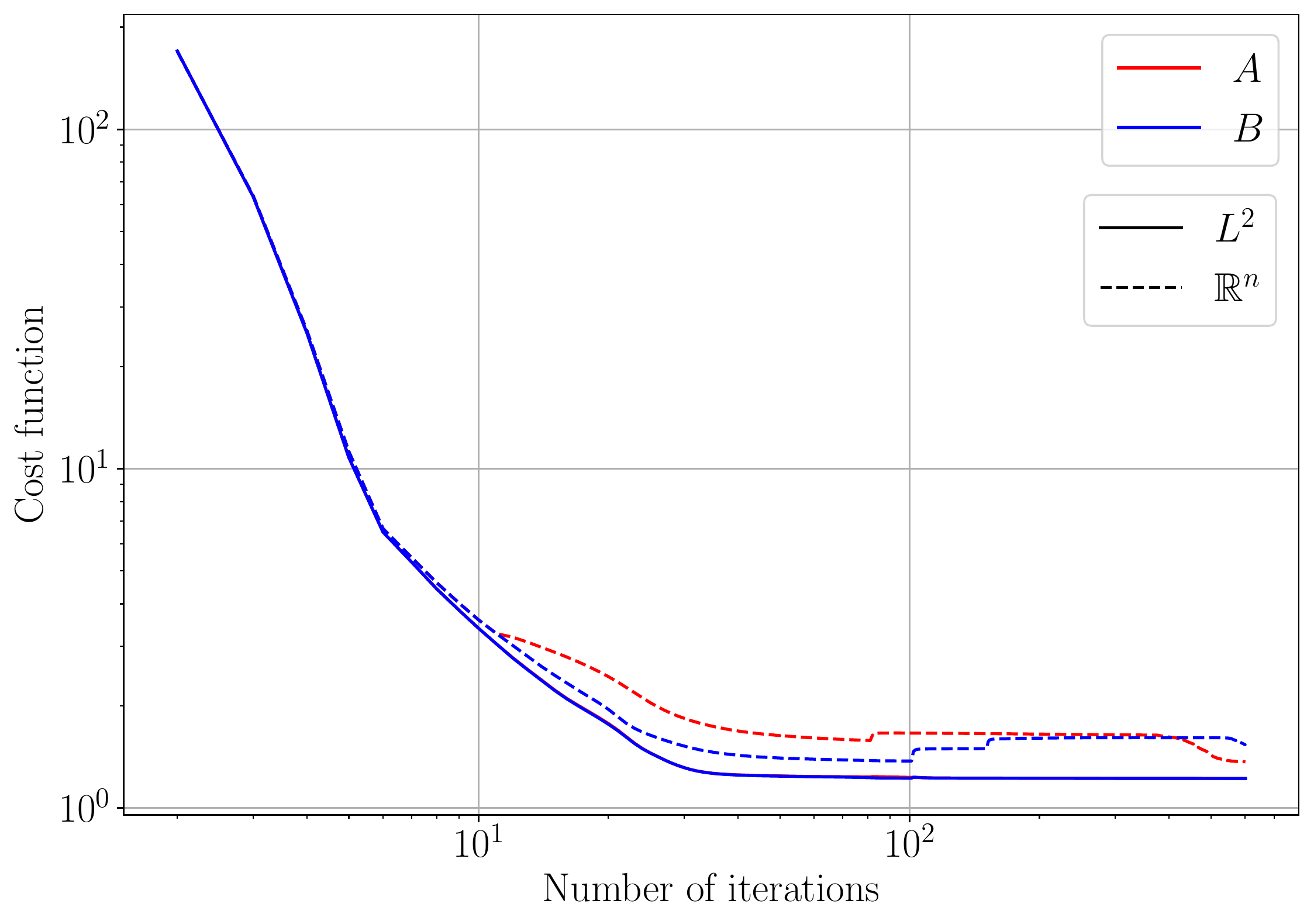}
    \caption{Cost function evolution for the compliance problem with AMR refinement.}
    \label{fig:compliance_amr_refine}
\end{figure}

Our first AMR example is a repeat of our cantilever example, 
but applying coarsening and refinement to the Table \ref{tab:ref_strategy} strategies. 
The optimized designs, shown in Tables \ref{fig:compliance_refinement}
and \ref{fig:compliance_coarsening} for refinement and coarsening respectively,
are different when using the $\mathbb{R}^n$ algorithm.
The cost function evolution in Figures \ref{fig:compliance_amr_coarsen}
and \ref{fig:compliance_amr_refine} further reflect their difference.
The meshes for the optimized designs illustrated in the Appendix, cf. Figures \ref{fig:compliance_amr_refinement_grid},
\ref{fig:compliance_amr_coarsening_grid}, highlight the mesh independence of the
the $L^2$ versus $\mathbb{R}^n$ algorithm.
The cost function evolution for coarsening, cf. Figure \ref{fig:compliance_amr_coarsen},
is not affected by the different refinement strategies even with the $\mathbb{R}^n$ optimization.
Most likely, this is because both $A$ and $B$ strategies are applied at similar iteration numbers,
after which the designs have almost converged.
It is important to highlight that the algorithm in $L^2$ still performs better.

\begin{table}
    \centering
    \begin{tabular}{|M{0.5cm}|M{5.6cm}|M{5.6cm}|}
        \hline
                                         & Optimization in $\mathbb{R}^n$ & Optimization in $L^2$ \\ \hline
        \rotatebox{90}{Strategy A}     &
        \includegraphics[scale=0.2]{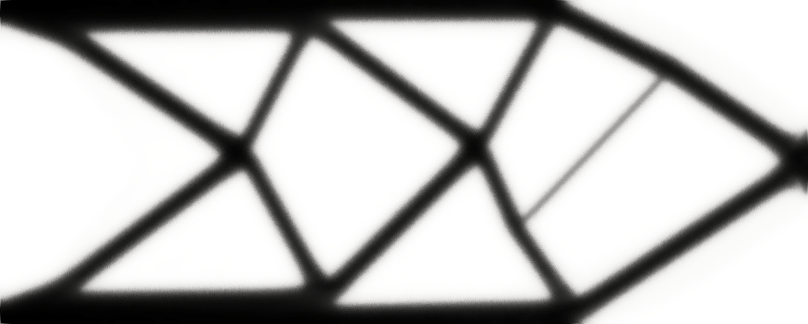}
                                         &
        \includegraphics[scale=0.2]{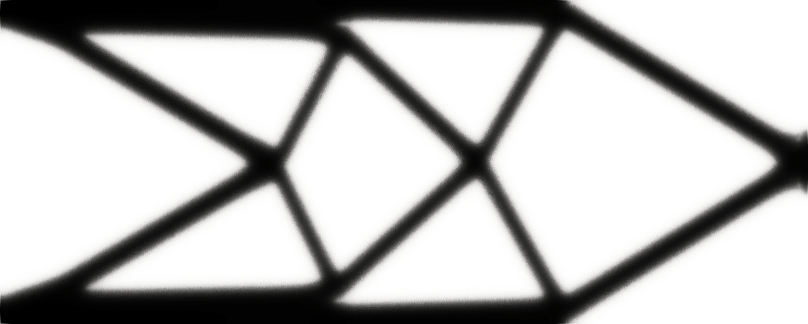}                        \\ \hline
        \rotatebox{90}{Strategy B} &
        \includegraphics[scale=0.2]{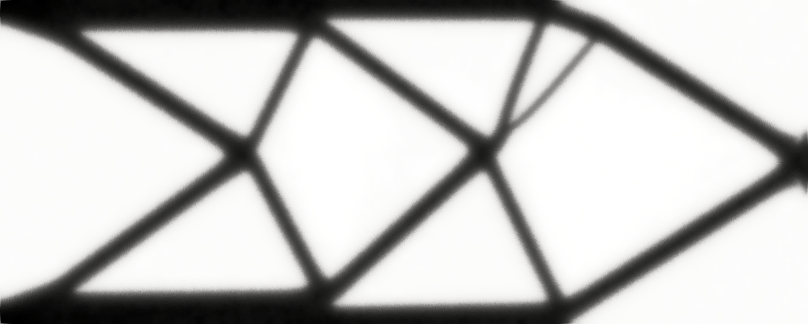}
                                         &
        \includegraphics[scale=0.2]{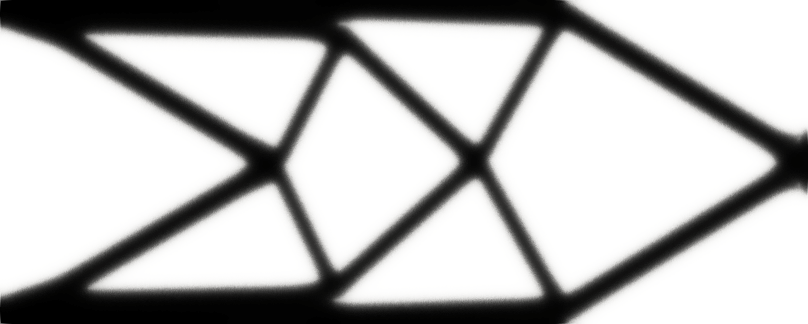}                    \\ \hline
    \end{tabular}
    \caption{Optimized designs for the compliance problem with AMR refinement only.}
    \label{fig:compliance_refinement}
\end{table}

\begin{table}
    \centering
    \begin{tabular}{|M{0.5cm}|M{5.6cm}|M{5.6cm}|}
        \hline
                                         & Optimization in $\mathbb{R}^n$ & Optimization in $L^2$ \\ \hline
        \rotatebox{90}{Strategy A}     &
        \includegraphics[scale=0.2]{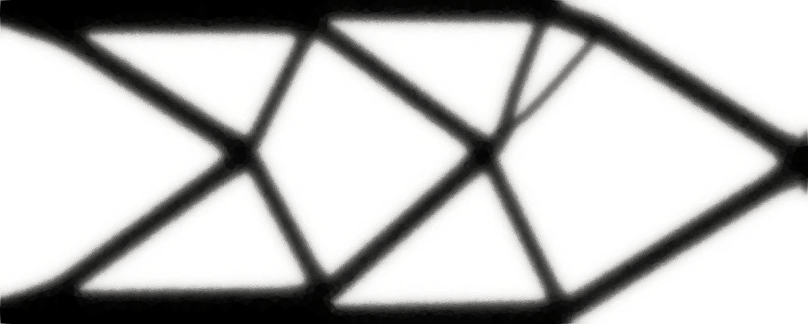}
                                         &
        \includegraphics[scale=0.2]{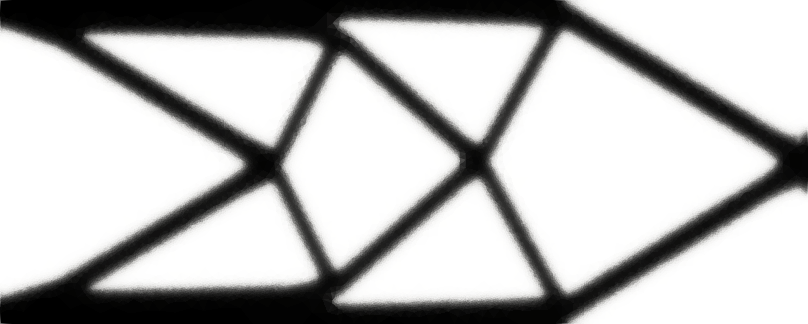}                        \\ \hline
        \rotatebox{90}{Strategy B} &
        \includegraphics[scale=0.2]{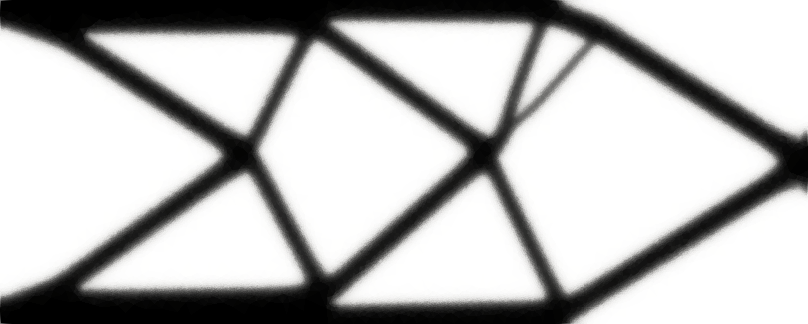}
                                         &
        \includegraphics[scale=0.2]{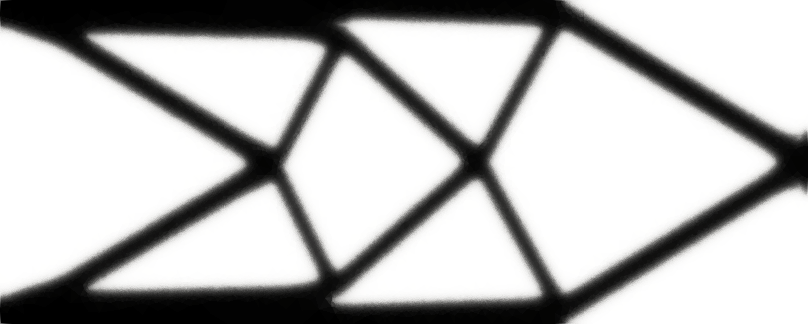}                    \\ \hline
    \end{tabular}
    \caption{Optimized designs for the compliance problem with AMR coarsening only.}
    \label{fig:compliance_coarsening}
\end{table}

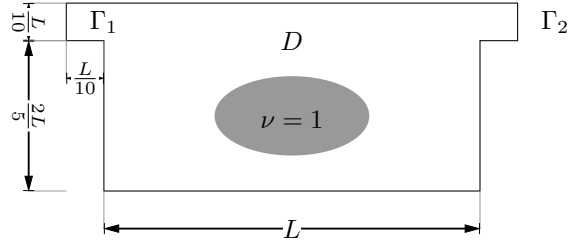
\begin{figure}[!h]
    \center
    \tikzset{>=latex}
    \begin{tikzpicture}[ every node/.style={inner sep=0pt}, spring/.style = {decorate,decoration={zigzag,amplitude=6pt,segment length=4pt}}, scale=5.0 ]
        \node (A) at (-0.1, 0.0) {};
        \node (B) at (-0.1, -0.1) {};
        \node (C) at (0, -0.1) {};
        \node (D) at (0.0, -0.5) {};
        \node (E) at (1.0, -0.5) {};
        \node (F) at (1.0, -0.1) {};
        \node (G) at (1.1, -0.1) {};
        \node (H) at (1.1, 0.0) {};
        \node (J) at (0.0, 0.0) {};
        \draw ($(A)$) -- ($(B)$) -- ($(C)$) -- ($(D)$) -- ($(E)$) -- ($(F)$) -- ($(G)$) -- ($(H)$) -- ($(A)$);
        \dimline[line style = {line width=0.7},
            extension start length=-0.1cm,
            extension end length=-0.1cm]{($(A) + (-0.1, 0)$)}{($(B) + (-0.1, 0.0)$)}{$\frac{L}{10}$};
        \dimline[line style = {line width=0.7},
            extension start length=0.1cm,
            extension end length=0.1cm]{($(D) + (-0.2, 0)$)}{($(B) + (-0.1, 0.0)$)}{$\frac{2L}{5}$};
        \dimline[line style = {line width=0.7},
            extension start length=0.1cm,
            extension end length=0.1cm]{($(E) + (0.0, -0.1)$)}{($(D) + (0.0, -0.1)$)}{$L$};
        \dimline[line style = {line width=0.7},
            extension start length=-0.1cm,
            extension end length=-0.1cm]{($(A) + (0.0, -0.2)$)}{($(J) + (0.0, -0.2)$)}{$\frac{L}{10}$};
        \node (GammaD) at ($(A) !.5! (B) + (0.1, 0)$) {$\Gamma_1$};
        \node (GammaF) at ($(G) !.5! (H) + (0.1, 0)$) {$\Gamma_2$};
        \node (Domain) at (0.5, -0.1) {$D$};
        \draw[ultra thick, draw=black, fill=black, opacity=0.4] (0.5, -0.3) ellipse (0.2cm and 0.1cm);
        \node (Domain) at (0.5, -0.31) {$\nu=1$};
    \end{tikzpicture}
    \caption{Thermal flow domain. $L$=1.}
    \label{fig:thermal_flow_domain}
\end{figure}

Next we solve a coupled thermal flow problem similar than in \cite{thermal_flow}
to demonstrate the algorithm's application to more complex physical phenomena.
The domain in Figure \ref{fig:thermal_flow_domain} is the cross section of a heat exchanger, where the design variable $\nu$ represents the volume fraction of a heat generating solid material that needs to be distributed to control the temperature and the fluid flow.
The goal of the optimization is to maximize the heat generated in the solid while keeping the pressure drop in the flow from the inlet $\Gamma_1$ to the outlet $\Gamma_2$ lower than a fixed value $P_{\text{drop}}$.
As such, the optimization problem reads

\newcommand{\ww}{\mathbf{w}}
\newcommand{\vv}{\mathbf{v}}
\newcommand{\VV}[1]{\mathbf{V}_{#1}}
\begin{subequations}\label{eq:thermal_flow_problem}
\begin{align}
    \underset{\nu \in V}{\text{max}}~ J(\nu)          & = \int_{D} \nu B (1 - T)~dV \,, \\
    \text{s.t.} ~\left(\ww, p, T \right) \in \VV{} \times Q \times W ~\text{satisfy} ~      \\
    F(\nu, \ww, p; \vv, q)  & = 0 \label{eq:flow_eq}\,, \\
    a_T(\ww; T, m) + c_T(\ww; T, m)                                & = 0     \label{eq:thermal_eq}                                                    \,,     \\
    \text{for all} ~\left(\vv, q, m \right) \in \VV{0} \times Q  \times W_0                \\
    G(\nu)                                             & = \int_{\Gamma_{1}} p ~dA - \int_{\Gamma_{2}} p ~dA \leq P_{\text{drop}}      \,,  
\end{align}
\end{subequations}
where $(\mathbf{w}, p, T)$ are the flow velocity, pressure and temperature and $(\mathbf{v}, q, m)$ are their admissible counterparts.
Equation \eqref{eq:flow_eq} is the weak form of the Navier-Stokes equation wherein
\begin{equation}
    \begin{aligned}
    F(\nu, \ww, p; \vv, q) & = \int_{D}\left( \left(\ww \cdot {\nabla} \ww \right) \cdot \vv + \frac{1}{Re} {\nabla} \ww	: {\nabla} \vv
                            + \frac{1}{Da} r(\nu) \ww \cdot \vv \right) ~dV\\
                            &+ \int_D \left( p \nabla \cdot \vv + q \nabla \cdot \ww \right)~dV                                                    \,,      \\
    \end{aligned}
    \label{eq:stokes_expanded}
\end{equation}
with Reynolds number $Re=1.0$ and Darcy number $Da=10^{-6}$.
The RAMP function \citep{ramp_penal} with $q_{\text{RAMP}} = 20.0$, i.e.
\begin{align}
    r(\nu) = \frac{\nu}{1+q_{\text{RAMP}}(1-\nu)}
\end{align}
is used to obtain discrete 0-1 designs.
At the inlet $\Gamma_1$, the Dirichlet condition $\ww = \ww_1$
is a horizontal parabolic profile with a maximum non-dimensional velocity $W_{\text{max}}=1$.
At the outlet $\Gamma_2$, a traction-free condition is imposed.
The remaining boundary $\Gamma \setminus(\Gamma_1 \cup \Gamma_2)$ has no-slip condition $\ww = 0$.

Equation \eqref{eq:thermal_eq} is the weak form of the advection-diffusion equation wherein
\newcommand{\nn}{\mathbf{n}}
\begin{equation}
    \begin{aligned}
        a_T(\nu; T, m) & = \int_{D} m \ww \cdot {\nabla} T + \frac{1}{Pe} {\nabla} T	: {\nabla} m  - \nu B (1 - T) m ~dV \,, \\
        \label{eq:temp_gls}
    \end{aligned}
\end{equation}
We assume the same thermal conductivities in the solid and and the fluid,
a Peclet number $Pe=10^4$.
The heat source in the solid $\nu B (1 - T)$ is proportional to the difference between a reference temperature 1 and the local temperature $T$, and the non-dimensional heat generation coefficient $B=0.01$ (more details on the heat source are in \cite{thermal_flow}).
The Galerkin Least Squares (GLS) stabilization term
\begin{equation}
    \begin{aligned}
        c_T(\ww; T, m)                & = \int_D \tau_{GLS}\mathcal{L}_T(T) \cdot \mathcal{L}_T(m) ~dV
    \end{aligned}
\end{equation}
stabilizes the otherwise highly oscillatory boundary layers due to the fluid velocity field,
where
\begin{equation}
    \begin{aligned}
        \tau_{GLS} = \beta_{GLS} \left( \frac{4 \ww \cdot \ww}{h^2} + \left(9 \frac{4 }{h^2Pe} \right)^2 \right)^{-0.5} \,,
    \end{aligned}
\end{equation}
$h$ is the element cell size and
\begin{align}
    \mathcal{L}_T(T) = \ww \cdot {\nabla} T + \frac{1}{Pe} \nabla T  - \nu B (1 - T)\,.
\end{align}
is the residual.
We use $\beta_{GLS} = 0.9$ in our examples.
We apply a constant temperature $T=0$ on $\Gamma_1$ and adiabatic boundary conditions over all surfaces with the exception of $\Gamma_1$ and $\Gamma_2$.
The finite element discretization of $(\mathbf{w}, p, T)$ uses linear Lagrange elements.
The cost function $J(\nu)$ aims to maximize the heat generation in the design domain
and the constraint $G(\nu) < P_{\text{drop}}=70.0$ limits the pressure drop in the system
and serves to regularize the problem as it imposes an upper bound on the fluid-solid
interface where $\nabla \nu \neq 0$.

\begin{figure}
    \centering
    \includegraphics[scale=0.5]{\MyPath/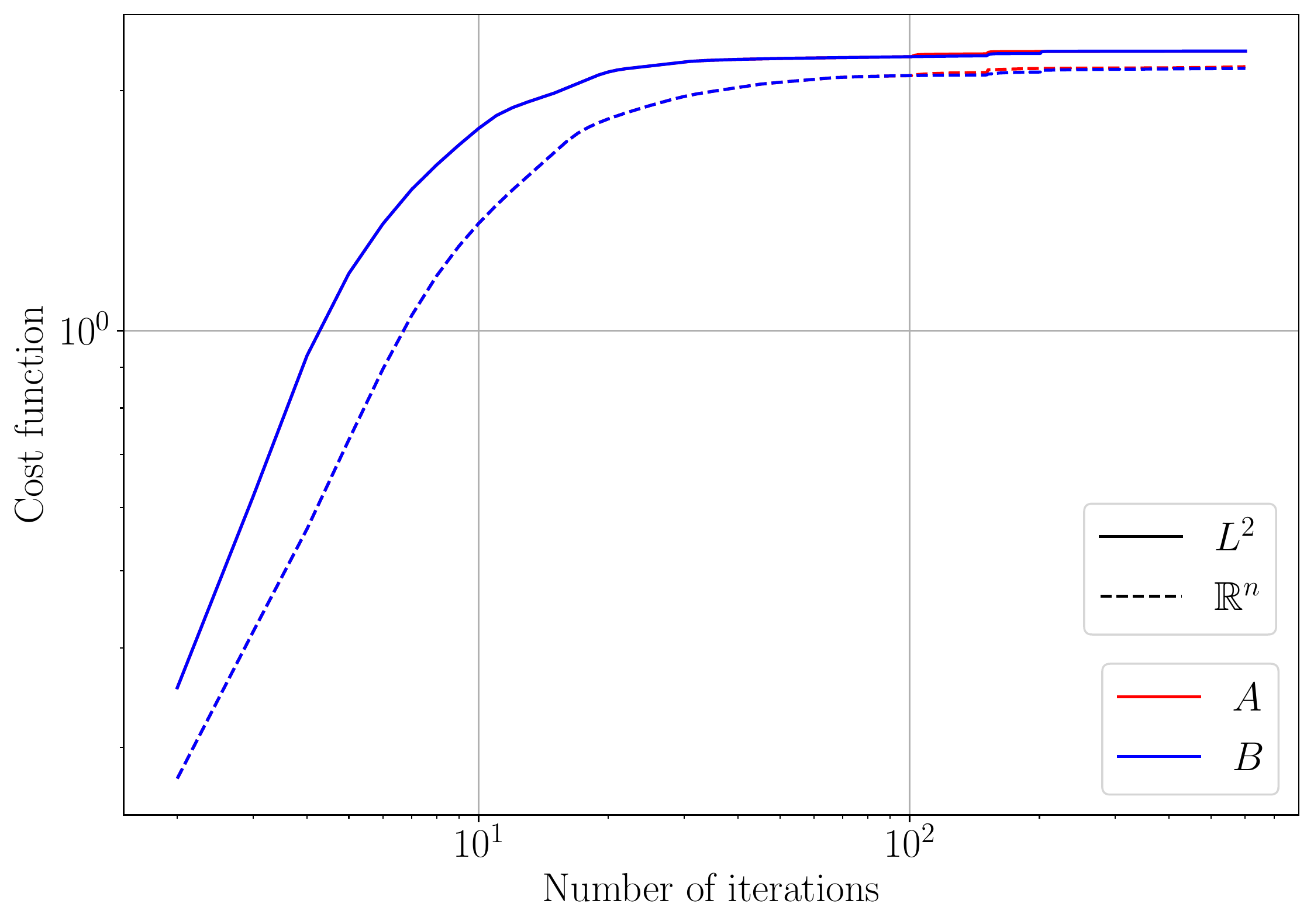}
    \caption{Cost function evolution for the thermal flow problem with AMR coarsening.}
    \label{fig:thermal_amr_refine}
\end{figure}

\begin{figure}
    \centering
    \includegraphics[scale=0.5]{\MyPath/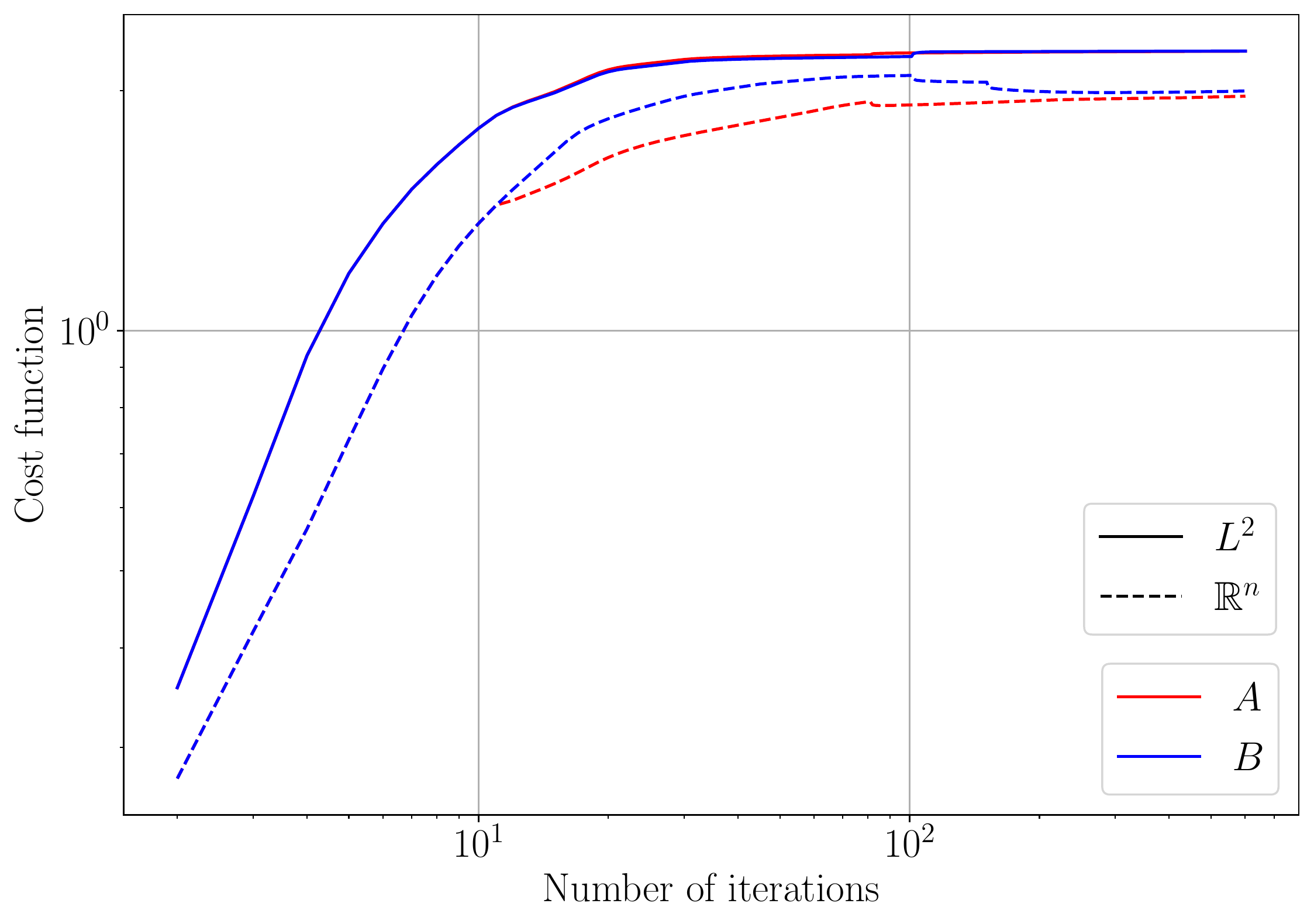}
    \caption{Cost function evolution for the thermal flow problem with AMR refinement.}
    \label{fig:thermal_amr_coarsen}
\end{figure}

Optimized designs are shown in Tables \ref{tab:thermal_flow_refinement} and \ref{tab:thermal_flow_coarsening} for refinement and coarsening strategies.
Again, the designs obtained using the $\mathbb{R}^n$ algorithm differ.
The cost function evolution in Figures \ref{fig:thermal_amr_refine} and \ref{fig:thermal_amr_coarsen} reflect this dependency as well.
The meshes for the optimized designs illustrated in the Appendix, cf. Figures \ref{fig:thermal_flow_amr_refinement_grid} and \ref{fig:thermal_flow_amr_coarsening_grid} highlight the mesh independence of the $L^2$ algorithm.

\begin{table}
    \centering
    \begin{tabular}{|M{0.5cm}|M{5.6cm}|M{5.6cm}|}
        \hline
                                         & Optimization in $\mathbb{R}^n$ & Optimization in $L^2$ \\ \hline
        \rotatebox{90}{Strategy A}     &
        \includegraphics[scale=0.2]{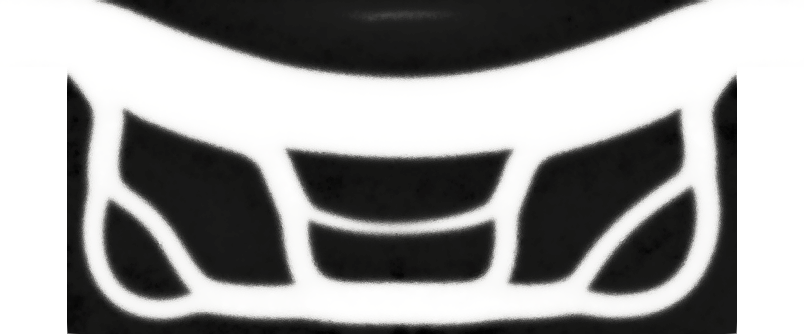}
                                         &
        \includegraphics[scale=0.2]{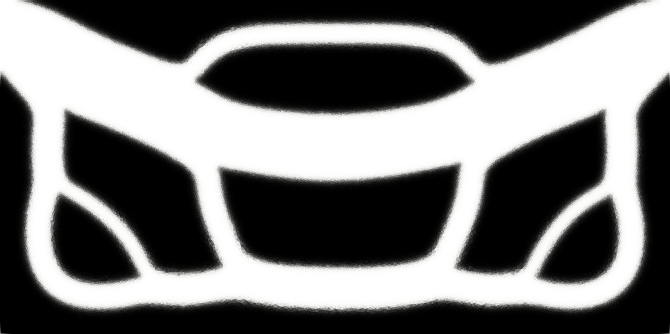}                        \\ \hline
        \rotatebox{90}{Strategy B} &
        \includegraphics[scale=0.2]{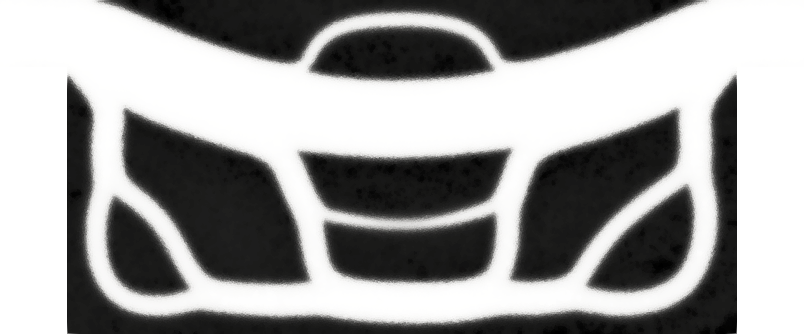}
                                         &
        \includegraphics[scale=0.2]{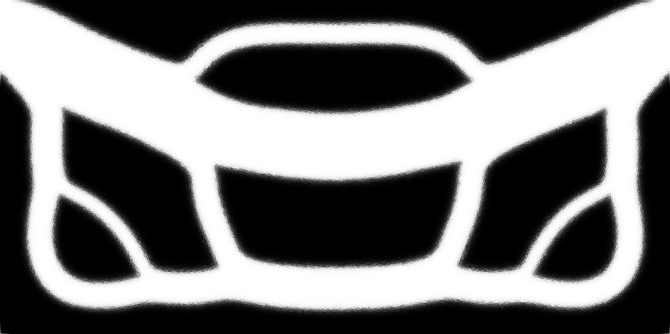}                    \\ \hline
    \end{tabular}
    \caption{Optimized designs for the thermal flow problem with AMR refinement only.}
    \label{tab:thermal_flow_refinement}
\end{table}

\begin{table}
    \centering
    \begin{tabular}{|M{0.5cm}|M{5.6cm}|M{5.6cm}|}
        \hline
                                         & Optimization in $\mathbb{R}^n$ & Optimization in $L^2$ \\ \hline
        \rotatebox{90}{Strategy A}     &
        \includegraphics[scale=0.2]{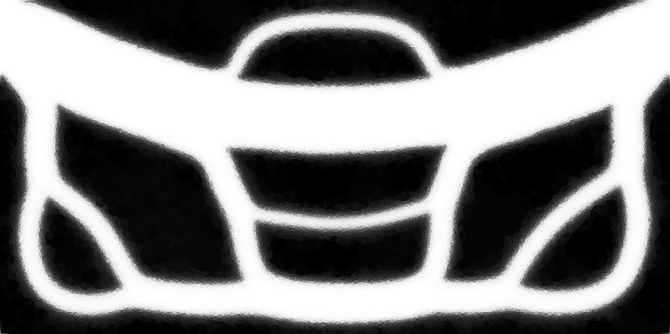}
                                         &
        \includegraphics[scale=0.2]{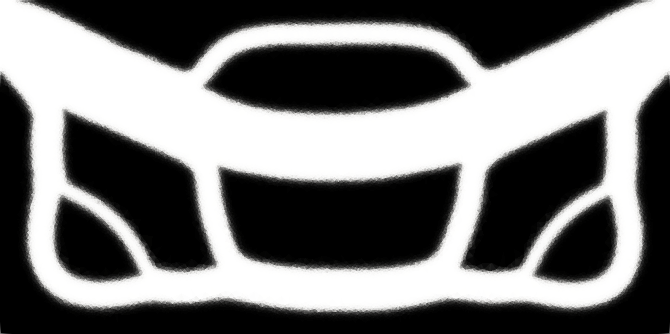}                        \\ \hline
        \rotatebox{90}{Strategy B} &
        \includegraphics[scale=0.2]{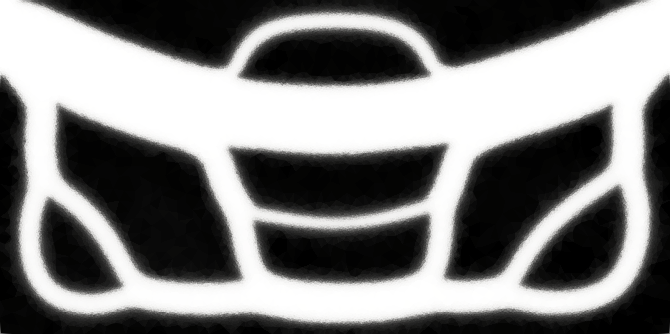}
                                         &
        \includegraphics[scale=0.2]{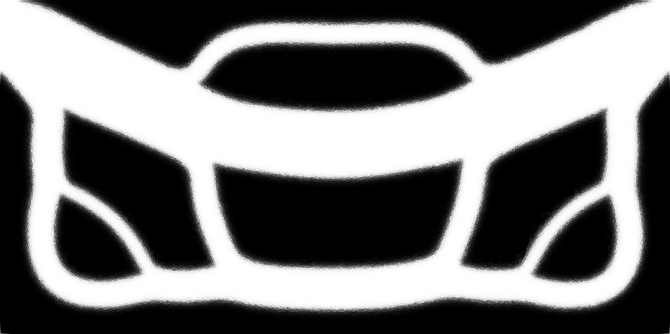}                    \\ \hline
    \end{tabular}
    \caption{Optimized designs for the thermal flow problem with AMR coarsening only.}
    \label{tab:thermal_flow_coarsening}
\end{table}

\section{Conclusion}
\label{sec:conclusions}
In this work, we presented the necessary mathematical concepts to understand
the relationship between the domain discretization and the NLP algorithm and applied them to the GCMMA algorithm.
\hl{Our $L^2$ GCMMA implementation is benchmarked with several problems in 
topology optimization and is able to obtain mesh independent designs starting 
from the same initial designs, while the original $\mathbb{R}^n$ NLP algorithm is not.
We first showed how the new algorithm solves ill-conditioned optimization problems where the $\mathbb{R}^n$ algorithm fails.
Then we showed its efficiency when solving large scale problems over uniform meshes.
Lastly, we illustrated the effectiveness of the $L^2$ algorithm when applying AMR during the optimization for two problems with different physics.}
\hl{For future work, }the algorithm can be extended to handle design fields in 
other common spaces in topology optimization such as $H^1$, i.e. for nodal design variables or for B-splines.
Lastly, a rigorous mathematical proof is necessary to ensure the NLP algorithm is 
mathematical sound for all corner cases.

\section{Replication of results}
\label{sec:replication}
The scripts used in this article are archived in \cite{firedrake_zenodo_2021_5526481} and require the pyMMAopt library \citep{miguel_salazar_de_troya_2021_4687573} 

\section{Acknowledgements}
This work was performed under the auspices of the U.S. Department of Energy by Lawrence Livermore National Laboratory under Contract DE-AC52-07NA27344. The author thanks the Livermore Graduate Scholar Program for its support.
On behalf of all authors, the corresponding author states that there is no conflict of interest. LLNL-JRNL-820905.

\bibliographystyle{elsarticle-num}
\bibliography{manuscript.bib}{}
\appendix
\section{Appendix}
\subsection{GCMMA}
\label{app:appendixA}
\normalsize
\rojo{We assume existence and uniqueness of the MMA subproblem \eqref{eq:mmasubproblem}} and solve it using a primal-dual interior-point NLP algorithm as in \cite{svanbergmma}, where a sequence of relaxed KKT conditions are solved using the Newton-Raphson method.
Omitting the iteration indices $(k, j)$ for clarity we solve \rojo{the MMA subproblem \eqref{eq:mmasubproblem} which is restated here as}
\begin{equation}
    \begin{aligned}
                       & \underset{
            \subalign{
        \nu            & \in V                                                                                                                                                                                                                       \\
        \boldsymbol{y} & \in \mathbb{R}^m                                                                                                                                                                                                            \\
        z              & \in \mathbb{R}                                                                                                                                                                                                              \\
        \boldsymbol{s} & \in \mathbb{R}^m
            }
        }{\text{min}}
                       &                  & \int_{D} \left(\frac{p_{o}}{U - \nu}+\frac{q_{o}}{\nu - L}\right) dV +a_0z +  \sum\limits_{i=1}^m \left( c_i y_i + \frac{1}{2}d_i y_i^2 \right) + r_0 \,,                                                \\
                       & \text{s.t.}
                       &                  & \int_{D} \left(\frac{p_{i}}{U - \nu}+\frac{q_{i}}{\nu - L}\right) dV - a_i z - y_i + r_i + s_i = 0, \; i = 1, \ldots, m \,,                                                                              \\
                       &                  &                                                                                                                                                           & \alpha \leq \nu \leq \beta \; \text{a.e.}\,, \\
                       &                  &                                                                                                                                                           & \boldsymbol{y} \geq 0 \,,                    \\
                       &                  &                                                                                                                                                           & z\geq 0 \,,                                  \\
                       &                  &                                                                                                                                                           & \boldsymbol{s}\geq 0 \,,
    \end{aligned}
    \label{eq:subprob_cost}
\end{equation}
where we have introduced the slack variable $\boldsymbol{s} \in \mathbb{R}^m$ to transform the  $\theta_i(\nu) \leq 0$ inequality constraints into $m$ equality constraints.

\rojo{
    To solve the above, we introduce log-barrier terms for the inequality constraints.
    The Lagrangian for the resulting penalized problem is}
\newcommand{\Lagr}{\mathcal{L}}
\newcommand{\logbarrier}[1]{\varepsilon \int_{D} \text{ln} (#1) ~dV}
\begin{equation}
    \begin{aligned}
        \Lagr (\nu,\mathbf{y},z,\lambda, \mathbf{s}) = &
        \int_{D} \left(\frac{p_{o}}{U - \nu}+\frac{q_{o}}{\nu - L}\right) dV +a_0z +  \sum\limits_{i=1}^m \left( c_i y_i +  \frac{1}{2}d_i y_i^2 \right) + r_0                                                                                          \\
                                                       & +\sum\limits_{i=1}^m \lambda_i \left( \int_{D} \left(\frac{p_{i}}{U - \nu}+\frac{q_{i}}{\nu - L}\right) dV - a_i z - y_i + r_i + s_i \right)                                                   \\
                                                       & - \epsilon \int_{D} \text{ln}( \nu - \alpha) ~dV - \epsilon \int_{D}  \text{ln} (\beta - \nu) dV - \epsilon \sumonm \text{ln} y_i - \epsilon \text{ln}(z) - \epsilon \sumonm \text{ln} s_i \,,
    \end{aligned}
    \label{eq:lagrasub}
\end{equation}
where $\lambda_i;~i=1\dots m$ are the Lagrange multipliers for the inequality constraints.

Stationarity of the Lagrangian is then solved sequentially for different values of the barrier parameter $\epsilon$, starting at $\epsilon=1$ and decreasing by a factor of 0.1 until $\epsilon \leq 10^{-5}$.
The MMA subproblem iteration loop terminates when the norm of the derivative of $\Lagr$ is less than $0.9\epsilon$.
Note that the bound constraints on the volume fraction field $\nu\in V$ are enforced via integrals \rojo{instead of via summations as in the original NLP algorithm.}

To keep the notation clean, we group the variables as $\chi=\left( \nu,\mathbf{y},z,\boldsymbol\lambda, \mathbf{s}\right)$ and express the function $\psi_j (\nu): V \rightarrow V,~ j = 0, \ldots, m $ as
\begin{equation}
    \psi_j (\nu) = \frac{p_{j}}{U - \nu} + \frac{q_{j}}{\nu - L},~ j = 0, \ldots, m  \,,
    \label{eq:psioperator}
\end{equation}
In this way, the stationary conditions on the Lagrangian become
\begin{equation}
    \begin{aligned}
        D_{\nu}\Lagr (\chi)[ \delta \nu]             & = \int_D \left( D\psi_0(\nu)[\delta \nu] + \sumonm \lambda_i D\psi_i(\nu)[\delta \nu] - \frac{\epsilon \delta \nu}{\nu - \alpha} + \frac{\epsilon \delta \nu}{\beta - \nu} \right) ~ dV & = 0 & ~~\forall \delta \nu \in V \,,            \\
        D_{y_i}\Lagr (\chi)[ \delta y_i]             & = \left( c_i + d_i y_i - \lambda_i - \frac{\epsilon}{y_i} \right) \delta y_i                                                                                                            & = 0 & ~~\forall \delta y_i  ~i=1...m \,,        \\
        D_z\Lagr (\chi)[ \delta z]                   & = (a_0 - \sumonm \lambda_i a_i - \frac{\epsilon}{z}) \delta z                                                                                                                           & = 0 & ~~\forall \delta z \,,                    \\
        D_{\lambda_i}\Lagr (\chi)[ \delta \lambda_i] & = \left(\int_D \psi_i(\nu) ~dV - a_i z - y_i +r_i + s_i \right) \delta \lambda_i                                                                                                        & = 0 & ~~ \forall \delta \lambda_i ~i=1...m  \,, \\
        D_{s_i}\Lagr (\chi)[ \delta s_i]             & = (\lambda_i - \frac{\epsilon}{s_i})\delta s_i                                                                                                                                          & = 0 & ~~\forall \delta s_i ~i=1...m  \,,
    \end{aligned}\label{eq:lagra1}
\end{equation}
Equation \eqref{eq:lagra1}.b is fulfilled for any $\delta y_i$, therefore the term in the parenthesis equals zero and similarly for Equations \eqref{eq:lagra1}.c-\eqref{eq:lagra1}.e.
From the last Equation \eqref{eq:lagra1}.e, we obtain $\lambda_i s_i = \epsilon\;\text{for}~i=1...m$.

We next introduce the variables
\begin{equation}
    \begin{aligned}
        \varepsilon & = \frac{\epsilon}{\nu - \alpha}  \,,    \\
        \eta        & = \frac{\epsilon}{\beta - \nu}  \,,     \\
        \mu_i       & = \frac{\epsilon}{y_i}\,, ~i=1...m  \,, \\
        \zeta       & = \frac{\epsilon}{z} \,,
    \end{aligned}
\end{equation}
where $\varepsilon, \eta \in V$ and $\mu_i, \zeta \in \mathbb{R}$ to generate the following system of nonlinear equations
\begin{equation}
    \begin{aligned}
        D_{\nu}\Lagr (\chi)[\delta \nu] & = \int_D \left( D\psi_0(\nu)[\delta \nu] + \sumonm \lambda_i D\psi_i(\nu)[\delta \nu] - \varepsilon\delta \nu + \eta\delta \nu \right) ~ dV & = 0 & ~~\forall \delta \nu \in V \,, \\
        D_{y_i}\Lagr (\chi)             & =  c_i + d_i y_i - \lambda_i - \mu_i                                                                                                        & = 0 & ~i=1...m \,,                   \\
        D_z\Lagr (\chi)                 & =  a_0 - \sumonm \lambda_i a_i  - \zeta                                                                                                     & = 0 & \,,                            \\
        D_{\lambda_i}\Lagr (\chi)       & = \int_D \psi_i(\nu) ~dV - a_i z - y_i +r_i + s_i                                                                                           & = 0 & ~i=1...m  \,,                  \\
        \lambda_i s_i                   & =\epsilon                                                                                                                                   &     & ~ ~i=1...m \,,                 \\
        (\nu - \alpha) \varepsilon      & =\epsilon   \,,                                                                                                                                                                    \\
        (\beta - \nu) \eta              & =\epsilon  \,,                                                                                                                                                                     \\
        \mu_i y_i                       & =\epsilon                                                                                                                                   &     & ~ ~i=1...m \,,                 \\
        \zeta z                         & =\epsilon  \,.
    \end{aligned}
    \label{eq:kkt_subproblem}
\end{equation}
We express the above in a more compact form $F(\upsilon)=0$, where $\upsilon=\left( \nu,\mathbf{y},z,\boldsymbol\lambda, \mathbf{s}, \varepsilon, \eta, \boldsymbol\mu, \zeta\right) $ and
$F:V \times \mathbb{R}^m \times \mathbb{R}\times \mathbb{R}^m\times \mathbb{R}^m \times V \times V \times \mathbb{R}^m \times \mathbb{R}
    \rightarrow V^* \times \mathbb{R}^m \times \mathbb{R}\times \mathbb{R}^m\times\mathbb{R}^m \times V \times V \times \mathbb{R}^m \times \mathbb{R}  $ is defined as
\begin{equation}
    F(\upsilon) =
    \left[	\begin{array}{l}
            \int_D \left( D\psi_0(\nu)[\delta \nu] +\textstyle \sumonm \lambda_i D\psi_i(\nu)[\delta \nu] - \varepsilon \delta \nu + \eta \delta \nu\right) ~dV \\
            c_i + d_i y_i - \lambda_i - \mu_i                                                                                                                   \\
            a_0 - \sumonm\lambda_i a_i   - \zeta                                                                                                                \\
            \int_D \psi_i(\nu) ~dV - a_i z - y_i +r_i + s_i ~                                                                                                   \\
            \lambda_i s_i -\epsilon                                                                                                                             \\
            \int_D ((\nu - \alpha) \varepsilon - \epsilon) \delta \varepsilon ~ dV                                                                              \\
            \int_D ((\beta - \nu) \eta - \epsilon ) \delta \eta ~ dV                                                                                            \\
            \mu_i y_i -\epsilon                                                                                                                                 \\
            \zeta z -\epsilon
        \end{array} \right]
    =
    \begin{bmatrix}
        \delta_{\nu}         \\
        \delta_{y_i}         \\
        \delta_{z}           \\
        \delta_{\lambda_i}   \\
        \delta_{s_i}         \\
        \delta_{\varepsilon} \\
        \delta_{\eta}        \\
        \delta_{\mu_i}       \\
        \delta_{\zeta}
    \end{bmatrix}
    \label{eq:newton_system}
\end{equation}
for all $\delta \nu \in V$, $\delta \varepsilon \in V$ and $\delta \eta \in V$.
As seen above, for discretization purposes, we enforce equations \eqref{eq:kkt_subproblem}.f and \eqref{eq:kkt_subproblem}.g weakly.

We solve $F(\upsilon)=0$ (for all $\delta \upsilon$) via Newton-Raphson by linearizing around the iterate $\upsilon^{(l)}$ and \rojo{requiring the update $\Delta \upsilon^{(l)}$} to satisfy
\begin{equation}
    D F(\upsilon^{(l)}) \Delta \upsilon^{(l)} = - F(\upsilon^{(l)})  \,,
\end{equation}
which we expand into
\begin{equation}
    \left[	\begin{array}{l}
            \int_D \left( D^2\psi_0(\nu)[\delta \nu,\Delta \nu] + \sumonm \lambda_i D^2\psi_i(\nu)[\delta \nu, \Delta \nu] +\sumonm D\psi_i(\nu)[\delta \nu] \Delta \lambda_i - \delta \nu \Delta \varepsilon + \delta \nu \Delta \eta \right) ~dV \\
            d_i \Delta y_i - \Delta \lambda_i - \Delta \mu_i                                                                                                                                                                                       \\
            - \sumonm a_i \Delta \lambda_i - \Delta \zeta                                                                                                                                                                                          \\
            \int_D D\psi_i(\nu)[\Delta \nu] ~dV - a_i \Delta z - \Delta y_i + \Delta s_i                                                                                                                                                           \\
            \Delta \lambda_i s_i + \lambda_i \Delta s_i                                                                                                                                                                                            \\
            \int_D \left( \Delta \nu \varepsilon  + (\nu - \alpha ) \Delta \varepsilon \right) \delta \varepsilon ~dV                                                                                                                              \\
            \int_D \left(- \Delta \nu  \eta  + (\beta - \nu) \Delta \eta \right) \delta \eta ~dV                                                                                                                                                   \\
            \Delta \mu_i y_i + \mu_i \Delta y_i                                                                                                                                                                                                    \\
            \Delta \zeta z  +    \zeta \Delta z
        \end{array} \right]
    =
    \begin{bmatrix}
        -\delta_{\nu}         \\
        -\delta_{y_i}         \\
        -\delta_{z}           \\
        -\delta_{\lambda_i}   \\
        -\delta_{s_i}         \\
        -\delta_{\varepsilon} \\
        -\delta_{\eta}        \\
        -\delta_{\mu_i}       \\
        -\delta_{\zeta}
    \end{bmatrix}
    \label{eq:newtonlinearsystem}
\end{equation}
for all $\delta \nu, \delta \varepsilon$ and $\delta \eta$ in $V$.

We proceed to discretize the fields $\nu, \varepsilon$ and $\eta$ in $V$ to be piecewise uniform over the finite elements as in Equation \eqref{eq:l2_disc}, \rojo{i.e. using the piecewise uniform basis functions $\mathscr{P} = \left\{\phi_1, ..., \phi_n \right\}$.}
So that, e.g.
\begin{equation}
    \begin{aligned}
        \Delta \nu(\mathbf{x}) = \boldsymbol{\Delta \mathbf{\nu}}^T \boldsymbol{\phi} (\mathbf{x})\,, \\
        \delta \nu (\mathbf{x})= \boldsymbol{\delta \mathbf{\nu}}^T \boldsymbol{\phi}(\mathbf{x})\,.
    \end{aligned}
\end{equation}
Substituting these expressions in Equations \eqref{eq:newton_system} and \eqref{eq:newtonlinearsystem} and using the arbitrariness of $\delta \nu, \delta \varepsilon$ and $\delta \eta$ yields the discretized residuals
\newcommand{\bsi}{\boldsymbol\phi}
\begin{equation}
    \begin{aligned}
        \boldsymbol{\delta_{\nu}}         & = \int_D \left( D\psi_0(\nu) +\textstyle \sumonm \lambda_i D\psi_i(\nu) - \varepsilon + \eta\right) \boldsymbol{\phi} ~dV   \,,                     \\
        \boldsymbol{\delta_{\varepsilon}} & =     (\boldsymbol{\nu} - \boldsymbol{\alpha}) \circ \boldsymbol{\varepsilon}  - \epsilon \boldsymbol{1} \,, \label{eq:epsi_hadamard} \footnotemark \\
        \boldsymbol{\delta_{\eta}}        & =     (\boldsymbol{\beta} - \boldsymbol{\nu} ) \circ \boldsymbol{\eta}  - \epsilon\boldsymbol{1}  \,,
    \end{aligned}
\end{equation}
\footnotetext{Note that for the choice of $\boldsymbol\phi$ being piecewise uniform over the elements $\int_D \left( \varepsilon (\nu - \alpha) - \epsilon \right) \boldsymbol\phi ~dV = \text{diag}(|\Omega_e|^1, |\Omega_e|^2, ...,) \left( \boldsymbol{\varepsilon} \circ (\boldsymbol{\nu} - \boldsymbol{\alpha}) - \epsilon \boldsymbol{1} \right)$. Equation \eqref{eq:epsi_hadamard}.b follows since $\text{diag}(|\Omega_e|^1, |\Omega_e|^2, ...,)$ is invertible.}
where $\boldsymbol{1}=(1,1,...,1)$ and the operator $\circ$ denotes the Hadamard product, i.e. the component wise multiplication between two vectors.
We also define
\begin{equation}
    \begin{aligned}
        \boldsymbol{\delta_y}         & =  \boldsymbol{c} + \boldsymbol{d} \circ \boldsymbol{y} - \boldsymbol{\lambda} - \boldsymbol{\mu} \,,                                            \\
        \delta_{z}                    & = a_0 - \boldsymbol{\lambda^T a} - \zeta \,,                                                                                                     \\
        \boldsymbol{\delta_{\lambda}} & = \left[\int_D \psi_1(\nu) ~dV, ... ,\int_D \psi_m(\nu) ~dV \right]^T - \boldsymbol{a} z - \boldsymbol{y} + \boldsymbol{r} + \boldsymbol{s}  \,, \\
        \boldsymbol{\delta_{s}}       & = \boldsymbol{\lambda} \circ \boldsymbol{s} - \epsilon \boldsymbol{1} \,,                                                                        \\
        \boldsymbol{\delta_{\mu}}     & = \boldsymbol{\mu} \circ \boldsymbol{y} - \epsilon \boldsymbol{1} \,,                                                                            \\
        \delta_{\zeta}                & = \zeta z - \epsilon \,.
    \end{aligned}
\end{equation}

In this way, the discretized update equations read
\newcommand{\ggb}{\mathbf{G}}
\begin{equation}
    \begin{aligned}
        \boldsymbol\Psi \boldsymbol\Delta \boldsymbol{\nu}
        + \ggb \boldsymbol{\Delta \lambda}
        - \mathbf{M} ( \boldsymbol{\Delta \varepsilon} - \boldsymbol{\Delta \eta})                                                                          & =- \boldsymbol{\delta_{\nu}}  \,,         \\
        \langle \boldsymbol{d} \rangle \boldsymbol{\Delta y} - \boldsymbol{\Delta \lambda} - \boldsymbol{\Delta \mu}                                        & =- \boldsymbol{\delta_{y}}\,,             \\
        - \mathbf{a}^T \boldsymbol{\Delta \lambda} - \Delta \zeta                                                                                           & =- \delta_z  \,,                          \\
        \ggb^T \boldsymbol{\Delta \nu} - \boldsymbol{\Delta y} - \boldsymbol{a} \Delta z + \boldsymbol{\Delta s}                                            & =- \boldsymbol{\delta_{\lambda}} \,,      \\
        \langle \boldsymbol s \rangle \Delta \boldsymbol\lambda  + \langle \boldsymbol\lambda \rangle \Delta \boldsymbol s                                  & =- \boldsymbol\delta_{\boldsymbol{s}} \,. \\
        \langle \boldsymbol\varepsilon \rangle \boldsymbol{\Delta \nu} + \langle \boldsymbol\nu - \boldsymbol\alpha \rangle \boldsymbol{\Delta \varepsilon} & =- \boldsymbol{\delta_{\varepsilon}} \,,  \\
        -\langle \boldsymbol\eta \rangle \boldsymbol{\Delta \nu} + \langle \boldsymbol\beta - \boldsymbol\nu \rangle \boldsymbol{\Delta \eta}               & =- \boldsymbol{\delta_{\eta}} \,,         \\
        \langle \boldsymbol\mu \rangle \Delta \boldsymbol y + \langle \boldsymbol y \rangle \Delta \boldsymbol\mu                                           & =-\boldsymbol\delta_{\boldsymbol\mu} \,,  \\
        \zeta \Delta z + z\Delta \zeta                                                                                                                      & =- \delta_{\zeta} \,,
    \end{aligned}
    \label{eq:discretized_residual}
\end{equation}
where the operator $\langle \cdot \rangle$ is a diagonal matrix, e.g. $\langle \boldsymbol \nu \rangle = \text{diag}(\nu_1, \nu_2, ..., \nu_n)$ and
$\boldsymbol \Psi$ and $\mathbf{G}$ are the block matrices
\begin{equation}
    \begin{aligned}
        \boldsymbol\Psi & = \int_D \left( D^2\psi_0(\nu) + \sumonm \lambda_i D^2\psi_i(\nu) \right) \bsi \bsi^T ~dV \,,                                       \\
        \mathbf{G}      & = \left[\int_D D\psi_1(\nu) \bsi ~dV, \int_D D\psi_2(\nu) \bsi ~dV, ...~, \int_D D\psi_m(\nu) \bsi ~dV \right]\label{eq:matrixG} \,
    \end{aligned}
\end{equation}
\rojo{where the integrals are calculated using one quadrature point per element.
    Each component of the matrix $\mathbf{G}$ is a column vector of dimensions $n \times 1$ due to the basis functions $\boldsymbol \phi$. \label{matrixG}
    Matrix $\boldsymbol \Psi$ has dimensions $n \times n$ and $\mathbf{G}$, $n \times m$.}
Ultimately, we are left with the update equation
\newcommand{\mm}{\mathbf{M}}
\begin{equation}
    \begin{pmatrix}
        \boldsymbol\Psi                        &                                &                 & \ggb                           &                                      & -\mm                                                & \mm                                                 &                            &    \\
                                               & \langle \mathbf{d} \rangle     &                 & -\mathbf{I}_m                  &                                      &                                                     &                                                     & -\mathbf{I}_m              &    \\
                                               &                                &                 & -\boldsymbol{a}^T              &                                      &                                                     &                                                     &                            & -1 \\
        \ggb^T                                 & -\mathbf{I}_m                  & -\boldsymbol{a} &                                & \mathbf{I}_m                         &                                                     &                                                     &                                 \\
                                               &                                &                 & \langle \boldsymbol{s} \rangle & \langle \boldsymbol   \lambda\rangle &                                                     &                                                     &                            &    \\
        \langle \boldsymbol\varepsilon \rangle &                                &                 &                                &                                      & \langle \boldsymbol\nu  - \boldsymbol\alpha \rangle &                                                     &                            &    \\
        -\langle \boldsymbol\eta \rangle       &                                &                 &                                &                                      &                                                     & \langle \boldsymbol\beta -  \boldsymbol\nu  \rangle &                            &    \\
                                               & \langle \boldsymbol\mu \rangle &                 &                                &                                      &                                                     &                                                     & \langle \mathbf{y} \rangle &    \\
                                               &                                & \zeta           &                                &                                      &                                                     &                                                     &                            & z
    \end{pmatrix}
    \begin{pmatrix}
        \boldsymbol\Delta \boldsymbol \nu         \\
        \boldsymbol\Delta \boldsymbol y           \\
        \Delta z                                  \\
        \boldsymbol\Delta \boldsymbol \lambda     \\
        \boldsymbol\Delta \boldsymbol s           \\
        \boldsymbol\Delta \boldsymbol \varepsilon \\
        \boldsymbol\Delta \boldsymbol \eta        \\
        \boldsymbol\Delta \boldsymbol \mu         \\
        \Delta \zeta
    \end{pmatrix}
    =
    \begin{pmatrix}
        -\boldsymbol{\delta_{\nu}}         \\
        -\boldsymbol{\delta_y}             \\
        -{\delta_z}                        \\
        -\boldsymbol{\delta_{\lambda}}     \\
        -\boldsymbol{\delta_{s}}           \\
        -\boldsymbol{\delta_{\varepsilon}} \\
        -\boldsymbol{\delta_{\eta}}        \\
        -\boldsymbol{\delta_{\mu}}         \\
        -{\delta_{\zeta}}
    \end{pmatrix} \,,
    \label{eq:hessian_matrix}
\end{equation}
where $\mathbf{I}_m$ is the $m\times m$ identity matrix.
The pointwise operations in Equations \eqref{eq:alphabeta}, \eqref{eq:LUupdate1} - \eqref{eq:LUrestrict} are carried out on the vector components of their corresponding discretized fields.
\newcommand{\bsdelta}[1]{\boldsymbol{\Delta #1}}

As in \cite{svanbergmma}, the linear system of equations in \eqref{eq:hessian_matrix} can be solved for
$\boldsymbol{\Delta \upsilon} = (\bsdelta{\nu}, \bsdelta{y}, \Delta z, \bsdelta{\lambda}, \bsdelta{s}, \bsdelta{\varepsilon}, \bsdelta{\eta}, \bsdelta{\mu}, \Delta \zeta)$
quickly given that all the block matrices are diagonal, although this only happens for the element-wise uniform discretization.

We limit the Newton step from violating the inequalities in Equations \eqref{eq:subprob_cost} as in \cite{svanbergmma}.
To do this, we first maximize the parameter $t$ such that
\begin{equation}
    \begin{aligned}
        t                                            & \leq 1 \label{eq:ls1}\,,                                                                                                                       \\
        \nu_l + t\Delta \nu_l - \alpha_l             & \geq 0.01\left(\nu_l - \alpha_l \right) ~ \text{for all  $l=1,...,n$} \,,                                                                      \\
        \beta_l - \left(\nu_l + t\Delta \nu_l\right) & \geq 0.01\left(\beta_l - \nu_l  \right) ~ \text{for all  $l=1,...,n$} \,,                                                                      \\
        \left( \mathbf{y}, z , \boldsymbol\lambda, \boldsymbol s,\boldsymbol\varepsilon, \boldsymbol\eta, \boldsymbol\mu, \zeta, \right)
        + t  \left( \boldsymbol\Delta \mathbf{y}, \Delta z , \boldsymbol\Delta \boldsymbol\lambda, \boldsymbol\Delta \boldsymbol s,\boldsymbol\Delta \boldsymbol\varepsilon, \boldsymbol\Delta \boldsymbol\eta, \boldsymbol\Delta \boldsymbol\mu, \Delta \zeta \right)
                                                     & \geq 0.01 \left( \mathbf{y}, z , \boldsymbol\lambda, \boldsymbol s, \boldsymbol\varepsilon, \boldsymbol\eta, \boldsymbol\mu, \zeta \right) \,.
    \end{aligned}
\end{equation}
where all inequalities are performed on the vector components and $n$ is the number of mesh elements.
As in \cite{svanbergmma}, we assign $\boldsymbol\upsilon \leftarrow \boldsymbol\upsilon + \tau \bsdelta{\upsilon}$ where $\tau$ is the largest of $t,~t/2,~t/4,...$ such that
\begin{align}
    \norm{F(\upsilon + \tau \Delta \upsilon)}_Q < \norm{F(\upsilon)}_Q \,.\label{eq:ls5}
\end{align}
We emphasize that $\norm{\cdot}_Q$ is not the ordinary \rojo{vector} norm, but the norm induced by the space $Q= V^* \times \mathbb{R}^m \times \mathbb{R}\times \mathbb{R}^m \times\mathbb{R}^m \times V \times V \times \mathbb{R}^m \times \mathbb{R} $, i.e.

\begin{equation}
    \norm{F(\upsilon)}_Q = \sqrt{ \norm{\delta_{\nu}}_{V^*_h}^2 +
    \norm{\delta_y}^2_{\mathbb{R}^m} +
    \lVert\delta_z \rVert^2 +
    \norm{\delta_{\lambda}}^2_{\mathbb{R}^m} +
    \norm{\delta_s}^2_{\mathbb{R}^m} +
    \norm{\delta_{\varepsilon}}_{V_h}^2 +
    \norm{\delta_{\eta}}_{V_h}^2 +
    \norm{\delta_{\mu}}^2_{\mathbb{R}^m} +
    \lVert\delta_{\zeta}\rVert^2} \,,
    \label{eq:norm_resi_mma}
\end{equation}
where, e.g.
\begin{equation}
    \begin{aligned}
        \norm{\delta_{\nu}}_{V_h^*}       & = \boldsymbol\delta_{\nu}^T \mathbf{M}^{-1} \boldsymbol\delta_{\nu} \,,             \\
        \norm{\delta_{\varepsilon}}_{V_h} & = \boldsymbol\delta_{\varepsilon}^T \mathbf{M} \boldsymbol\delta_{\varepsilon}  \,.
    \end{aligned}
    \label{eq:norm_partial_resi}
\end{equation}
To recover the original NLP algorithm by \cite{svanbergmma} posed in the \rojo{sequence} space $\mathbb{R}^m$, we merely replace $\mathbf{M}$ with the identity matrix in the norm calculations and Equation \eqref{eq:hessian_matrix}.

The discretized Fr\'echet derivative $\mathbf{D}\boldsymbol\theta$ needs to be passed to the NLP algorithm.
However, the GCMMA requires the discretized gradient $\boldsymbol{\nabla\theta}$ cf. Equations \eqref{eq:pfunc}, \eqref{eq:qfunc} and \eqref{eq:kkt_stopping}.a which is calculated using the Riesz map \eqref{eq:rieszmap} as
\begin{equation}
    \begin{aligned}
        \boldsymbol{\nabla \theta} & = \Phi_h^{-1}(\mathbf{D}\boldsymbol{\theta}) \,,    \\
                                   & = \mathbf{M}^{-1} \mathbf{D}\boldsymbol{\theta} \,.
    \end{aligned}
\end{equation}
The gradient components are used to evaluate Equations \eqref{eq:pfunc} and \eqref{eq:qfunc}, i.e.
\begin{equation}
    \begin{aligned}
        p^{(k,j)}_{i,l} = (U_l^{(k)} - \nu_l^{(k,j)})^2\left(1.001\left(\nabla \theta_{i,l}(\nu^{(k,j)})\right)^+ +0.001\left(\nabla \theta_{i,l}(\nu^{(k,j)})\right)^- +\frac{\rho_i^{(k,j)}}{\nu_{\text{max}} - \nu_{\text{min}}}\right) \,, \\
        q^{(k,j)}_{i,l} = (\nu_l^{(k)} - L_l^{(k,j)})^2\left(0.001\left(\nabla \theta_{i,l}(\nu^{(k,j)})\right)^+ +1.001\left(\nabla \theta_{i,l}(\nu^{(k,j)})\right)^- +\frac{\rho_i^{(k,j)}}{\nu_{\text{max}} - \nu_{\text{min}}}\right) \,.
    \end{aligned}
\end{equation}
where the subscript $l$ in $p^{(k,j)}_{i,l}$ and $q^{(k,j)}_{i,l}$ corresponds to the mesh element $l$ and in $\nabla \theta_{i,l}$ it corresponds to the component of the gradient $\boldsymbol{\nabla \theta}_i$.

The termination criteria of the original problem \eqref{eq:minproblem}, i.e. not the subproblem \eqref{eq:mmasubproblem}, is derived from its KKT conditions.
Notably, here we follow the convergence metric of \cite{gcmma} and monitor the gradient of the Lagrangian
\begin{align}
    \Lagr(\nu, \lambda, \tau^+, \tau^-) & = \theta_0(\nu) + \sumonm \lambda_i \theta_i(\nu) +  \int_D \tau^-(\nu_{\text{min}} - \nu)~dV
    + \int_D \tau^+(\nu - \nu_{\text{max}})~dV \,.
\end{align}
Upon defining $\omega = (\nu, \lambda, \tau^-, \tau^+)$, the KKT conditions read
\begin{equation}
    \begin{aligned}
        \nabla_\nu    \Lagr(\omega)    & =\nabla \theta_0(\nu) + \sumonm \lambda_i \nabla \theta_i(\nu) - \tau^- + \tau^+ =0 \,,\phantom{i = 1, \ldots, m \,,}                                  \\
        \theta_i(\nu)                  & \leq 0                                                                                                                & \; ~ i = 1, \ldots, m\,,\notag \\
        \lambda_i \theta_i(\nu)        & = 0                                                                                                                   & \; ~ i = 1, \ldots, m\,,\notag \\
        \lambda_i                      & \geq 0                                                                                                                & \; ~ i = 1, \ldots, m\,,\notag \\
        \nu_{\text{min}} \leq \nu      & \leq \nu_{\text{max}} \,,\notag                                                                                                                        \\
        \tau^-(\nu_{\text{min}} - \nu) & = 0 \,,                                                                                                                                                \\
        \tau^+(\nu - \nu_{\text{max}}) & = 0 \,,                                                                                                                                                \\
        \tau^+, \tau^-                 & \geq 0\notag \,.
    \end{aligned}
    \label{eq:kkt_stopping}
\end{equation}
Wherever the bound constraints are active, i.e. wherever $\nu = \nu_{\text{min}}$, or $\nu = \nu_{\text{max}}$, their corresponding Lagrange multipliers are $\tau^- =  \nabla \theta_0(\nu) + \sum_{i=1}^m \lambda_i \nabla \theta_i(\nu) $ or $\tau^+ = - \left( \nabla \theta_0(\nu) + \sum_{i=1}^m \lambda_i \nabla \theta_i(\nu) \right)$ respectively.
And since $\tau^- \geq 0 $, we have $\nabla \theta_0(\nu) + \sum_{i=1}^m \lambda_i \nabla \theta_i(\nu) \geq 0 $ and similarly since $\tau^+ \geq 0 $ we have $\nabla \theta_0(\nu) + \sum_{i=1}^m \lambda_i \nabla \theta_i(\nu) \leq 0 $.
Using these inequalities in the complementary slackness Equations \eqref{eq:kkt_stopping}.f - \eqref{eq:kkt_stopping}.g, transforms them to
\begin{equation}
    \begin{aligned}
        \lambda_i (\theta_i(\nu))^-                                                                           & = 0\,,\; i = 1, \ldots, m  \,,   \\
        (\theta_i(\nu))^+                                                                                     & = 0\,, \; i = 1, \ldots, m \,,   \\
        \lambda_i                                                                                             & \geq 0\,,\; i = 1, \ldots, m \,, \\
        (\nu_{\text{min}} - \nu)\left( \nabla \theta_0(\nu) + \sumonm\lambda_i \nabla \theta_i(\nu) \right)^+ & = 0 \,,                          \\
        (\nu_{\text{max}} - \nu)\left( \nabla \theta_0(\nu) + \sumonm\lambda_i \nabla \theta_i(\nu) \right)^- & = 0 \,,                          \\
        \tau^+, \tau^-                                                                                        & \geq 0
    \end{aligned}
    \label{eq:kkt_stopping_final}
\end{equation}
where the Lagrange multipliers $\lambda_i$ are obtained from the solution of the convex approximation problem \eqref{eq:mmasubproblem}.
The norm of the KKT conditions is given by the norms of the left hand sides of Equations \eqref{eq:kkt_stopping_final}.a-\eqref{eq:kkt_stopping_final}.b and \eqref{eq:kkt_stopping_final}.d-\eqref{eq:kkt_stopping_final}.e where each norm is taken in its corresponding space, i.e., Equations \eqref{eq:kkt_stopping_final}.a-\eqref{eq:kkt_stopping_final}.b in $\mathbb{R}^m$ and Equations \eqref{eq:kkt_stopping_final}.d-\eqref{eq:kkt_stopping_final}.e in $V$.
The primal-dual NLP algorithm for the subproblem is outlined in Algorithm \ref{alg:subproblem}
\begin{algorithm}[!h]
    \caption{MMA subproblem algorithm outline.}
    \begin{algorithmic}[1]
        \State $\mathbf{Input:}$ Starting point for $k=1$: $\nu^{(k)} = \frac{\alpha + \beta}{2}, y^{(k)}_i = s^{(k)}_i = \lambda^{(k)}_i = 1.0, ~ \mu^{(k)}_i = \text{max}\{1.0, 0.5 c_i\}~ \text{for} ~ i=1,...,m, ~ z^{(k)}=\zeta^{(k)}=1.0, $
        $\varepsilon^{(k)} = \text{max} \{ 1.0, \frac{1.0}{\nu^{(k)} - \alpha} \},\eta^{(k)} = \text{max}\{ 1.0, \frac{1.0}{\beta - \nu^{(k)}} \}  $ and $\epsilon = 1$.
        \While{$\epsilon > 10^{-5}$}
        \State Calculate $\boldsymbol{\Delta \upsilon}$ in Equation \eqref{eq:hessian_matrix}.
        \State Calculate step length $\tau$ in Equation \eqref{eq:ls5}.
        \State Let $\boldsymbol{\upsilon}^{(k+1)} = \boldsymbol{\upsilon}^{(k)} + \tau \boldsymbol{\Delta \upsilon}^{(k)}$.
        \If{$\norm{F(\upsilon)}_Q<0.9 \epsilon$ }
        \State $\epsilon = 0.1 \epsilon$
        \EndIf
        \EndWhile
    \end{algorithmic}
    \label{alg:subproblem}
\end{algorithm}

\rojo{As seen in the calculation of the search direction, cf. Equation \eqref{eq:hessian_matrix}, the residual norm, cf. Equation \eqref{eq:norm_resi_mma} and the convergence metric, cf. Equation \eqref{eq:kkt_stopping_final}, the $L^2$ NLP algorithm cannot be obtained by merely scaling the derivatives that feed the $\mathbb{R}^n$ NLP algorithm by the mass matrix, e.g. be replacing $D\theta_i$ with $\nabla \boldsymbol \theta = \mathbf{M}^{-1} \boldsymbol{D\theta}$.
    Obtaining a mesh independent NLP algorithm requires starting from the infinite dimensional formulation and using the tools presented in this paper to obtain the proper discretization.}
Summarizing all the changes made to the GCMMA algorithm:
\begin{itemize}
    \item The design variable $\nu$, its lower and upper bounds, the moving asymptotes $L$ and $U$ and the subproblem bounds $\alpha$ and $\beta$ are functions in the Hilbert space $L^2$.
    \item The gradients $\nabla \theta_i(\nu)$ for $i=0,...,m$ are used in the convex approximation, i.e. Equation \eqref{eq:convex_approx} instead of derivatives $D\theta_i(\nu)$.
    We obtain the gradients $\nabla \theta_i(\nu)$ by applying the Riesz map as in Equation \eqref{eq:rieszmap_grad}.
    \item The convex approximation in Equation \eqref{eq:convex_approx} is built with an integral over the domain, instead of a summation over the vector of design variables.
    \item Similarly, the summations in the global convergence mechanism, i.e. Equations \eqref{eq:rho_mechanism} and \eqref{eq:d_mechanism}, are replaced with integrals over the domain.
    \item The norms used to check for convergence of convex approximation subproblem, i.e. Equation \eqref{eq:norm_resi_mma} and the convergence metric, i.e. \eqref{eq:kkt_stopping_final} are taken in the appropriate spaces
\end{itemize}

\begin{figure}[h!]
    \centering
    \includegraphics[width=\BeamL]{\MyPath/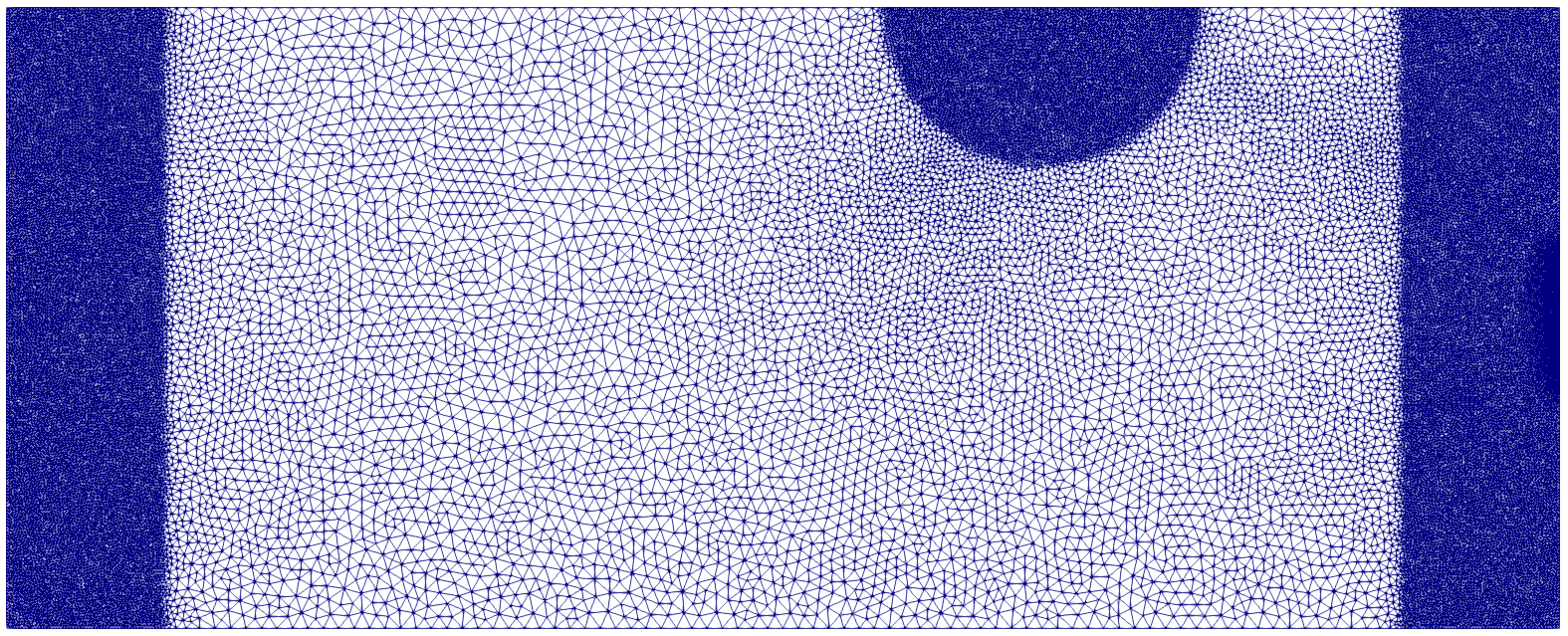}
    \caption{Non-uniform mesh for the compliance problem.}
    \label{fig:amr_compliance}
\end{figure}

\begin{figure}[h!]
    \centering
    \includegraphics[width=\BeamL]{\MyPath/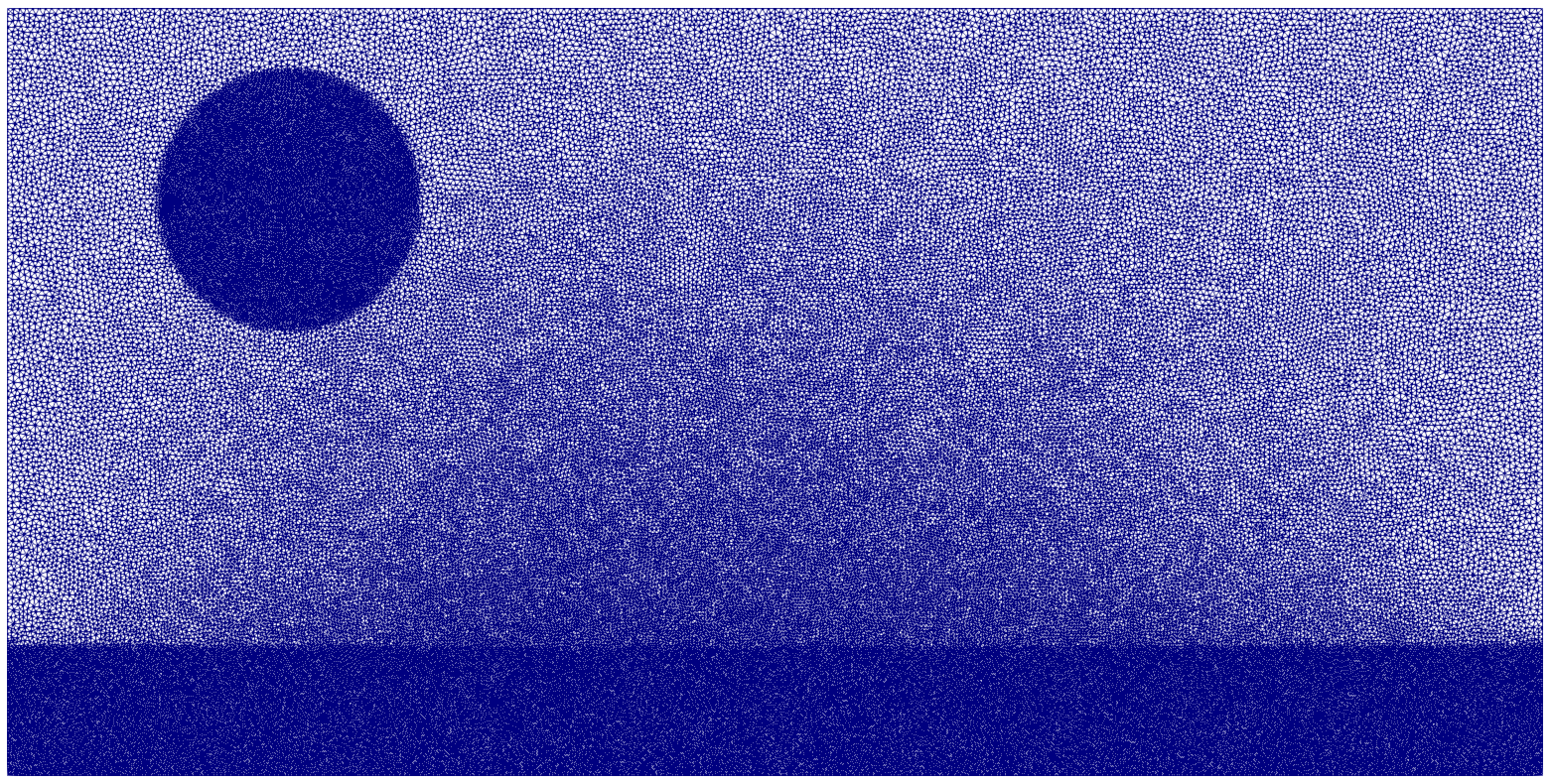}
    \caption{Non-uniform mesh for the mechanism problem.}
    \label{fig:amr_mechanism}
\end{figure}

\begin{figure}[h!]
    \centering
    \includegraphics[width=\BeamL]{\MyPath/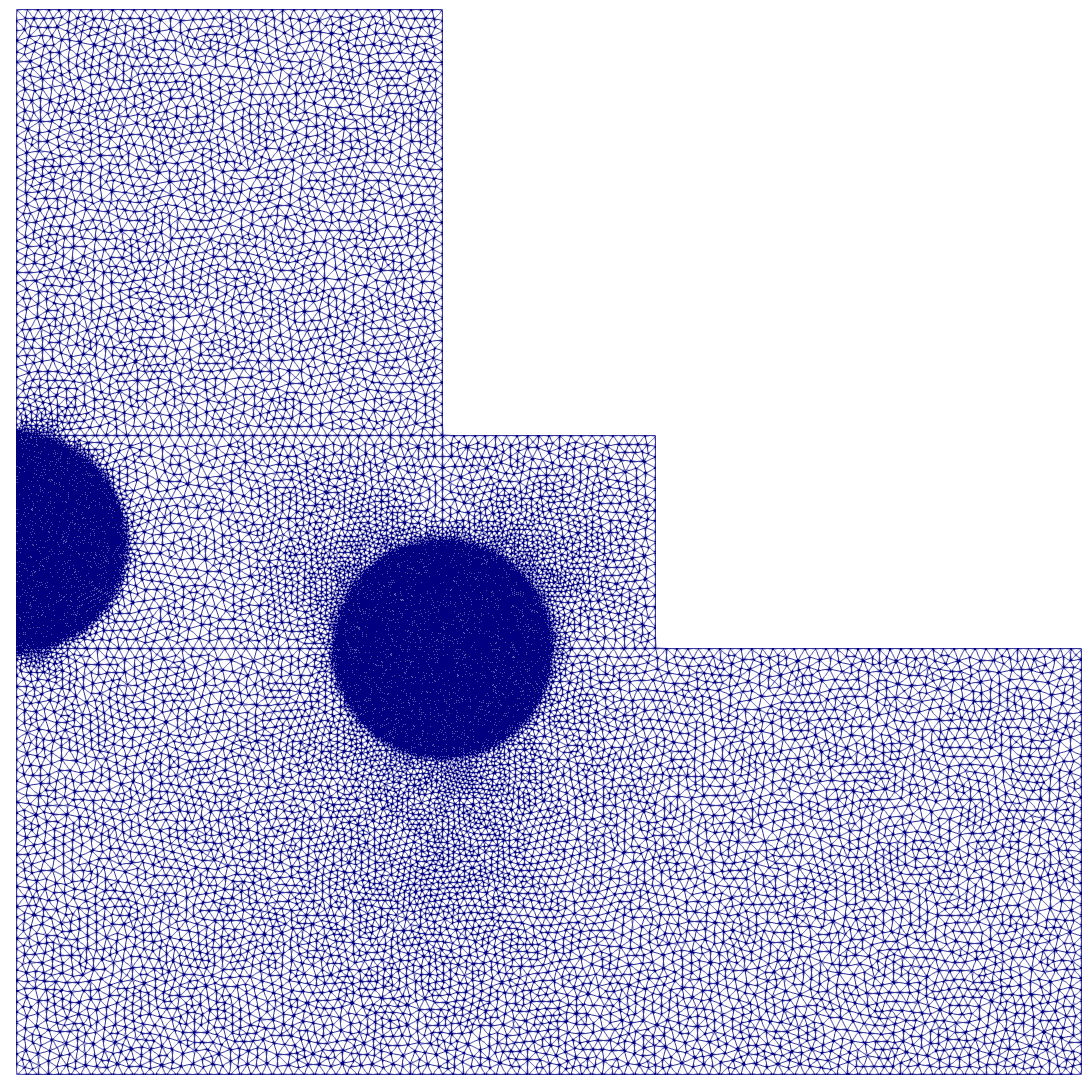}
    \caption{Non-uniform mesh for the stress constrained problem.}
    \label{fig:amr_stress}
\end{figure}

\subsection{\hl{Adaptive mesh refinement}}
\label{sec:amr_appendix}
We apply AMR using an element-based error quantity that measures the jump of the volume fraction $\hat\nu$ ($\nu$ for the thermal problem) across each element $K$ edges:
\begin{align}
    e_K = \frac{\sum_{\mathcal{E} \in K} \int_{\mathcal{E}} \llbracket \hat\nu \rrbracket ~dS}{\sum_{\mathcal{E} \in K} \int_{\mathcal{E}} ~dS}\,,
\end{align}
where $\mathcal{E}$ are the element $K$ edges.
Note that $e_K$ is properly weighted to account for the element size.

We perform either only refinement or coarsening during the AMR operation using Triangle
\citep{triangle_amr} and MeshPy \citep{meshpy}.
During refinement, we use MeshPy to update the element $\Omega_k$ area $a_K$ to
\begin{align}
    a_K = \text{max} \left( \frac{a_K}{\left(1 + 100 \frac{e_K}{\underset{K\in \mathcal{T}}{\text{max}}e_K}\right)^2}, a_{\text{min}} \right) \,.
\end{align}
$\mathcal{T}$ is the mesh discretization.
The lower area limit $a_{\text{min}}$ prevents very small elements and is equal to 0.01 and $6\times 10^{-6}$ for the compliance and thermal-flow problem respectively.

The coarsening requires a different criteria because MeshPy does not have such functionality.
Instead, we remesh the domain from scratch, letting the largest element size to be
$E_{\mathcal{T}} \left(~\underset{K\in \mathcal{T}}{\text{max }} a_K \right)$,
where the expansion factor $E_{\mathcal{T}}$ controls the coarsening.
We use $E_{\mathcal{T}}=10.0$ and 2.0 for the compliance and thermal-flow problem respectively.
Since we want to keep the design geometry refined, we define the high error region
\begin{align}
    D_R = \left\{ \mathbf{x} \in D ~|~ e_{K(\mathbf{x})} > e_{50}\right\}
\end{align}
where $e_{50}$ is the 50-th of the element error distribution.
This region is refined until the element areas are smaller than the smallest element
from the previous mesh.
$K(\mathbf{x})$ refers to those elements from the previous mesh that contain $\mathbf{x}$,

We remark here that the beam and thermal flow AMR examples start from the same respective mesh
and either refine it or coarsen it. This ensures that the initial mesh is fine enough
to approximate the infinite dimensional solution of the state variables.
An initially too coarse mesh will introduce a large error in the state variables
and the cost and constraint function derivatives which will affect the optimized design.

To ensure continuity of the optimization algorithm at iteration $k$,
the design variable $\nu_k$, the previous iterations $\nu_{k-1}$, $\nu_{k-2}$ and the MMA moving asymptotes
$L_{k}$ and $U_{k}$ are projected onto the refined or coarsened mesh.

\begin{table}
    \centering
    \begin{tabular}{|M{0.5cm}|M{5.6cm}|M{5.6cm}|}
        \hline
                                         & Optimization in $\mathbb{R}^n$ & Optimization in $L^2$ \\ \hline
        \rotatebox{90}{Strategy A}     &
        \includegraphics[scale=0.05]{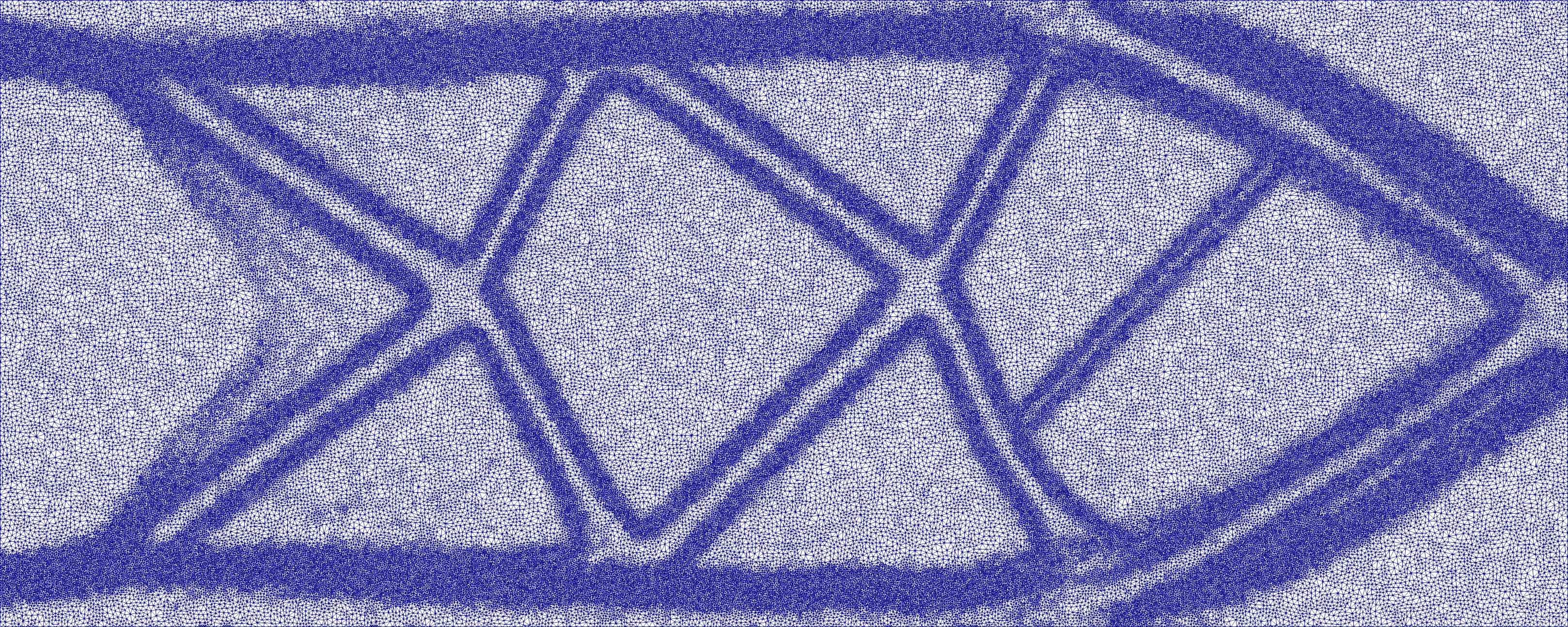}
                                         &
        \includegraphics[scale=0.05]{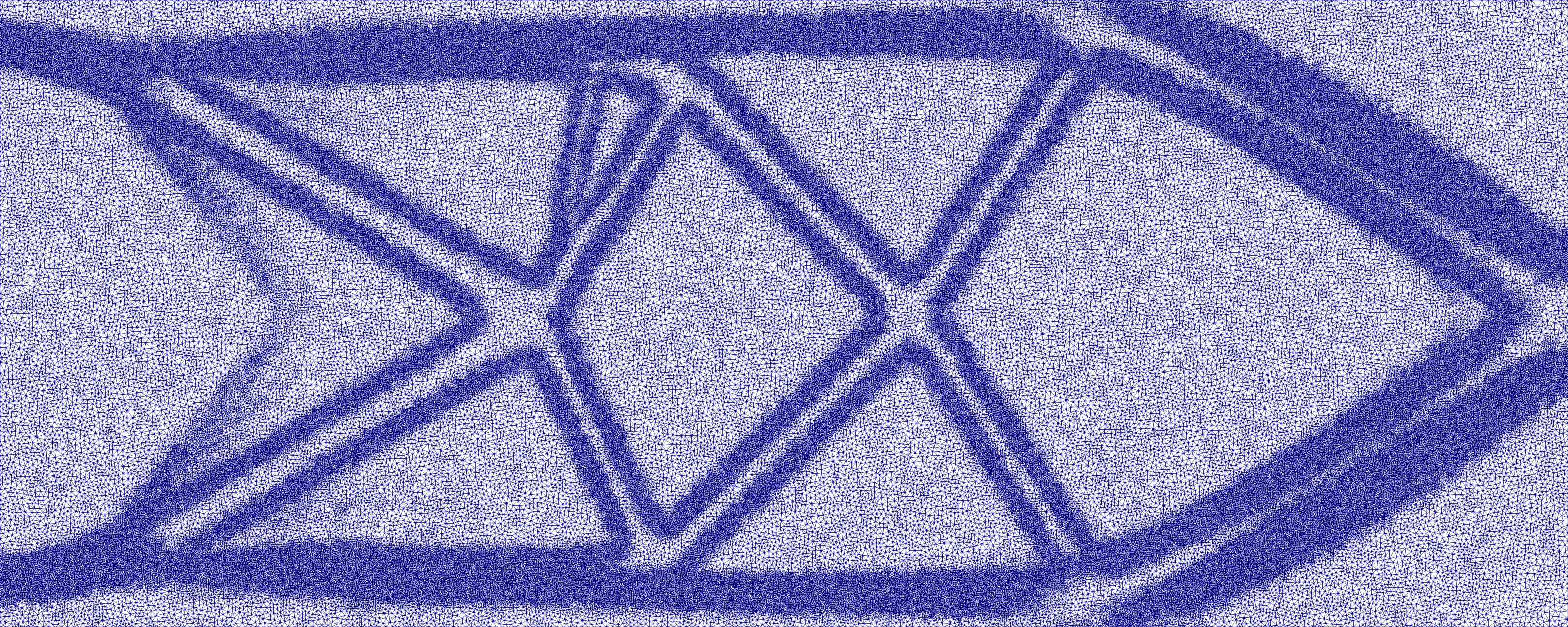}                        \\ \hline
        \rotatebox{90}{Strategy B} &
        \includegraphics[scale=0.05]{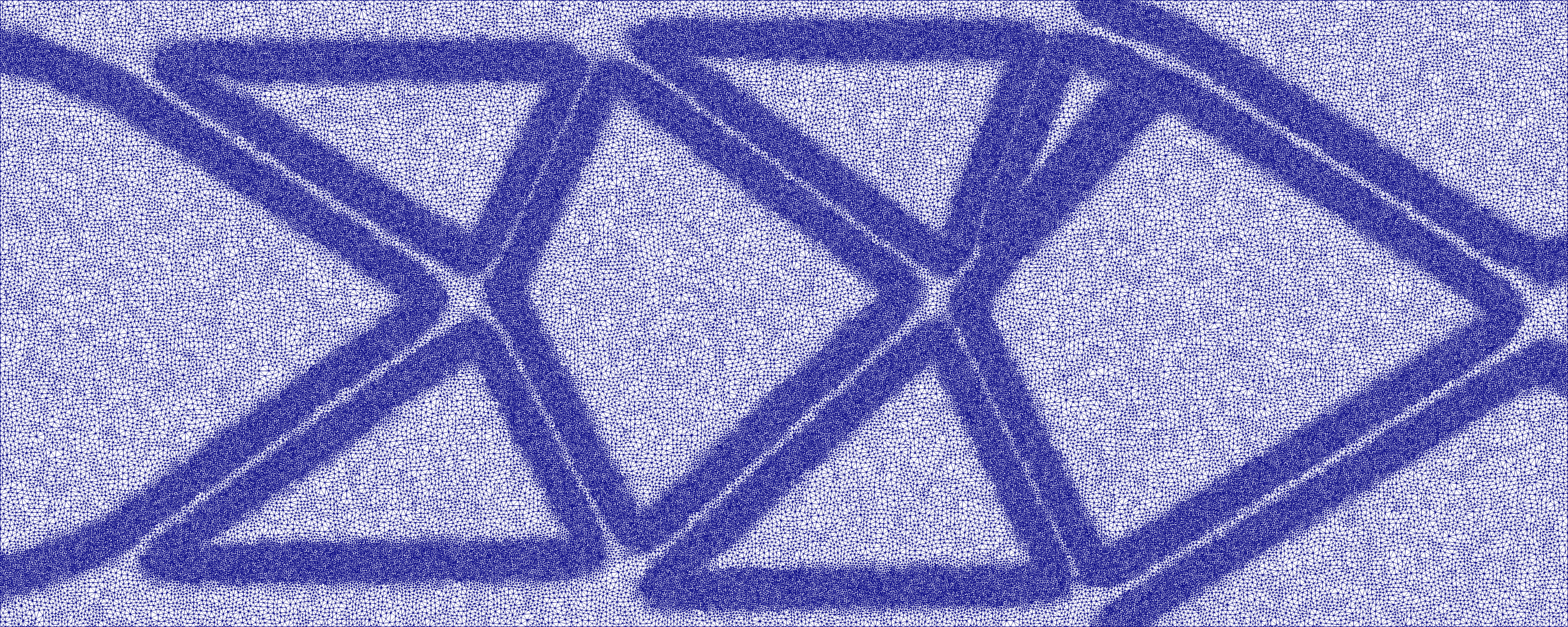}
                                         &
        \includegraphics[scale=0.05]{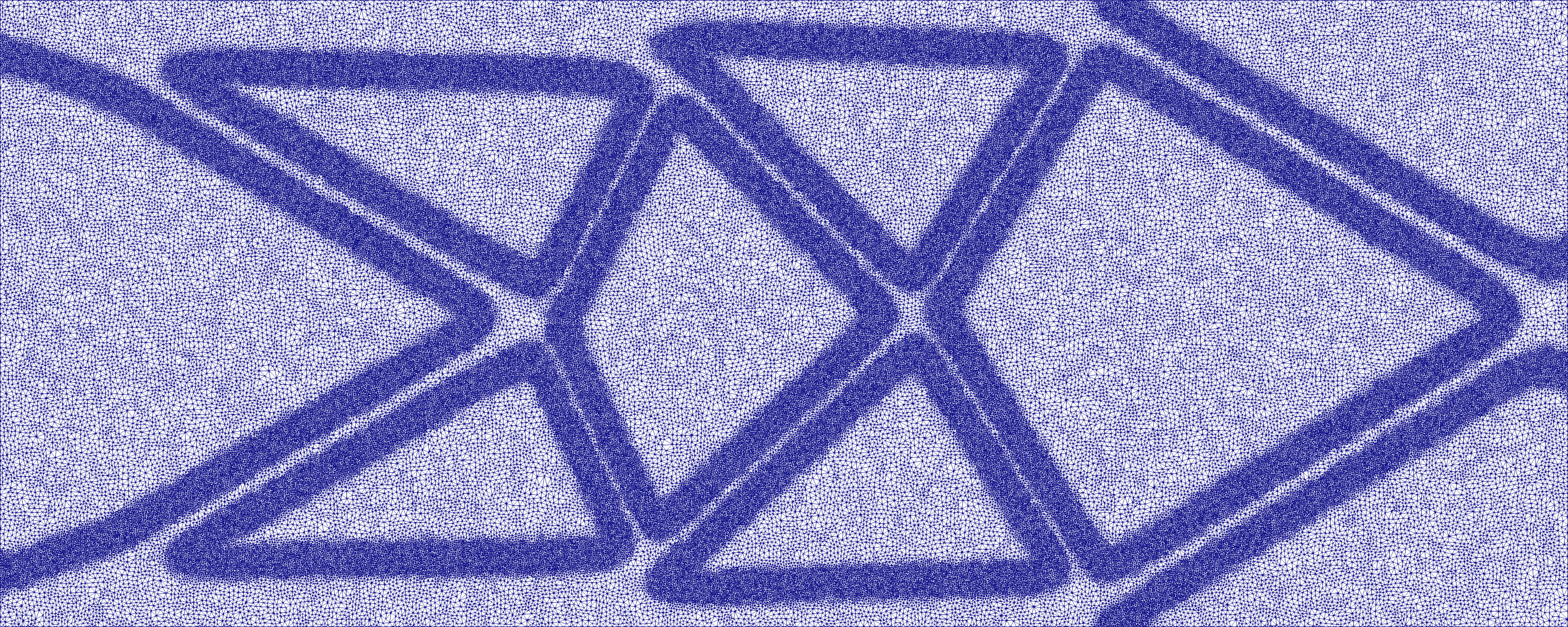}                    \\ \hline
    \end{tabular}
    \caption{Meshes for optimized designs for the AMR compliance problem with refinement.}
    \label{fig:compliance_amr_refinement_grid}
\end{table}

\begin{table}
    \centering
    \begin{tabular}{|M{0.5cm}|M{5.6cm}|M{5.6cm}|}
        \hline
                                         & Optimization in $\mathbb{R}^n$ & Optimization in $L^2$ \\ \hline
        \rotatebox{90}{Strategy A}     &
        \includegraphics[scale=0.05]{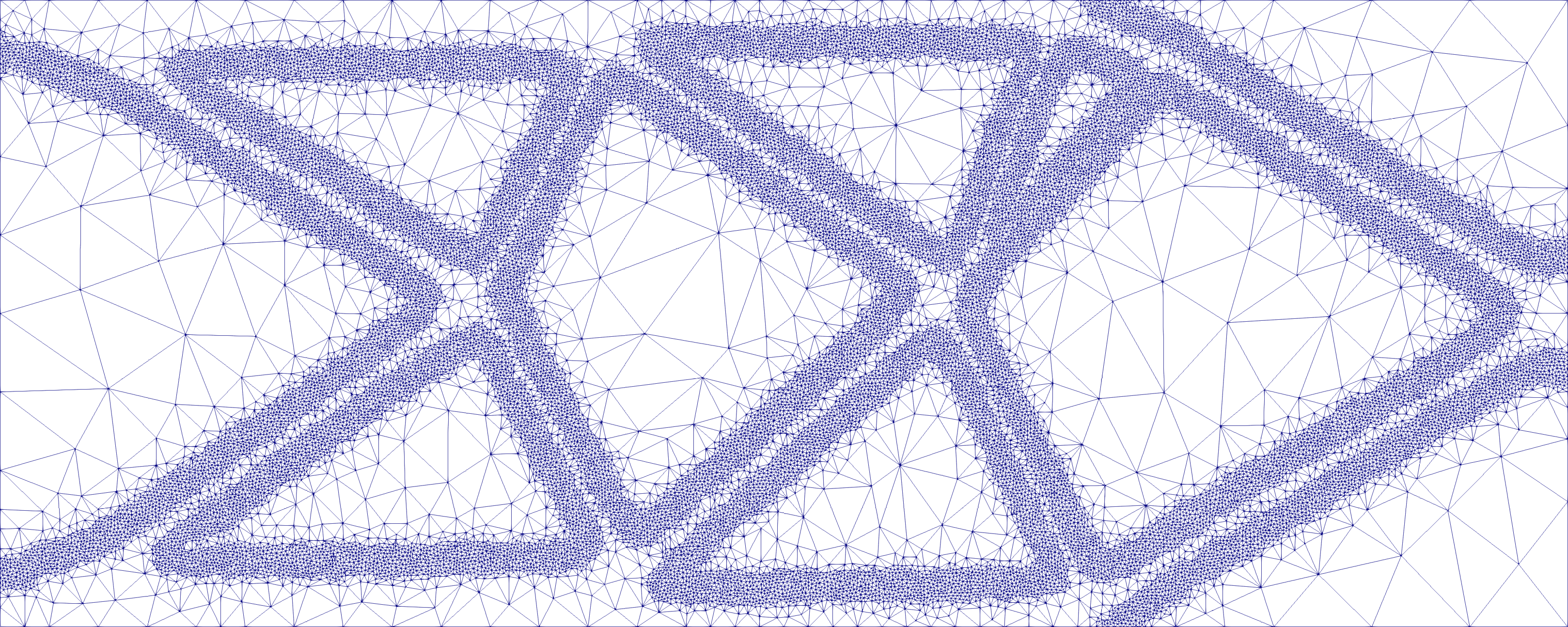}
                                         &
        \includegraphics[scale=0.05]{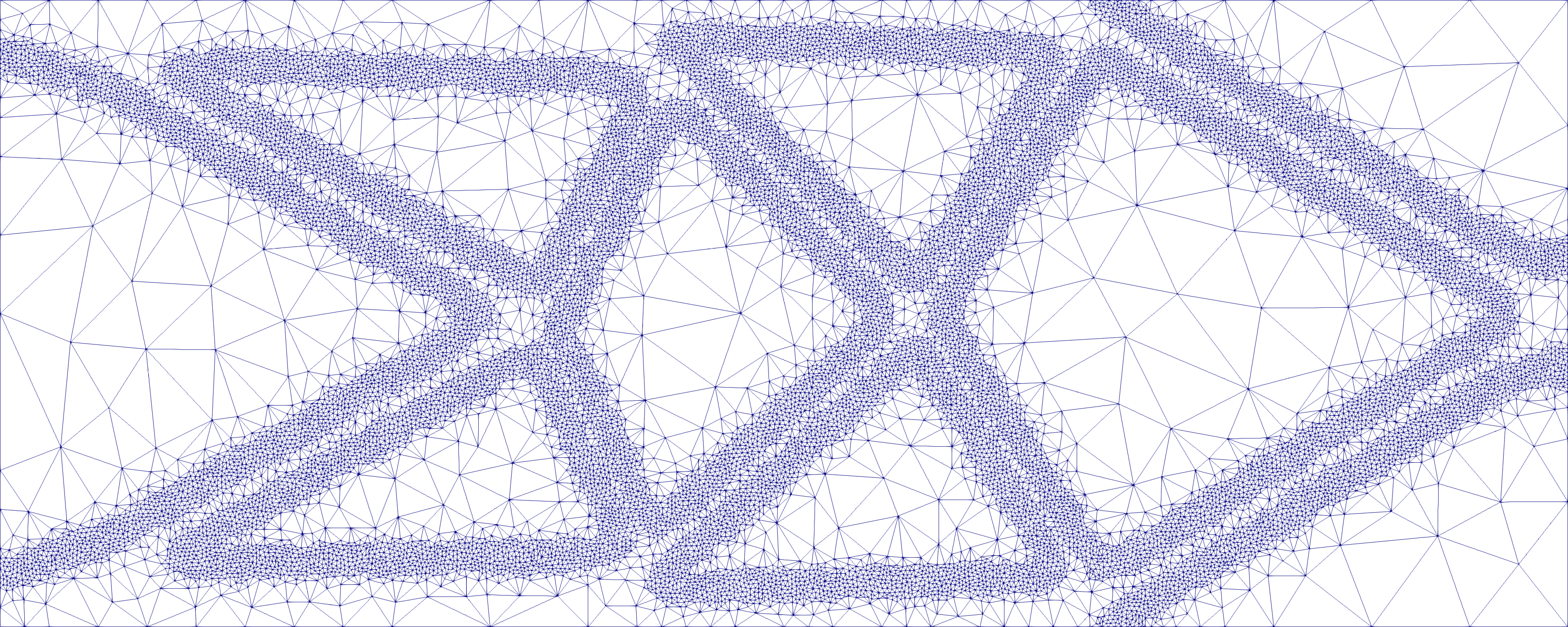}                        \\ \hline
        \rotatebox{90}{Strategy B} &
        \includegraphics[scale=0.05]{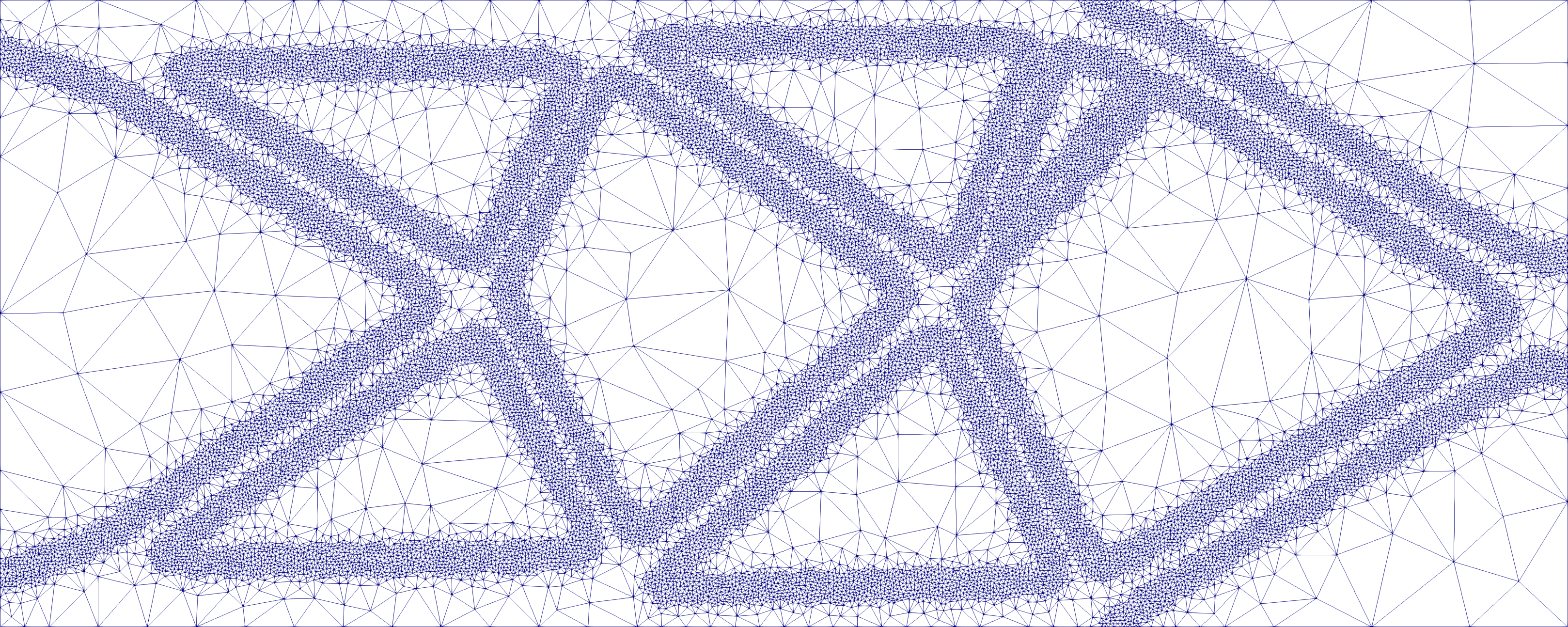}
                                         &
        \includegraphics[scale=0.05]{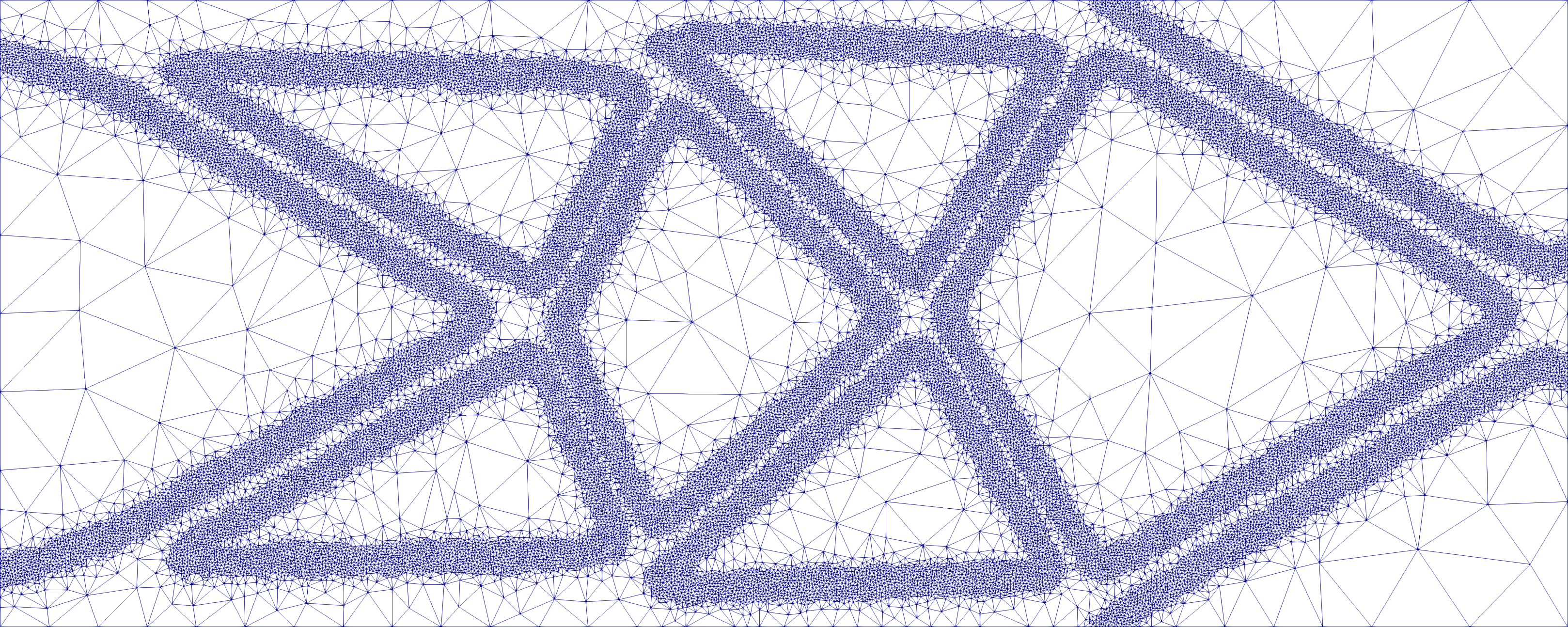}                    \\ \hline
    \end{tabular}
    \caption{Meshes for optimized designs for the AMR compliance problem with coarsening.}
    \label{fig:compliance_amr_coarsening_grid}
\end{table}

\begin{table}
    \centering
    \begin{tabular}{|M{0.5cm}|M{5.6cm}|M{5.6cm}|}
        \hline
                                         & Optimization in $\mathbb{R}^n$ & Optimization in $L^2$ \\ \hline
        \rotatebox{90}{Strategy A}     &
        \includegraphics[scale=0.05]{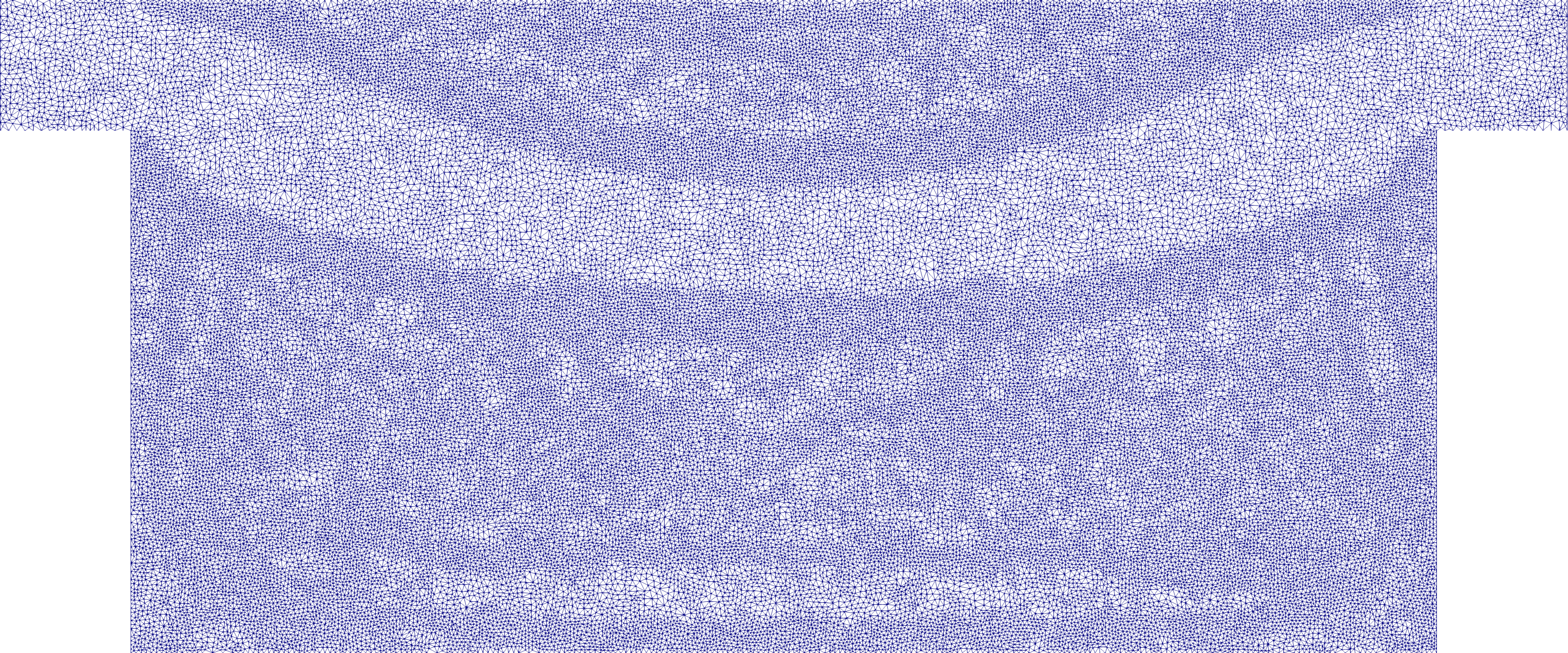}
                                         &
        \includegraphics[scale=0.05]{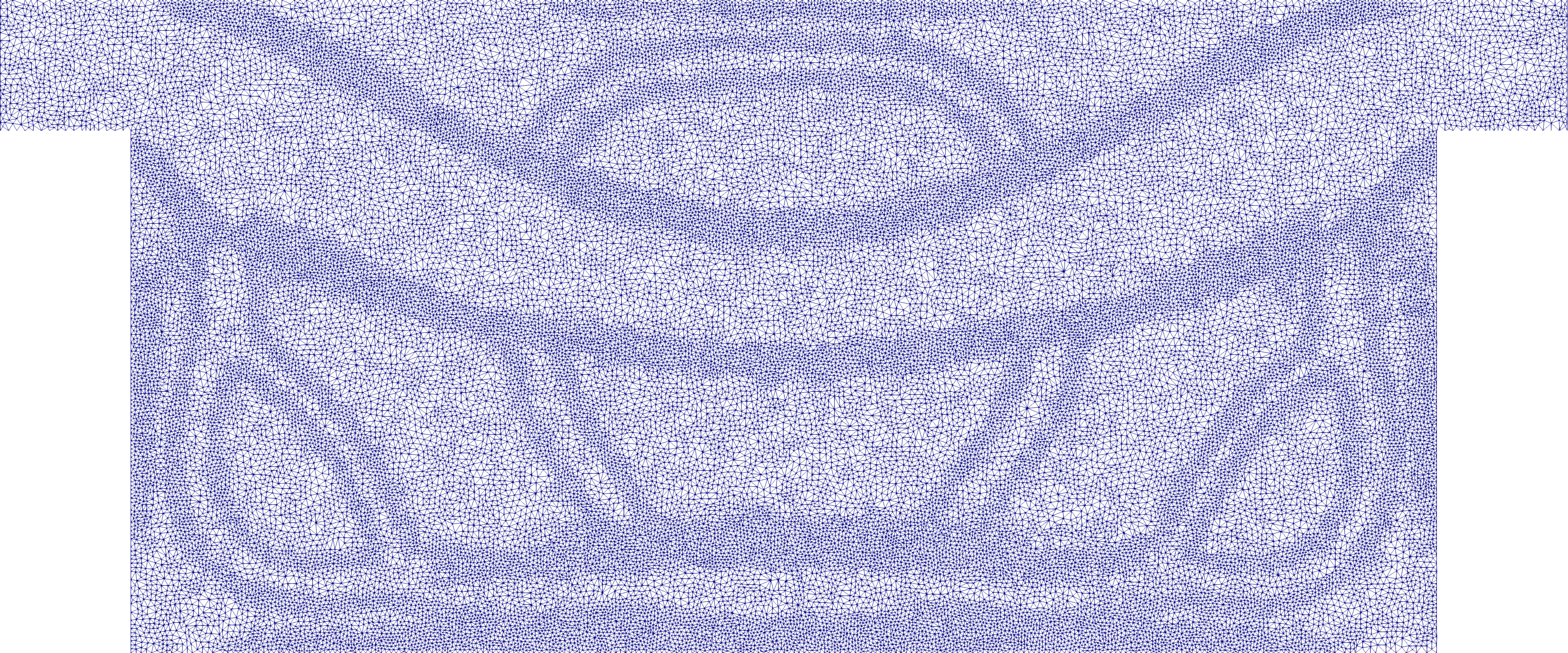}                        \\ \hline
        \rotatebox{90}{Strategy B} &
        \includegraphics[scale=0.05]{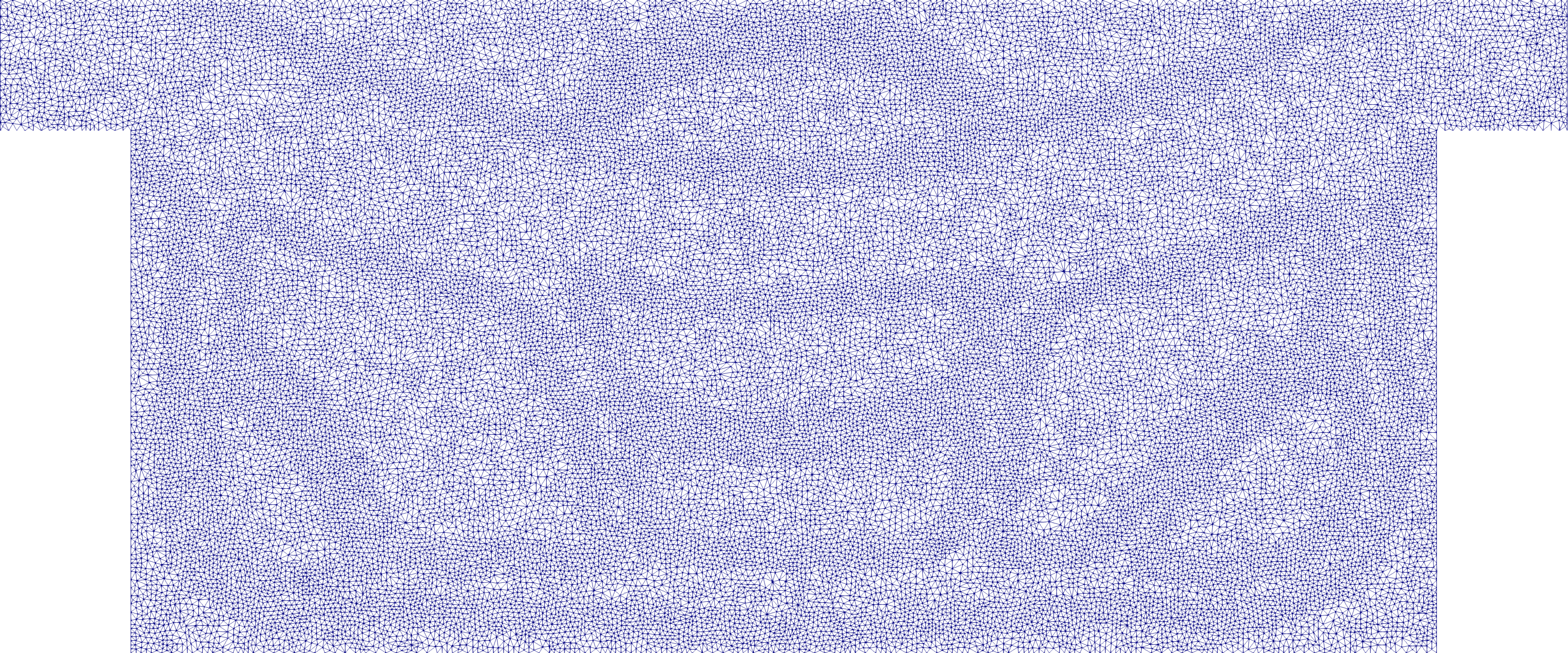}
                                         &
        \includegraphics[scale=0.05]{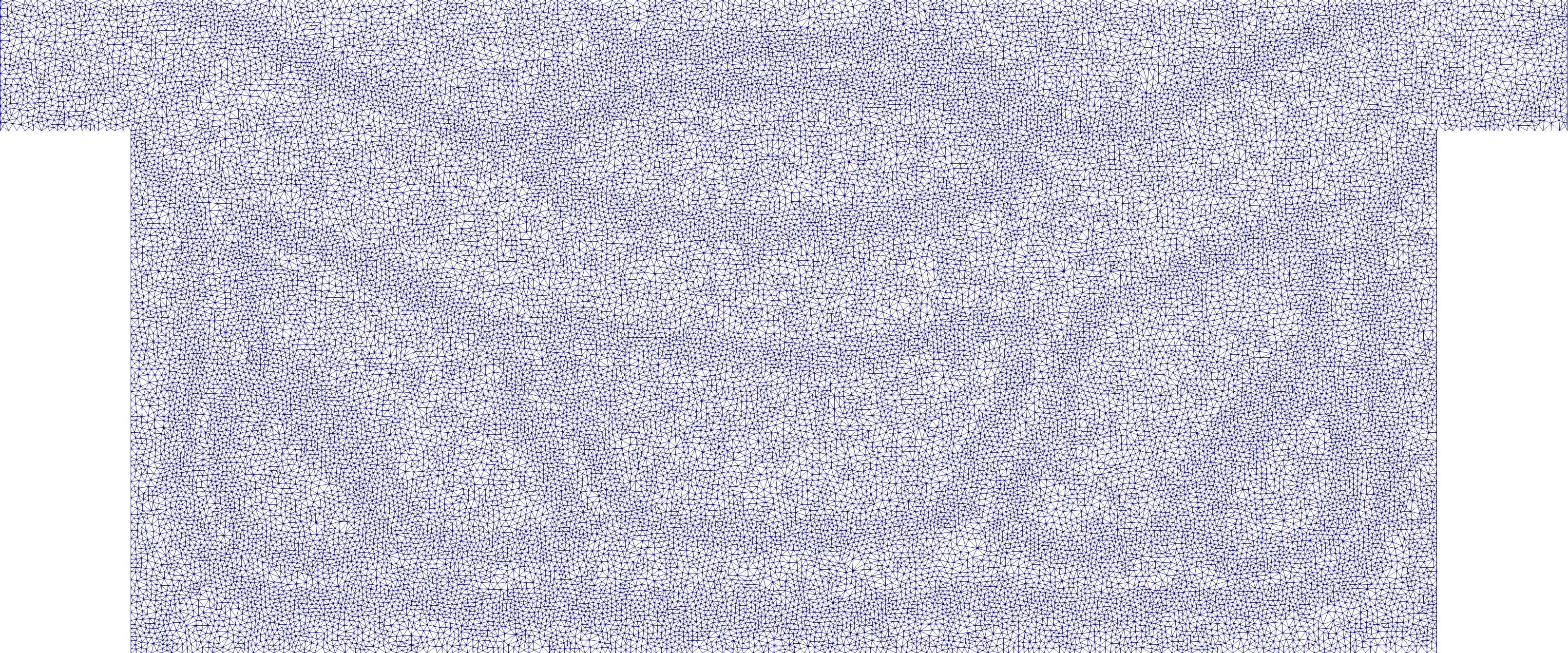}                    \\ \hline
    \end{tabular}
    \caption{Meshes for optimized designs for the AMR thermal flow problem with refinement.}
    \label{fig:thermal_flow_amr_refinement_grid}
\end{table}

\begin{table}
    \centering
    \begin{tabular}{|M{0.5cm}|M{5.6cm}|M{5.6cm}|}
        \hline
                                         & Optimization in $\mathbb{R}^n$ & Optimization in $L^2$ \\ \hline
        \rotatebox{90}{Strategy A}     &
        \includegraphics[scale=0.05]{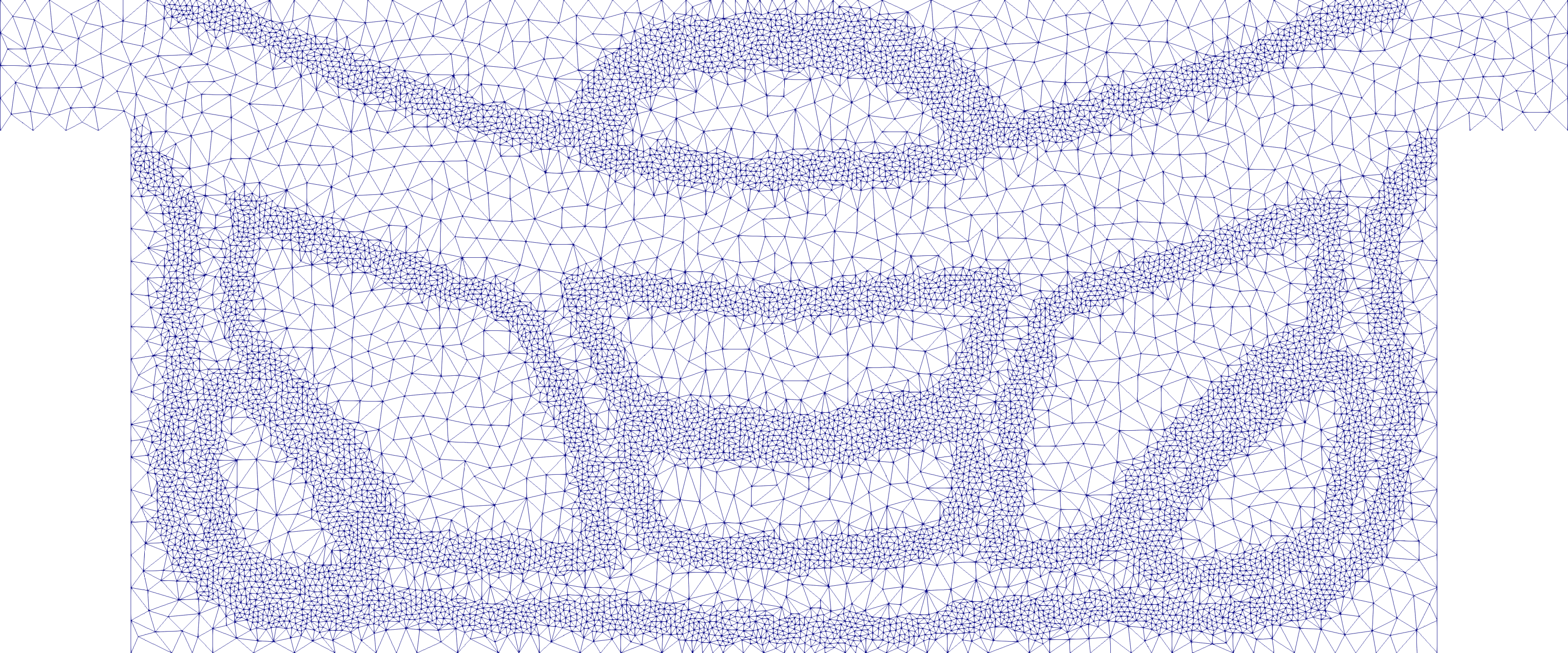}
                                         &
        \includegraphics[scale=0.05]{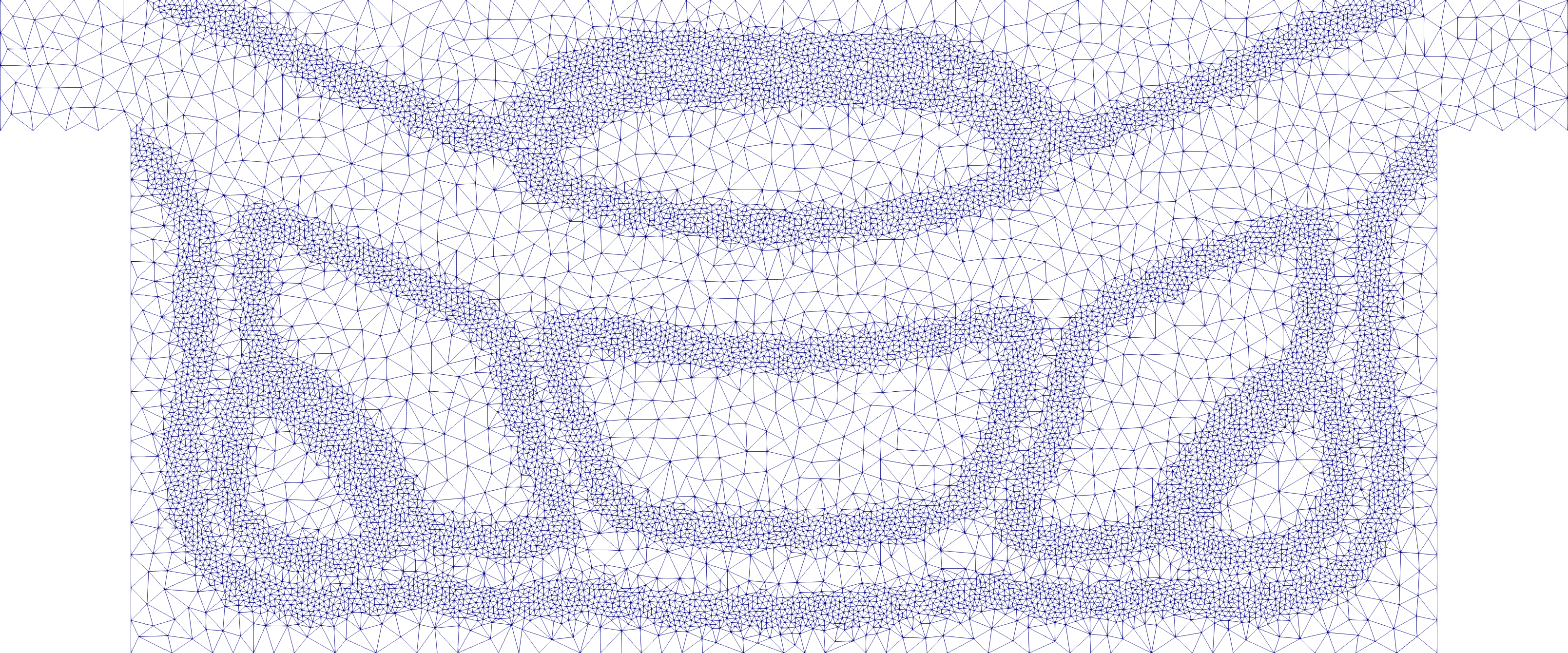}                        \\ \hline
        \rotatebox{90}{Strategy B} &
        \includegraphics[scale=0.05]{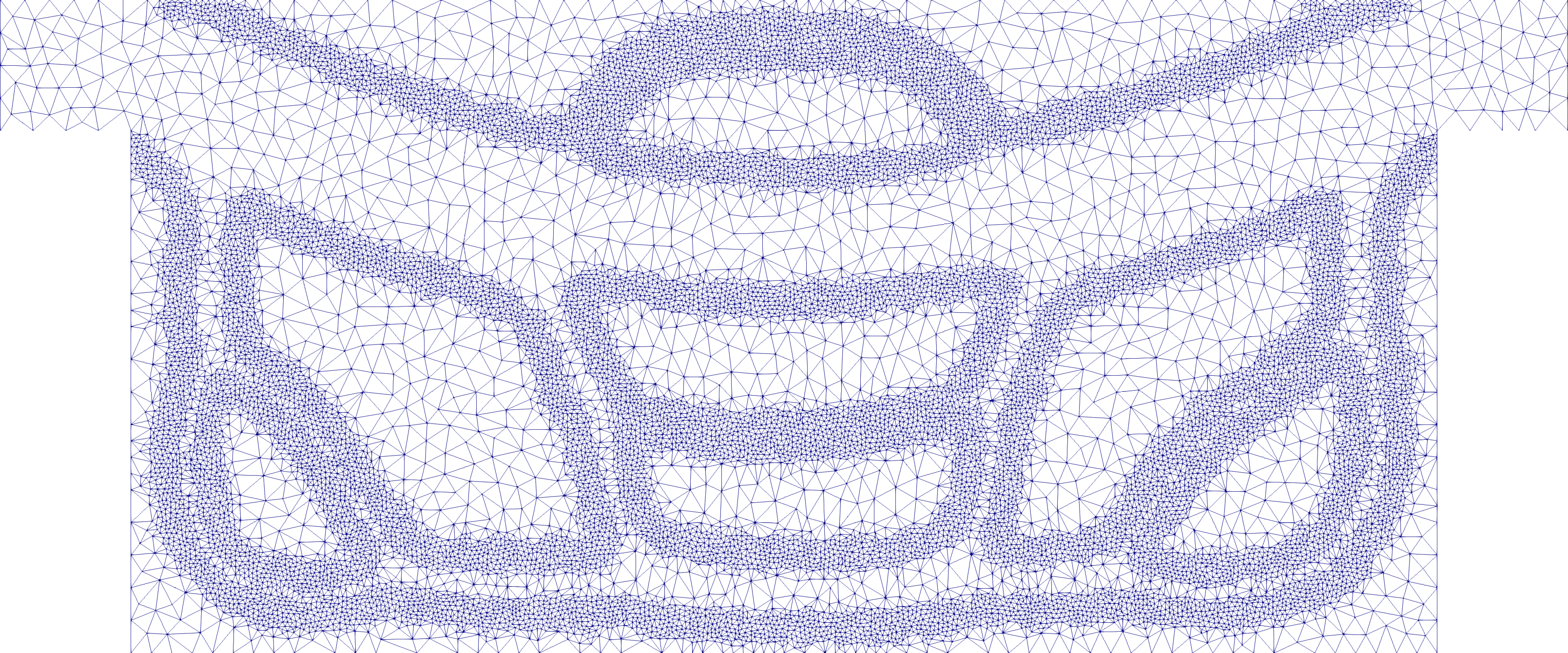}
                                         &
        \includegraphics[scale=0.05]{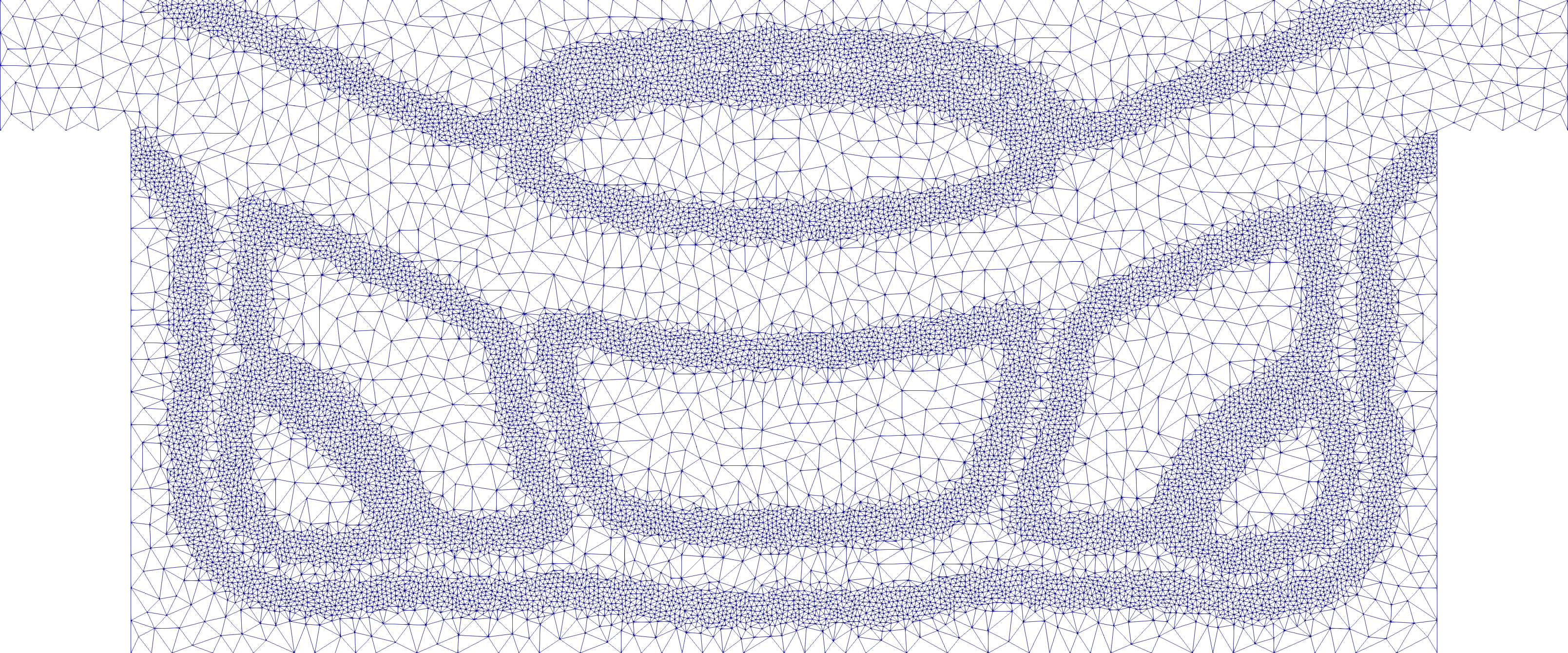}                    \\ \hline
    \end{tabular}
    \caption{Meshes for optimized designs for the AMR thermal flow problem with coarsening.}
    \label{fig:thermal_flow_amr_coarsening_grid}
\end{table}

\end{document}